\documentclass[12pt,reqno]{amsart}
\usepackage[T1]{fontenc}
\usepackage{lmodern}
\usepackage[margin=1.05in]{geometry}
\usepackage{amsmath,amssymb,amsthm,mathtools}
\usepackage{enumitem,microtype,setspace,xcolor}
\usepackage[colorlinks=true,linkcolor=blue!55!black,citecolor=blue!55!black,urlcolor=blue!55!black]{hyperref}
\hypersetup{pdftitle={On the generators of separable unital C*-algebras},pdfauthor={George A. Elliott, Chun Guang Li, Zhuang Niu}}
\numberwithin{equation}{section}
\newtheorem{theorem}{Theorem}[section]
\newtheorem{proposition}[theorem]{Proposition}
\newtheorem{lemma}[theorem]{Lemma}
\newtheorem{corollary}[theorem]{Corollary}

\theoremstyle{definition}
\newtheorem{definition}[theorem]{Definition}
\newtheorem{remark}[theorem]{Remark}
\newtheorem{convention}[theorem]{Convention}
\newcommand{\C}{\mathbb C}

\newcommand{\R}{\mathbb R}
\newcommand{\T}{\mathrm T}

\newcommand{\sa}{\mathrm{sa}}

\newcommand{\rank}{\operatorname{rank}}

\newcommand{\im}{\operatorname{Im}}
\newcommand{\re}{\operatorname{Re}}
\newcommand{\Her}{\operatorname{Her}}

\newcommand{\Tr}{\operatorname{Tr}}

\setlist[enumerate]{label=\textup{(\roman*)},leftmargin=*}
\title[On the generators of separable unital $\mathrm C^*$-algebras]{On the generators of separable unital $\mathrm C^*$-algebras}
\author{George A. Elliott}
\address{Department of Mathematics, University of Toronto, Toronto,
Ontario M5S 2E4, Canada}
\email{elliott@math.toronto.edu}
\author{Chun Guang Li}
\address{School of Mathematics and Statistics, Northeast Normal
University, Changchun 130024, P. R. China}
\email{licg864@nenu.edu.cn}
\author{Zhuang Niu}
\address{Department of Mathematics and Statistics, University of
Wyoming, Laramie, WY 82071, USA}
\email{zniu@uwyo.edu}
\date{}
\begin{document}
\begin{abstract}
Let $A$ be a separable unital $\mathrm C^*$-algebra with nonempty
tracial state space and strict comparison by ordinary tracial states.
If $\mathcal Z$ embeds unitally into $A$, then $A$ is singly generated.
In particular, every simple separable unital $\mathrm C^*$-algebra $A$
with $\mathrm T(A)\ne\varnothing$ and strict comparison is singly
generated.
\end{abstract}
\maketitle
\enlargethispage{2pt}
{\footnotesize\tableofcontents}

\section{Introduction}

In \cite{Kadison-list}, among other questions, Kadison asked whether
every von Neumann algebra acting on a separable Hilbert space is
singly generated (see \cite{Ge_2003} for a survey of Kadison's problem
list).  Many von Neumann algebras are known to be singly generated,
and the problem was reduced to the $\mathrm{II}_1$ factor case
\cite{Willig_1974}.  The recent preprint \cite{OpenAI2026} resolves
that remaining case.

The analogous question for $\mathrm C^*$-algebras asks which separable
$\mathrm C^*$-algebras are singly generated.  Commutative examples such
as $\mathrm C(S^2)$ show that separability alone is insufficient;
see \cite{Nagisa-04}.  Among separable $\mathrm C^*$-algebras, stable
algebras and unital UHF-absorbing, $\mathcal Z$-absorbing, and
approximately divisible algebras are known to be singly generated
\cite{Topping_1968,OlsenZame,Li-Sh-AD,ThielWinter}.
An explicit single generator for the reduced group algebra
$\mathrm C_r^*(\mathbb F_2)$ was constructed in
\cite{ELNF2Generators2026}.
Single generation has also been established for Villadsen's algebras
of the first type \cite{LNR-26}, while simple unital AH algebras
which are not singly generated were constructed in \cite{ELNcounter}.
Beyond $\mathrm C^*$-algebras, Hayes \cite{Hayes2026} constructed a
$\mathrm W^*$-bundle over a compact metrizable space whose fibers are $\mathrm{II}_1$ factors
with separable predual, but no finite family of continuous sections
generates every fiber as a von Neumann algebra.

Related results were recently obtained by Bouwen, Castillejos, and
Curda \cite{BouwenCastillejosCurda}.  They prove that simple separable
tracially AF algebras have generator rank at most one, and obtain
single-generation criteria using comparison and diffuse
Cantor-spectrum subalgebras.

Here we prove single generation for separable unital algebras with
nonempty trace space, stabilized strict comparison by ordinary
tracial states, and a unital embedding of the Jiang--Su algebra
$\mathcal Z$.  This embedding supplies the
matrix cones and hereditary coefficient algebras used in the proof.
There is no nuclearity assumption or restriction on the ideal
structure or the extreme trace boundary.

We specify the comparison convention, since different formulations
need not agree for nonsimple algebras.

\begin{definition}\label{def:comparison}
Let $A$ be unital with $\T(A)\ne\varnothing$.  We say that $A$ has
\emph{strict comparison by tracial states} if, for
$a,b\in(A\otimes\mathcal K)_+$,
\[
 \mathrm d_\tau(a)<\mathrm d_\tau(b)\quad(\tau\in\T(A))
 \quad\Longrightarrow\quad a\precsim b.
\]
All traces here are ordinary traces, extended using the
unnormalized matrix trace, and
$\mathrm d_\tau(a)=\lim_m(\tau\otimes\Tr)(a^{1/m})$.
The relation $a\precsim b$ is Cuntz subequivalence.
We distinguish this condition from comparison by the full cone of
lower semicontinuous traces in Section~\ref{sec:structural}.
\end{definition}

For $X=\T(A)$ put
\[
 \|a\|_{2,X}=\sup_{\tau\in X}\tau(a^*a)^{1/2}.
\]
We write $\mathcal M_A$ for the bounded uniform tracial completion.
Its construction first identifies elements at zero $2$-seminorm and
then completes norm-bounded balls.  The canonical map
$A\to\mathcal M_A$ is trace preserving and may have a kernel.
Generation throughout the paper is norm generation in the original
$\mathrm C^*$-algebra, without automatically adjoining a unit.

\begin{theorem}\label{thm:main}
Let $A$ be a separable unital $\mathrm C^*$-algebra with
$\T(A)\ne\varnothing$ and strict comparison by tracial states.
If $\mathcal Z$ embeds unitally into $A$, then $A$ is singly generated.
\end{theorem}

The embedding hypothesis can be verified from intrinsic
Cuntz-semigroup regularity.  Recall that $A$ is \emph{pure} if
$\mathrm{Cu}(A)$ is almost unperforated and almost divisible.
A unital pure algebra of stable rank one admits a unital
$\mathcal Z$-embedding by \cite[Proposition~7.6.5]{APT2018}.
When its lower semicontinuous $2$-quasitraces are traces, comparison
by tracial states follows from comparison on the full trace cone;
the quotient argument establishing this implication is included below.

\begin{corollary}[Corollary~\ref{cor:intrinsic}]
\label{cor:intro-intrinsic}
Let $A$ be a separable unital pure $\mathrm C^*$-algebra of stable
rank one.  Suppose that every lower semicontinuous $2$-quasitrace
on $A$ is an ordinary trace.  Then $A$ is singly generated.
\end{corollary}

In particular, the conclusion holds for separable unital nowhere
scattered algebras of stable rank one with strict comparison by the
full ordinary trace cone; see Corollary~\ref{cor:nowhere-scattered}.
Here nowhere scattered means that no nonzero ideal-quotient is
elementary.  Pureness in this case follows from
\cite[Theorem~5.12]{APTV2024}.

The simple tracial comparison theorem is a further consequence.

\begin{corollary}[Corollary~\ref{cor:simple-comparison}]
\label{cor:intro-simple}
Let $A$ be a simple separable unital $\mathrm C^*$-algebra with
$\T(A)\ne\varnothing$ and strict comparison by tracial states.
Then $A$ is singly generated.
\end{corollary}

In the infinite-dimensional simple case, the embedding hypothesis
of Theorem~\ref{thm:main} follows from comparison.  Indeed,
Lemma~\ref{lem:simple-state-comparison-quasitraces} passes to
comparison by the full ordinary trace cone.  Thiel's theorem
\cite[Corollary~B]{Thiel2026} then gives pureness, and Lin's theorem
\cite[Corollary~4.5]{Lin2024} yields stable rank one.
Proposition~\ref{prop:structural-z} supplies the required unital
$\mathcal Z$-embedding.
The finite-dimensional case is elementary.

For nonsimple algebras, the unital $\mathcal Z$-embedding in
Theorem~\ref{thm:main} supplies noncommutative matrix structure
uniformly over the trace space.  Strict comparison alone does not
imply single generation in this generality.  For example, if $X$
is the geometric realization of the nonplanar finite graph $K_5$,
then
$\mathrm{Cu}(\mathrm C(X))\cong
 \mathrm{Lsc}(X,\overline{\mathbb N})$
by \cite[Theorem~3.4]{APS2011}.
Thus $\mathrm C(X)$ has strict comparison by tracial states, but it
is not singly generated: a generator would give a continuous
embedding $X\to\mathbb C$.

The proof begins with a completely positive contractive order-zero
map $\phi:\mathrm M_n\to\mathcal Z\subseteq A$ such that
$\tau(\phi(1_n))$ is arbitrarily close to one, uniformly over
$\mathrm T(A)$.  A cutdown of one diagonal entry
defines a full hereditary sub-$\mathrm C^*$-algebra $D$.
The defining positive element has constant positive rank across
the trace space, so the normalized trace space of $D$ is affinely
homeomorphic to $\mathrm T(A)$.  Its strictly positive Cuntz class
is almost divisible, and $D$ inherits strict comparison.
Proposition~1.4 and Theorem~1.3 of \cite{EHSS2026} therefore give
real rank zero of the uniform tracial completion of $D$, together
with realization and comparison of projections by their continuous
affine trace profiles.

These regularity properties allow us to construct a unital simple
AF sub-$\mathrm C^*$-algebra in the completion, with the same
trace simplex as $D$ and a uniform tracial averaging property.
Its tracial closures provide the auxiliary inclusions for our
adaptation of the relative-generation method of \cite{OpenAI2026}.
For an irreducible inclusion of finite factors, that method
produces a unitary which generates the ambient factor together
with the subfactor.  Here the estimates must hold uniformly over
the entire trace simplex and ultimately yield norm generation
in $D$.

Comparison first captures the AF generators in the uniform
$2$-closure of a finitely generated sub-$\mathrm C^*$-algebra of
$D$.  Finite-word estimates and spectral clipping then give
unitary perturbations small in operator norm.  A Baire argument
that allows the trace to vary gives relative generation
simultaneously over the trace space.  A further comparison
argument stores the remaining errors in actual elements of $D$,
upgrading tracial density to norm generation and giving seven
selfadjoint generators for $D$.

The matrix cone provided by $\phi$ then compresses these seven generators
to two selfadjoint generators of $A$.  To recover the hereditary
bimodules in this final step, each connecting element is chosen
so that its initial square is full in the initial hereditary
algebra.  This fullness condition makes the compression work
without simplicity.

Real rank zero gives a density theorem for the intrinsic class.
For a $\mathrm C^*$-algebra $A$, write
\[
 \operatorname{Gen}(A)=\{g\in A:\mathrm C^*(g)=A\}.
\]
We use Thiel's generator rank $\operatorname{gr}$ \cite{Thiel2021}.

\begin{theorem}[Theorem~\ref{thm:rr0-intrinsic-density}]
\label{thm:intro-rr0-intrinsic-density}
Let $A$ be a separable unital pure $\mathrm C^*$-algebra of real
rank zero and stable rank one.  Suppose that every lower
semicontinuous $2$-quasitrace on $A$ is an ordinary trace.
Then $\operatorname{Gen}(A)$ is a norm-dense $G_\delta$ subset of
$A$.  Equivalently, $\operatorname{gr}(A)\leq1$.
\end{theorem}

The density argument uses finite spectral refinements and single
generation of all nonzero projection corners.  In a separable
real-rank-zero algebra, elements that generate an entire rectangular
corner as a bimodule are dense in that corner.  These connecting
elements recover every off-diagonal block and allow the argument
to accommodate arbitrary ideals.
The simple tracial comparison result follows as a corollary.

\begin{corollary}[Corollary~\ref{cor:rr0-simple-density}]
\label{cor:rr0-density}
Let $A$ be a simple separable unital $\mathrm C^*$-algebra with
$\T(A)\ne\varnothing$, real rank zero, and strict comparison by
tracial states.  Then $\operatorname{Gen}(A)$ is a norm-dense
$G_\delta$ subset of $A$.
\end{corollary}

Thiel proved the nonsimple density conclusion for separable
$\mathcal Z$-stable algebras of real rank zero
\cite[Theorem~5.3]{Thiel2022Generators}.
Our intrinsic theorem uses hereditary Cuntz-semigroup regularity
and the single-generation theorem above.

The matrix-cone compression also controls topological free entropy
dimension.  Its scalar spectral labels and sparse connecting blocks
give an operator-norm microstate covering estimate.  Increasing the
matrix size while keeping the seven-generator coefficient bound
fixed makes the upper bound approach one.
We use
\[
 \delta_{\mathrm{top}}(g)
 =\delta_{\mathrm{top}}(\re g,\im g),\qquad
 \chi_{\mathrm{top}}(g)
 =\chi_{\mathrm{top}}(\re g,\im g)
\]
as in \cite{VoiculescuTopological}.  Recall that a separable unital
$\mathrm C^*$-algebra is $\mathrm{MF}$ if it embeds unitally into a
quotient
$\prod_j\mathrm M_{q_j}/\bigoplus_j\mathrm M_{q_j}$.
Such an algebra has a tracial state.

\begin{theorem}[Theorem~\ref{thm:entropy-infimum}]
\label{thm:intro-entropy}
Let $A$ be a separable unital $\mathrm{MF}$ $\mathrm C^*$-algebra
with strict comparison by tracial states and a unital embedding of
$\mathcal Z$.  Then
\[
 \inf_{\substack{g\in A\\\mathrm C^*(g)=A}}
 \delta_{\mathrm{top}}(g)=1.
\]
For every $\eta>0$, there is an invertible generator $g=uh$, with
$u\in\mathcal U(A)$ and $h=|g|\ge1_A$, such that
\[
 1\le\delta_{\mathrm{top}}(g)<1+\eta,\qquad
 \chi_{\mathrm{top}}(g)=-\infty.
\]
\end{theorem}

The infimum statement does not assert that a generator of dimension
exactly one exists.

Note that, for $n\in\{2,3,\ldots,\infty\}$, the reduced free
group $\mathrm C^*$-algebras $\mathrm C_r^*(\mathbb F_n)$ are simple
\cite{Powers1975} and $\mathrm{MF}$
\cite[Corollary~8.4]{HaagerupThorbjornsen2005}.
For finite $n$, strict comparison follows from
\cite[Theorem~A]{AGKEP2025}, using Robert's selflessness
\cite{RobertSelfless2025} and hence the work of Dykema and
R\o rdam \cite{DykemaRordam2000}; for $n=\infty$, see
\cite[Proposition~6.3.2]{Robert2012}.
Thus Corollary~\ref{cor:intro-simple} and
Theorem~\ref{thm:intro-entropy} apply.

Applying Corollary~\ref{cor:intro-simple} to a suitable weakly dense
separable sub-$\mathrm C^*$-algebra of a $\mathrm{II}_1$ factor also
recovers the finite-factor conclusion of \cite{OpenAI2026}; see
Remark~\ref{rem:factor-connection}.

Single generation is substantially weaker than strict comparison
even among simple separable unital algebras with nonempty trace
space.  The first-type Villadsen algebras in \cite{LNR-26}
include examples without strict comparison
\cite[Section~8]{TomsWinterVilladsen2009}.
Our transformation-group theorem \cite{ELNTransformation2026}
also applies beyond strict comparison, including to the crossed
products of Giol and Kerr \cite[Theorem~2.2]{GiolKerr2010}.

Section~\ref{sec:structural} records comparison, $\mathcal Z$-embeddings,
and the regularity of the coefficient tracial completions.
Section~\ref{sec:external-af-core} constructs the AF algebra, and
Section~\ref{sec:quantization} establishes uniform finite-word
estimates.  Section~\ref{sec:audited-nonunital-capture} proves
comparison storage and finite capture.
Section~\ref{sec:uniform-relative-generation} proves relative
generation by a finite perturbation and a Baire argument; the
perturbation is adapted from \cite{OpenAI2026}.
Section~\ref{sec:seven-generators} obtains a uniform finite generator
bound, and Section~\ref{sec:cone-compression} gives the matrix-cone
compression.  Section~\ref{sec:main-proof} proves
Theorem~\ref{thm:main} and its intrinsic and simple corollaries.
Section~\ref{sec:rr0-density} proves
Theorem~\ref{thm:intro-rr0-intrinsic-density} and
Corollary~\ref{cor:rr0-density}.  Section~\ref{sec:entropy} proves
the entropy estimates.

\medskip
\noindent\textbf{Acknowledgment.}
The work was carried out using ChatGPT.  The first author is
supported by an NSERC Discovery Grant, the second author is
supported by an NNSF Grant of China (No.~12671150), and the third
author is supported by a Simons Foundation Grant
(MP-TSM-00002606).

\section{Comparison, Jiang--Su embeddings, and tracial completions}
\label{sec:structural}\label{sec:structural-inputs}

We extend tracial states to matrix algebras using the unnormalized
matrix trace.  Unless explicitly stated otherwise, strict comparison
means the stabilized comparison by ordinary tracial states in
Definition~\ref{def:comparison}.  The following distinction is relevant
when simplicity is not assumed.

\begin{convention}\label{def:full-trace-comparison}
Write $\mathrm T_{\mathrm{lsc}}(A)$ for the cone of all lower
semicontinuous $[0,\infty]$-valued traces on $A_+$, and
$\mathrm{QT}_{\mathrm{lsc}}(A)$ for the corresponding cone of
$2$-quasitraces.  Densely finite weights are not the only weights
included in these cones.  A $2$-quasitrace is a homogeneous positive
function that is additive on commuting positive elements, satisfies
$q(x^*x)=q(xx^*)$, and has a compatible such extension to
$\mathrm M_2(A)$; additivity on all positive elements makes it a trace.
All these weights are extended to $A\otimes\mathcal K$.

We say that $A$ has \emph{strict comparison by the full ordinary trace
cone} if, for $a,b\in(A\otimes\mathcal K)_+$, the existence of
$0<\gamma<1$ with
\[
 \mathrm d_\rho(a)\leq\gamma\,\mathrm d_\rho(b)
                  \qquad(\rho\in\mathrm T_{\mathrm{lsc}}(A))
\]
implies $a\precsim b$.  This is the full-cone convention of
\cite[Definition~3.1]{NgRobert2016}.  It differs from the pointwise
strict comparison by normalized tracial states used in our main
theorem.  We keep the two conventions separate.
\end{convention}

Recall that a Cu-semigroup $S$ is \emph{almost unperforated} if
$(k+1)x\leq ky$ implies $x\leq y$, for $x,y\in S$ and $k\geq1$.
It is \emph{almost divisible} if, for $x'\ll x$ and $n\geq1$, there is
$z\in S$ such that
\[
 nz\leq x,\qquad x'\leq(n+1)z.
\]
An element $x$ is called almost divisible if this condition holds for
that element.  A $\mathrm C^*$-algebra is \emph{pure} if its Cuntz
semigroup is almost unperforated and almost divisible.  Almost
unperforation is equivalent to strict comparison by the full cone of
lower semicontinuous $2$-quasitraces, with the uniform-ratio convention
above; see \cite[Section~3, before Definition~3.1]{NgRobert2016}.
A unital algebra has \emph{stable rank one} if its invertible elements
are norm dense.

\begin{lemma}\label{lem:simple-state-comparison-quasitraces}
Let $A$ be a simple unital $\mathrm C^*$-algebra with
$\mathrm T(A)\ne\varnothing$ and strict comparison by ordinary
tracial states.  Then $A$ has strict comparison by the full
ordinary trace cone, and every lower semicontinuous $2$-quasitrace
on $A$ is an ordinary trace.
\end{lemma}

\begin{proof}
Suppose that $a,b\in(A\otimes\mathcal K)_+$ and $0<\gamma<1$
satisfy
\[
 \mathrm d_\rho(a)\le\gamma\mathrm d_\rho(b)
            \qquad(\rho\in\mathrm T_{\mathrm{lsc}}(A)).
\]
If $b=0$, faithfulness of a tracial state gives $a=0$.
Otherwise every $\tau\in\mathrm T(A)$ is faithful, so
$\mathrm d_\tau(b)>0$.
Fix $\varepsilon>0$.  Finite-matrix localization gives
$c\in\mathrm M_m(A)_+$ with $[(a-\varepsilon)_+]\le[c]$.
Thus $\mathrm d_\tau((a-\varepsilon)_+)\le m<\infty$.
If $\mathrm d_\tau(b)<\infty$, then
\[
 \mathrm d_\tau((a-\varepsilon)_+)
 \le\mathrm d_\tau(a)
 \le\gamma\mathrm d_\tau(b)<\mathrm d_\tau(b).
\]
If $\mathrm d_\tau(b)=\infty$, finiteness of the left side gives
the same strict inequality.
Comparison by tracial states gives $(a-\varepsilon)_+\precsim b$.
Letting $\varepsilon\downarrow0$ proves $a\precsim b$, as required
by Convention~\ref{def:full-trace-comparison}.
Ng--Robert's theorem \cite[Theorem~3.6\textup{(i)}]{NgRobert2016}
now gives the assertion about quasitraces.
\end{proof}

\begin{lemma}\label{lem:intrinsic-state-comparison}
Let $A$ be a nonzero unital pure $\mathrm C^*$-algebra of stable rank one.
Suppose that every lower semicontinuous $2$-quasitrace on $A$ is a
trace.  Then $\mathrm T(A)\ne\varnothing$, and $A$ has strict
comparison by ordinary tracial states in the sense of
Definition~\ref{def:comparison}.
\end{lemma}

\begin{proof}
Every nonzero unital quotient of $A$ has stable rank one and is
therefore stably finite.  By \cite[Theorem~2.4]{Handelman1981}, it
admits a unital homomorphism into a finite $\mathrm{AW}^*$-factor.
Pulling back the canonical normalized quasitrace gives a normalized
lower semicontinuous $2$-quasitrace by
\cite[Proposition~II.1.6 and Corollary~II.1.10]{BlackadarHandelman1982}.
Its pullback to $A$ is a trace by
hypothesis, so the quotient has a tracial state.  Taking the quotient
to be $A$ first shows that $\mathrm T(A)$ is nonempty.

Suppose that $a,b\in(A\otimes\mathcal K)_+$ satisfy
\[
 \mathrm d_\tau(a)<\mathrm d_\tau(b)
                         \qquad(\tau\in\mathrm T(A)).
\]
The element $b$ must be full.  Otherwise its ideal is
$I\otimes\mathcal K$ for a proper ideal $I$ of $A$, and a tracial
state on the nonzero unital quotient $A/I$ pulls back to a tracial
state of $A$ annihilating $b$.  This contradicts the displayed strict
inequality.

Fix $\varepsilon>0$.  Finite-matrix localization gives
$c\in\mathrm M_m(A)_+$ and $\eta>0$ such that
\[
 [(a-\varepsilon)_+]\leq[(c-\eta)_+],\qquad [c]\leq[a].
\]
For completeness, put $p_m=1_A\otimes\sum_{j=1}^m e_{jj}$ and use
$a^{1/2}p_ma^{1/2}\to a$ in norm.  Cutdown comparison and the
Cuntz equivalence of the corresponding cutdowns of $vv^*$ and
$v^*v$, with $v=a^{1/2}p_m$, give these assertions with
$c=p_map_m$ for sufficiently large $m$.
Choose a continuous function $f$ on $[0,\|c\|]$ with
$0\leq f\leq1$, $f(0)=0$, and $f=1$ on $[\eta,\|c\|]$.
The finite continuous affine function
$F(\tau)=\tau(f(c))$ satisfies
\[
 \mathrm d_\tau((a-\varepsilon)_+)
 \leq F(\tau)\leq\mathrm d_\tau(c)
 \leq\mathrm d_\tau(a)<\mathrm d_\tau(b).
\]
Lower semicontinuity and compactness of $\mathrm T(A)$ yield a
number $\delta>0$ with
$F(\tau)+\delta\leq\mathrm d_\tau(b)$ for every $\tau$.
If $M=\max F>0$, take $\gamma=M/(M+\delta)$; if $M=0$, take
$\gamma=1/2$.  Then $0<\gamma<1$ and
\begin{equation}\label{eq:intrinsic-uniform-ratio}
 \mathrm d_\tau((a-\varepsilon)_+)
 \leq\gamma\,\mathrm d_\tau(b)
                      \qquad(\tau\in\mathrm T(A)).
\end{equation}

Every lower semicontinuous trace $\rho$ with
$0<\rho(1_A)<\infty$ is a scalar multiple of a tracial state, so
\eqref{eq:intrinsic-uniform-ratio} holds for it as well.  A weight
with $\rho(1_A)=0$ is zero.  If $\rho(1_A)=\infty$, fullness of
$b$ and compactness of $[1_A]$ give $[1_A]\leq k[b]$ for some
$k$, hence $\mathrm d_\rho(b)=\infty$.  The same inequality therefore
holds for every $\rho\in\mathrm T_{\mathrm{lsc}}(A)$.
Pureness, the hypothesis on quasitraces, and full-cone comparison
now give $(a-\varepsilon)_+\precsim b$.
Letting $\varepsilon\downarrow0$ proves $a\precsim b$.
\end{proof}

\begin{proposition}\label{prop:structural-z}\label{prop:audited-z-embedding}\label{prop:z-embedding}
Let $A$ be a nonzero unital $\mathrm C^*$-algebra of stable rank one.
If $\mathrm{Cu}(A)$ is almost unperforated and $[1_A]$ is almost
divisible, then $\mathcal Z$ embeds unitally in $A$.
In particular this holds for every nonzero unital pure algebra of
stable rank one.
\end{proposition}

\begin{proof}
By \cite[Paragraph~7.6.4]{APT2018}, almost unperforation is the
strict comparison hypothesis used in
\cite[Proposition~7.6.5]{APT2018}.  The implication
\textup{(2)}$\Rightarrow$\textup{(4)} of that proposition gives a
unital homomorphism $\mathcal Z\to A$.  It is injective because
$\mathcal Z$ is simple.  The criterion does not require simplicity
of $A$ or identification of its Cuntz semigroup with a space of
affine functions.
\end{proof}

\begin{proposition}\label{prop:z-cone}
Let $A$ be a unital $\mathrm C^*$-algebra with
$\mathrm T(A)\ne\varnothing$, containing a fixed unital copy of
$\mathcal Z$.  For every $n\geq1$ and $\eta>0$ there is a c.p.c.\
order zero map $\phi:\mathrm M_n\to A$ with range in this copy and
\[
 \tau(\phi(1_n))>1-\eta\qquad(\tau\in\mathrm T(A)).
\]
Write $t=\phi(1_n)$.  If
$0<\delta<\tau_{\mathcal Z}(t)$, then
$d=(\phi(e_{11})-\delta)_+$ is full in $A$ and
\[
 \mathrm d_\tau(d)=c_\delta>0\quad(\tau\in\mathrm T(A)),
 \qquad nc_\delta\geq\tau_{\mathcal Z}(t)-\delta.
\]
The class $[d]$ is almost divisible in
$\mathrm{Cu}(\operatorname{Her}_A(d))$.
\end{proposition}

\begin{proof}
There is a c.p.c.\ order zero map
$\psi:\mathrm M_n\to\mathcal Z$ with
$\tau_{\mathcal Z}(\psi(1_n))>1-\eta$; see
\cite[opening paragraph of Section~6]{SatoWhiteWinter2015}.
Compose with the fixed embedding.  Every tracial state of $A$
restricts to the unique tracial state on this copy of $\mathcal Z$,
which gives the first assertion.
The supporting-homomorphism description of order zero maps
\cite[Theorem~2.3 and Corollary~3.2]{WinterZacharias2009} gives
\[
 \mathrm d_\tau((\phi(e_{11})-\delta)_+)
   =\frac1n\mathrm d_\tau((t-\delta)_+).
\]
The ranks are independent of $\tau$, since both cutdowns belong to
$\mathcal Z$.  If $\mu$ is the spectral distribution of $t$ for its
trace, then
\[
 \tau_{\mathcal Z}(t)
 =\int_0^1 s\,\mathrm d\mu(s)
 \leq\delta+\mu((\delta,1])=\delta+nc_\delta.
\]
Thus $d\ne0$.  Simplicity of $\mathcal Z$ makes $d$ full in
$\mathcal Z$, hence full in $A$ because the embedding is unital.

The description
$\mathrm{Cu}(\mathcal Z)=\mathbb N_0\sqcup(0,\infty]$
\cite[Paragraph~7.3.2]{APT2018}, with compact
integers distinguished from soft integers, shows that for every
$x\in\mathrm{Cu}(\mathcal Z)$ and $k\geq1$ there is a class $y$
with $ky\leq x\leq(k+1)y$.
For a soft $x=r$, use the soft class $r/k$ (and $\infty$ when
$r=\infty$).  For a nonzero compact integer $r$, choose a soft
number $s$ with $r/(k+1)<s<r/k$.  The mixed order rules give the
same inequalities; the case $x=0$ is immediate.
Apply this to $x=[d]$ and pass the inequalities to
$\mathrm{Cu}(A)$.  The full hereditary inclusion
$\operatorname{Her}_A(d)\subseteq A$ induces an isomorphism on
Cuntz semigroups, so the asserted almost divisibility follows.
\end{proof}

We next recall the form of regularity of a uniform tracial
completion that we use.  A \emph{tracially complete $\mathrm C^*$-algebra} is
a pair $(\mathcal M,X)$ in which $\mathcal M$ is unital,
$X\subseteq \mathrm T(\mathcal M)$ is nonempty, compact, and convex, the
uniform $2$-seminorm
\[
 \|x\|_{2,X}=\sup_{\tau\in X}\tau(x^*x)^{1/2}
\]
is a norm, and the norm-closed unit ball of $\mathcal M$ is
complete for this norm.  The pair is \emph{factorial} if $X$ is a
face of $\mathrm T(\mathcal M)$: whenever
$t\rho+(1-t)\sigma\in X$, with $0<t<1$ and
$\rho,\sigma\in \mathrm T(\mathcal M)$, both $\rho$ and $\sigma$ belong
to $X$.  Equivalently, $\pi_\tau(\mathcal M)''$ is a factor for
every extreme point $\tau$ of $X$.  The pair is of type
$\mathrm{II}_1$ if $\pi_\tau(\mathcal M)''$ is a type
$\mathrm{II}_1$ von Neumann algebra for every $\tau\in X$.

For a factorial pair of type $\mathrm{II}_1$, we use the term
\emph{regular} in the sense of \cite[Theorem~1.3]{EHSS2026}.
The equivalent formulation needed here consists of real rank zero
of $\mathcal M$ and the following two projection properties:
\begin{enumerate}
\item for projections $p,q\in\mathcal M$, one has
$p\precsim_{\mathrm{MvN}}q$ if and only if
$\tau(p)\leq\tau(q)$ for all $\tau\in X$;
\item for every continuous affine $f:X\to[0,1]$, there is a
projection $p\in\mathcal M$ with $\tau(p)=f(\tau)$ for every
$\tau\in X$.
\end{enumerate}
Here $p\precsim_{\mathrm{MvN}}q$ means that some partial isometry
$v\in\mathcal M$ satisfies $v^*v=p$ and $vv^*\leq q$.

If $B$ is a $\mathrm C^*$-algebra with compact nonempty normalized trace space,
its \emph{uniform tracial completion} $\overline B^{\,\mathrm T(B)}$ is
the quotient of the algebra of norm-bounded,
$\|\cdot\|_{2,\mathrm T(B)}$-Cauchy sequences by the ideal of
$\|\cdot\|_{2,\mathrm T(B)}$-null sequences.  Every $\tau\in \mathrm T(B)$
extends to this completion by taking ordinary limits along such
sequences.  We denote the resulting compact set of traces by
$X_B$.  The canonical map from $B$ to its completion need not be
injective: its kernel is the common trace kernel.  Tracial
calculations with elements of $B$ use their canonical images; this
identification does not assert that their original norms are
preserved.

\begin{proposition}\label{prop:audited-regular-completion}\label{prop:regular-completion}
Let $B$ be a separable $\mathrm C^*$-algebra such that
$\mathrm T(B)$ is nonempty and compact.  Suppose that $B$ has strict
comparison by ordinary tracial states, the class of a strictly
positive element is almost divisible, and $B$ has no nonzero
finite-dimensional representations.  Then
\[
 (\mathcal M_B,X_B):=(\overline B^{\,\mathrm T(B)},X_B)
\]
is a regular factorial tracially complete $\mathrm C^*$-algebra of
type $\mathrm{II}_1$.  Moreover, $\mathrm T(\mathcal M_B)=X_B$.
\end{proposition}

\begin{proof}
The completion with its designated traces is factorial and
tracially complete by \cite[Section~2.4.2]{EHSS2026}.
Bounded uniform $2$-density identifies the corresponding tracial
GNS von Neumann algebras before and after completion.
The finite GNS von Neumann algebras of $B$ have no type I summand:
a nonzero type I summand, disintegrated into finite-dimensional
factor representations, would give a nonzero finite-dimensional
representation of the separable algebra $B$.
Thus the pair is of type $\mathrm{II}_1$.
Proposition~1.4 of \cite{EHSS2026} applies to the stated comparison
and almost-divisibility hypotheses and gives regularity.
The equality of the full trace space with the designated traces
follows from \cite[Theorem~1.3]{EHSS2026}.
Neither simplicity nor nuclearity is required here.
\end{proof}

\begin{lemma}[Constant-rank hereditary trace spaces]
\label{lem:audited-constant-rank-corner}
Let $A$ be unital with $\mathrm T(A)\ne\varnothing$, and let
$d\in A_+$ be full and satisfy $\mathrm d_\tau(d)=c$ for every $\tau\in \mathrm T(A)$,
where $c>0$ is constant.  Set $D=\operatorname{Her}_A(d)$.
Then
\[
 \mathrm T(A)\longrightarrow \mathrm T(D),\qquad
                   \tau\longmapsto c^{-1}\tau|_D
\]
is an affine homeomorphism.  In particular $\mathrm T(D)$ is compact.  If
$\mathcal M_A$ is the bounded uniform tracial completion of $A$, the
support of the image of $d$ is a projection $p\in\mathcal M_A$, and the normalized
completion of $D$ identifies with $p\mathcal M_Ap$.
\end{lemma}

\begin{proof}
Fullness of $d$ gives a finite family $\xi_i=x_i d^{1/2}\in A$ such
that $\sum_i\xi_i\xi_i^*=1$.  To justify this form, choose a finite
sum $\sum_i z_i d w_i$ within $1/2$ of $1$.  Writing it as $XY$, with
$X=(z_i d^{1/2})_i$ and $Y=(d^{1/2}w_i)_i^{\mathsf T}$, shows that
$XX^*=\sum_i z_i d z_i^*$ is invertible: use
$(XY)(XY)^*\le\|Y\|^2XX^*$.  Left normalization by $(XX^*)^{-1/2}$
gives the required $\xi_i$.

For a bounded trace $\sigma$ on $D$, define
\[
                 \tau_\sigma(a)=\sum_i\sigma(\xi_i^*a\xi_i),
                 \qquad a\in A_+.
\]
Every summand is in $D$, so this defines a bounded positive functional.
Inserting the finite frame on both sides of $x\in A$ yields
\[
\begin{split}
 \tau_\sigma(x^*x)
 &=\sum_{i,j}\sigma((\xi_j^*x\xi_i)^*(\xi_j^*x\xi_i))\\
 &=\sum_{i,j}\sigma((\xi_j^*x\xi_i)(\xi_j^*x\xi_i)^*)
 =\tau_\sigma(xx^*).
\end{split}
\]
Thus $\tau_\sigma$ is a trace.  For $a\in D_+$, the elements
$a^{1/2}\xi_i=a^{1/2}x_i d^{1/2}$ lie in $D$, and cyclicity gives
$\tau_\sigma(a)=\sigma(a)$.  Any trace extending $\sigma$ must have
the displayed formula, by $\sum_i\xi_i\xi_i^*=1$, so the extension
is unique.

For a tracial state $\tau$ on $A$, the norm of $\tau|_D$ is
$\mathrm d_\tau(d)=c$.  Conversely, if $\sigma\in \mathrm T(D)$ and
$\tau_\sigma=\lambda\tau$ for $\tau\in \mathrm T(A)$, then
$1=\|\sigma\|=\lambda c$, so $\lambda=1/c$.  This proves the
bijection.  The restriction map is continuous and affine, so its
compact domain makes it an affine homeomorphism onto $\mathrm T(D)$.

Let $(e_k)$ be an increasing positive contractive approximate unit
obtained by functional calculus of $d$.  The continuous functions
$\tau\mapsto\tau(e_k)$ increase pointwise to the constant $c$;
Dini's theorem gives uniform convergence.  For $l\ge k$,
\[
 \|e_l-e_k\|_{2,X}^2\le\sup_{\tau\in X}\tau(e_l-e_k)
 \le\sup_{\tau\in X}(c-\tau(e_k))\longrightarrow0,
 \qquad X=\T(A).
\]
Hence the bounded sequence has a limit $p\in\mathcal M_A$.
The same Dini argument applies to $\tau(e_k^2)\uparrow c$, so
$\|e_k-e_k^2\|_{2,X}^2\le\sup_\tau\tau(e_k-e_k^2)\to0$.
Continuity of multiplication on bounded balls gives $p^2=p$.
In each tracial GNS representation this limit is the support of
$d$, and $\tau(p)=c$ at every designated trace.  On $D$ the
normalized $2$-norm is precisely $c^{-1/2}$ times the ambient one.
The unit ball of $D$ is uniformly $2$-dense in that of
$p\mathcal M_Ap$: if $a_k\in A$ are contraction approximants to
$x\in(p\mathcal M_Ap)_1$, the elements $e_k a_k e_k\in D$ have the
same norm bound and converge to $x$ after taking sufficiently fine
approximants.  This identifies the two completions.
\end{proof}

\begin{remark}[Projection calculus in the hereditary completion]
If the ambient completion is regular, its projection-calculus
properties pass to $p\mathcal M_Ap$ with the normalized designated
traces.  Bounded-ball completeness passes to this corner, as does
real rank zero.
An ambient partial isometry whose initial and final projections are
bounded by $p$ lies in $p\mathcal M_Ap$.  To realize a normalized
profile $f$ in the corner, realize the ambient profile $cf$ and use
comparison to move the resulting projection below $p$.
All trace formulas are rescaled by the single constant $c$.
The next corollary establishes regularity directly for the
hereditary completions used below.
\end{remark}

\begin{corollary}\label{cor:audited-regular-hereditary-completion}
Let $A$ be a separable unital $\mathrm C^*$-algebra with nonempty
tracial state space and strict comparison by ordinary tracial
states.  Suppose that $A$ contains a unital copy of $\mathcal Z$.
For every nonzero $h\in\mathcal Z_+$, put
$D=\operatorname{Her}_A(h)$.  Then
$(\overline D^{\,\mathrm T(D)},X_D)$ is regular, factorial, and of
type $\mathrm{II}_1$, and its full tracial state space is $X_D$.
\end{corollary}

\begin{proof}
The element $h$ is full in $A$, and
$\mathrm d_\tau(h)=\mathrm d_{\tau_{\mathcal Z}}(h)=c>0$ is
constant on $\mathrm T(A)$.  Lemma~\ref{lem:audited-constant-rank-corner}
identifies $\mathrm T(D)$ affinely homeomorphically with
$\mathrm T(A)$ by $\tau\mapsto c^{-1}\tau|_D$.
Comparison in a hereditary subalgebra agrees with comparison in the
ambient algebra.  The trace identification therefore shows that
$D$ has strict comparison by ordinary tracial states.
The last part of the proof of Proposition~\ref{prop:z-cone}, applied
to $[h]\in\mathrm{Cu}(\mathcal Z)$, shows that the strictly positive
class $[h]\in\mathrm{Cu}(D)$ is almost divisible.

There are no nonzero finite-dimensional representations of $A$,
since their restriction to the unital copy of $\mathcal Z$ would
be a nonzero finite-dimensional representation of $\mathcal Z$.
We check the same assertion for $D$ using the unitality of $A$.
Choose the finite frame $\xi_i$ from the proof of
Lemma~\ref{lem:audited-constant-rank-corner}, so
$\sum_i\xi_i\xi_i^*=1_A$ and $\xi_i^*\xi_j\in D$.
The matrix
$P=(\xi_i^*\xi_j)_{i,j}\in\mathrm M_m(D)$ is a projection, and
\[
 a\longmapsto(\xi_i^*a\xi_j)_{i,j}
\]
is a unital isomorphism from $A$ onto
$P\mathrm M_m(D)P$.
If $\sigma:D\to\mathrm M_r$ were a nonzero representation,
$\sigma^{(m)}(P)$ would be nonzero.  Indeed the entries of $P$
generate $D$ as an ideal: for $a,b\in D$,
\[
 ab=\sum_{i,j}(a\xi_i)(\xi_i^*\xi_j)(\xi_j^*b),
\]
and all three factors of each summand belong to $D$.
The linear span of such products $ab$ is dense in $D$, so a
representation annihilating all entries of $P$ is zero.
The displayed corner isomorphism would then give a nonzero
finite-dimensional representation of $A$, a contradiction.
Thus $D$ has no such representations.  Apply
Proposition~\ref{prop:audited-regular-completion}.
\end{proof}
\section{An external AF core with uniform tracial averaging}
\label{sec:external-af-core}

Let $(\mathcal M,X)$ be a unital factorial tracially complete
$\mathrm C^*$-algebra of type $\mathrm{II}_1$.  We assume that $X$ is a
metrizable Choquet simplex, that the unit ball of $\mathcal M$ is
separable for
\[
 \|x\|_{2,X}=\sup_{\tau\in X}\tau(x^*x)^{1/2},
\]
and that $(\mathcal M,X)$ is regular.  We use the following precise
consequences of regularity:
\begin{enumerate}
\item $\mathcal M$ has real rank zero;
\item every continuous affine $f:X\to[0,1]$ is the profile
      $\tau\mapsto\tau(p)$ of a projection $p\in\mathcal M$;
\item projections with equal profiles are Murray--von Neumann
      equivalent, and profile inequality gives Murray--von Neumann
      subequivalence.
\end{enumerate}
These are the regular projection conclusions of
\cite[Theorem~1.3 and Definition~3.2]{EHSS2026}.  In particular, if $p$ is a projection
and $0\le f\le\widehat p$, where $\widehat p(\tau)=\tau(p)$,
there is $q\le p$ with $\widehat q=f$.  This assertion includes
profiles which vanish at some points of $X$: first realize $f$ in
$\mathcal M$ and then move the resulting projection into $p$ by
comparison.  This non-strict version is needed below.

If $\mathcal M$ is the bounded uniform tracial completion of a
norm-separable algebra $A$, separability of its unit ball for
$\|\cdot\|_{2,X}$ is automatic.  Here the canonical homomorphism
from $A$ to its completion need not be injective.  No real-rank
assumption on $A$ is imposed in this section.

\begin{lemma}[A dimension group of affine profiles]
\label{lem:external-profile-group}
There is a countable dense rational vector subspace
$F\subseteq\operatorname{Aff}(X)$ containing the constant function $1$.
With
\[
 F_+=\{0\}\cup\{f\in F:\min_{\tau\in X}f(\tau)>0\},
\]
$(F,F_+,1)$ is a simple dimension group with order unit, and
$X\to S(F,1)$, $\tau\mapsto(f\mapsto f(\tau))$, is an affine
homeomorphism.  Consequently there is a simple unital
infinite-dimensional AF algebra $B_0$ with this scaled ordered
$\mathrm K_0$-group and with $\mathrm T(B_0)$ identified with $X$.
\end{lemma}

\begin{proof}
Choose a countable dense subset of $\operatorname{Aff}(X)$ and take
its rational linear span together with $1$.  The strict cone is
proper and unperforated, $1$ is an order unit, and every nonzero
positive element is an order unit.  To check interpolation, suppose
$a_i\le b_j$ for $i,j\in\{1,2\}$.  If $a_i=b_j$ for one pair,
that element interpolates.  Otherwise the four differences
$b_j-a_i$ are uniformly strictly positive.  Strict affine
interpolation on the Choquet simplex gives a continuous affine
function lying strictly between the two lower and two upper
functions; a sufficiently close uniform approximation from $F$
interpolates in $F$.

An order-unit state $s$ on $F$ satisfies
$|s(f)|\le\|f\|_\infty$: for every rational
$q>\|f\|_\infty$, one has $q1\pm f\in F_+$.
Thus $s$ extends to a continuous unital linear functional on
$\operatorname{Aff}(X)$.  It is positive, since a nonnegative
continuous affine function can be approximated by strictly positive
functions in $F$.  States on $\operatorname{Aff}(X)$ are precisely
evaluations at points of $X$.  The evaluation map is continuous and
bijective from the compact space $X$ to the Hausdorff state space,
and hence is an affine homeomorphism.
Dimension-group realization and AF classification
supply $B_0$; see \cite{EHS1980,Elliott1976,Goodearl1986}.
It is simple because every nonzero positive group element is an
order unit.  It is infinite dimensional because $F$ contains the
constant rational functions and hence is divisible and noncyclic.
\end{proof}

\begin{lemma}[Profile realization and relative conjugacy]
\label{lem:external-fd-realization}
Let $D=\bigoplus_\alpha \mathrm M_{n_\alpha}$ be finite dimensional.
If strictly positive continuous affine functions $f_\alpha$ satisfy
$\sum_\alpha n_\alpha f_\alpha=1$, there is a unital embedding
$D\to\mathcal M$ whose minimal projections in the $\alpha$th
summand have profile $f_\alpha$.

Two unital homomorphisms $D\to\mathcal M$ with the same minimal
projection profiles are unitarily conjugate.  Moreover, two finite
spectral refinements of a fixed $D\subseteq\mathcal M$, with the
same spectral labels and corresponding profiles inside every old
minimal corner, are conjugate by a unitary in $D'\cap\mathcal M$.
\end{lemma}

\begin{proof}
Repeated projection realization inside the remaining projection
produces a partition of $1$ with the prescribed profiles.  Equal
profiles give partial isometries and hence matrix units in each
summand.  Strict positivity ensures that no summand is killed.

For conjugacy, write the homomorphisms as $\phi$ and $\psi$.
Choose $v_\alpha$ with
\[
 v_\alpha^*v_\alpha=\phi(e^\alpha_{11}),\qquad
 v_\alpha v_\alpha^*=\psi(e^\alpha_{11}).
\]
Then
\[
 u=\sum_{\alpha,i}
 \psi(e^\alpha_{i1})v_\alpha\phi(e^\alpha_{1i})
\]
is a unitary and $u\phi(d)u^*=\psi(d)$ for $d\in D$.
For the relative assertion, perform the corresponding partition
conjugacies inside one old minimal corner of each matrix summand,
and replicate them using the old matrix units.  The resulting
unitary commutes with $D$.
\end{proof}

Fix a unital finite-dimensional exhaustion
$B_0=\overline{\bigcup_l B_{0,l}}$.  Projection realization and
Lemma~\ref{lem:external-fd-realization} construct a unital embedding
$\iota:B_0\to\mathcal M$ preserving the prescribed trace profiles.
Indeed, embed successive stages with their $F$-profiles, and conjugate
the next-stage embedding so that it agrees with the preceding one.
The union extends to $B_0$.  Its map is injective because $B_0$ is
simple and the map is unital.

\begin{lemma}[Finite-spectrum approximation with prescribed profiles]
\label{lem:external-profile-refinement}
Let $D=\iota(B_{0,l})\subseteq\mathcal M$, and let
$c=c^*\in D'\cap\mathcal M$ have finite spectrum.  Given
$\varepsilon>0$, there is a finite-spectrum selfadjoint
$\widetilde c\in D'\cap\mathcal M$ such that
\[
 \|\widetilde c-c\|_{2,X}<\varepsilon,
 \qquad \|\widetilde c\|\le\|c\|,
\]
and, inside every minimal corner of $D$, its spectral partition has
strictly positive profiles in $F$.
There are $L>l$ and
$c_0\in D'\cap\iota(B_{0,L})$ with the same spectral labels and
corresponding minimal-corner profiles as $\widetilde c$.
\end{lemma}

\begin{proof}
It suffices to treat one minimal corner $p=e^\alpha_{11}$ and
replicate the construction over its matrix summand.  In that corner
write
\[
 pcp=\sum_{i=1}^r\lambda_i p_i,\qquad
 \sum_i p_i=p.
\]
The nonzero $p_i$ need not have strictly positive profiles.
Put $R=\|c\|$.  If $R=0$ the conclusion is immediate.
Choose $s\in(0,1)$.  Realize projections $p_i'\le p_i$ with
\[
 \widehat{p_i'}=(1-s)\widehat{p_i}.
\]
For $e=p-\sum_i p_i'$ one has $\widehat e=s\widehat p$.
Partition $e=\sum_i z_i$ with
$\widehat{z_i}=(s/r)\widehat p$, and put $q_i=p_i'+z_i$.
Then $\sum_iq_i=p$ and
\[
 \widehat{q_i}\ge(s/r)\widehat p>0.
\]
The old and new spectral elements agree on $p-e$; their difference
is supported on $e$.  Thus
\[
 \left\|\sum_i\lambda_i p_i-
             \sum_i\lambda_i q_i\right\|_{2,\tau}
 \le 2R\sqrt{s\tau(p)}.
\]

Next choose $t\in(0,1)$ and write
$\mu=\min_{i,\tau}\tau(q_i)>0$.
Since $\widehat p\in F$, density of $F$ supplies functions
$g_i\in F$ with $\sum_i g_i=\widehat p$ and
\[
 \|g_i-\widehat{q_i}\|_\infty<t\mu/2
 \quad(1\le i\le r).
\]
For completeness, approximate the first $r-1$ profiles sufficiently
closely and define the last by subtraction from $\widehat p$.
Every $g_i$ is strictly positive and
\[
 h_i:=g_i-(1-t)\widehat{q_i}>0,\qquad
 \sum_i h_i=t\widehat p.
\]
Shrink $q_i$ to $q_i'\le q_i$ with
$\widehat{q_i'}=(1-t)\widehat{q_i}$.
Partition $p-\sum_iq_i'$ into orthogonal projections $w_i$ with
$\widehat{w_i}=h_i$.  The projections
$\widetilde p_i=q_i'+w_i$ then partition $p$ and have profiles
$g_i\in F_+\setminus\{0\}$.
This second modification has $2$-norm error at most
$2R\sqrt{t\tau(p)}$.

Perform this procedure in one minimal corner of each summand and
replicate by the matrix units.  Choosing the same sufficiently small
$s,t$ in all finitely many summands gives a total error at most
$2R(\sqrt s+\sqrt t)$, since the replicated old minimal corners
partition the unit.  The spectral labels are unchanged, proving
the norm bound.

In $\mathrm K_0(B_0)=F$, the functions $g_i$ form a positive decomposition
of the old minimal projection class.  This decomposition is realized
by projections at a later finite-dimensional stage.  More explicitly,
represent all the positive classes at a common later stage and pass
farther until their sum equality holds there.  The equality is then
an equality of nonnegative integer multiplicity vectors, so the
old projection can be partitioned accordingly inside that stage.
Replication over the old matrix units realizes the entire refinement.
This supplies $c_0$ at a common stage $L>l$.
\end{proof}

\begin{theorem}[External AF core]
\label{thm:external-af-core}
Under the standing hypotheses there is a unital simple
infinite-dimensional AF algebra
\[
 B=\overline{\bigcup_nD_n}\subseteq\mathcal M,
\]
where the $D_n$ are increasing unital finite-dimensional algebras,
such that restriction $X\to \mathrm T(B)$ is an affine homeomorphism and
\begin{equation}\label{eq:external-af-averaging}
 \operatorname{dist}_{2,X}
 \bigl(\mathbb E_{D_n}(x),Z(D_n)\bigr)\longrightarrow0
 \qquad(x\in\mathcal M).
\end{equation}
Here $\mathbb E_D(x)=\int_{\mathcal U(D)}uxu^*\,du$.
For every $\rho\in X$, if
\[
 P_\rho=\pi_\rho(\mathcal M)'',\qquad
 N_\rho=\pi_\rho(B)'',
\]
then
\begin{equation}\label{eq:external-af-centres}
 N_\rho'\cap P_\rho=Z(N_\rho)=Z(P_\rho).
\end{equation}
The algebra $N_\rho$ is hyperfinite with separable predual and has
only type $\mathrm{II}_1$ factor fibres.  In particular the inclusion
is irreducible when $\rho$ is extreme in $X$.
\end{theorem}

\begin{proof}
Choose a countable $2$-dense sequence $(a_j)$ in the selfadjoint unit
ball of $\mathcal M$, a schedule $j(n)$ visiting each index
infinitely often, and positive errors $\varepsilon_n\to0$.
Starting with the trace-preserving embedding $\iota_1:B_0\to\mathcal M$
constructed above, build embeddings $\iota_n$ and stage indices
$l_n\nearrow\infty$ with
\[
 \iota_{n+1}|_{B_{0,l_n}}=\iota_n|_{B_{0,l_n}}.
\]
Put $D_n=\iota_n(B_{0,l_n})$.
The algebra $D_n'\cap\mathcal M$ is a finite direct sum of corners
of $\mathcal M$ and hence has real rank zero.  Choose a
finite-spectrum selfadjoint contraction $c$ there with
\[
 \|\mathbb E_{D_n}(a_{j(n)})-c\|<\varepsilon_n/2.
\]
Apply Lemma~\ref{lem:external-profile-refinement} with error
$\varepsilon_n/2$, obtaining $\widetilde c$ and $c_0$ at a later
stage $L$.  Relative conjugacy gives
$w\in\mathcal U(D_n'\cap\mathcal M)$ with
$wc_0w^*=\widetilde c$.
Set $\iota_{n+1}=\operatorname{Ad}(w)\iota_n$ and choose
$l_{n+1}\ge\max\{L,l_n+1,n+1\}$.  This fixes $D_n$ and places
$\widetilde c$ in $D_{n+1}$.

Eventually constant values on each abstract stage define a unital
embedding $\iota_\infty:B_0\to\mathcal M$.
All stages preserve all profiles in $F$, so restriction identifies
$X$ with $\mathrm T(B)$ for its range $B$.
For $m\ge n+1$, averaging absorption $\mathbb E_{D_m}\mathbb E_{D_n}=\mathbb E_{D_m}$
and $2$-contractivity imply
\[
 \|\mathbb E_{D_m}(a_{j(n)})-\mathbb E_{D_m}(\widetilde c)\|_{2,X}
 <\varepsilon_n.
\]
Since $\widetilde c\in D_m$, its last average lies in $Z(D_m)$.
The repeated dense schedule and $2$-contractivity prove
\eqref{eq:external-af-averaging}.

Fix $\rho$ and suppress representations from the notation.
Let $\mathbb E_N:P_\rho\to N_\rho$ be the trace-preserving normal
conditional expectation.  For $x\in N_\rho'\cap P_\rho$ put
$y=x-\mathbb E_N(x)$.  Then $y$ commutes with $N_\rho$, and
$\rho(yb)=0$ for every $b\in N_\rho$.
For $a\in\mathcal M$,
\[
 \rho(ya)=\rho(y\mathbb E_{D_n}(a)).
\]
By \eqref{eq:external-af-averaging}, Cauchy--Schwarz makes the
right side tend to zero, since its central approximants belong to
$B$.  Density in $L^2(P_\rho,\rho)$ gives $y=0$.
Thus $N_\rho'\cap P_\rho=Z(N_\rho)$.

To identify the centres, let $h\in Z(N_\rho)_+$ and $L=\|h\|$.
The positive tracial functional
$\varphi(b)=\rho(hb)$ on $B$ is dominated by $L\rho|_B$.
The affine homeomorphism $X\cong \mathrm T(B)$ extends to a homogeneous
additive bijection from the cone of designated positive traces on
$\mathcal M$ to the positive tracial functionals on $B$.
Extend $\varphi$ and $L\rho|_B-\varphi$ through that bijection.
Their extensions sum to $L\rho$, so the extension of $\varphi$ is
dominated by $L\rho$ on $\mathcal M$.
Finite von Neumann Radon--Nikodym theory therefore gives
$k\in Z(P_\rho)$, $0\le k\le L$, with
$\widetilde\varphi(a)=\rho(ka)$; domination ensures normality in
this GNS representation, and traciality makes the density central.
The relative commutant equality places $k$ in $Z(N_\rho)$.
Restriction to $B$ and $L^2$ density in $N_\rho$ imply $k=h$.
Thus $Z(N_\rho)\subseteq Z(P_\rho)$; the reverse inclusion follows
from $N_\rho'\cap P_\rho=Z(N_\rho)$.

Hyperfiniteness follows from the finite-dimensional exhaustion of
$B$, and separability of the predual follows from separability of
$B$.  A finite-dimensional factor representation of the simple
infinite-dimensional algebra $B$ would be injective and is therefore
impossible.  Hence its finite factorial fibres are all of type
$\mathrm{II}_1$.  Finally, extremality of $\rho$ corresponds to
extremality of its restriction to $B$, so both fibre algebras are
factors by \eqref{eq:external-af-centres} and the stated factoriality.
\end{proof}

\begin{remark}\label{rem:external-core-quantization-inputs}
Put $\mathcal R=\overline B^{\,2,X}\subseteq\mathcal M$, with
completion on bounded balls.  Since $X\cong\T(B)$, this is the
uniform tracial completion of $B$.  The simple AF algebra $B$ has
strict comparison by traces, and its divisible ordered group $F$
gives a divisible unit.  Thus its completion has the regular projection
calculus of \cite[Proposition~1.4 and Theorem~1.3]{EHSS2026}.
For each $\rho\in X$, its von Neumann closure is $N_\rho$.
When $\mathcal M$ is the bounded uniform tracial completion of $A$,
bounded-ball density of the canonical image of $A$ makes that image
$L^2(P_\rho,\rho)$-dense.  Consequently its von Neumann closure is
$P_\rho$, canonically identified with the tracial GNS von Neumann
algebra of $A$ at the restricted trace.  This identification does
not require the canonical map from $A$ to be injective.
The averaging conclusion is in uniform $2$-norm; all traces used
below are the designated extensions in $X$.
\end{remark}

\section{Uniform quantization and finite word cancellation}
\label{sec:quantization}

This section records all the uniform quantization arguments used below.
In particular, the constants in the small-mean estimates are allowed to
depend on the cardinality of the finite coefficient family.
No compactness or dimension assumption on the extreme boundary is used.

\paragraph{Standing external AF-completion hypotheses.}
For this section abbreviate $M=\mathcal M$ and $R=\mathcal R$.
Let $(M,X)$ be the factorial tracially complete algebra from the
preceding section.  In particular, $X$ is a compact convex face of
$\mathrm T(M)$, and its traces define a faithful uniform $2$-norm.  We use the
external AF algebra supplied by Theorem~\ref{thm:external-af-core}:
\[
 B=\overline{\bigcup_n D_n}\subseteq M,
 \qquad R=\overline B^{\,2,X}\subseteq M,
\]
where $B$ is unital, simple, and infinite dimensional, the $D_n$ are
increasing unital finite-dimensional subalgebras, and the latter
completion is taken on bounded balls.  Restriction $X\to \mathrm T(B)$ is
an affine homeomorphism.  We write
\[
 \|x\|_{2,u}=\|x\|_{2,X}
   =\sup_{\tau\in X}\tau(x^*x)^{1/2},
 \qquad \|x\|_{1,\tau}=\tau(|x|).
\]
The following inputs have been established above; they are recorded
here to specify precisely what the arguments use.
\begin{enumerate}
\item[(AF)] If $\mathbb E_n$ is averaging over the compact group
 $U(D_n)$, then
 \[
  \operatorname{dist}_{2,X}(\mathbb E_n(x),Z(D_n))\longrightarrow0
  \qquad(x\in M).
  \tag{AF}
 \]
\item[(C)] The unit balls of $M$ and $R$ are complete for $\|\cdot\|_{2,u}$,
 and the unit ball of $B$ is dense in that of $R$ in this metric.
 The uniform $2$-norm is faithful on $R$.
 Each function $\tau\mapsto\tau(x)$, $x\in M$, is continuous affine on $X$.
\item[(P)] The algebra $R$ has real rank zero.  If $e\in R$ is a projection
 and $f\in\operatorname{Aff}(X)$ satisfies
 \[
     0<f(\tau)<\tau(e)\qquad(\tau\in X),
 \]
 there is a projection $g\le e$ in $R$ with $\tau(g)=f(\tau)$ for every
 $\tau\in X$.  Projections in $R$ with identical trace profiles are
 Murray--von Neumann equivalent in $R$.
\item[(F)] For $\rho\in X$, set
 $P_\rho=\pi_\rho(M)''$ and $N_\rho=\pi_\rho(B)''$.
 Then
 \[
  N_\rho'\cap P_\rho=Z(N_\rho)=Z(P_\rho).
  \tag{F}
 \]
 In the common central decomposition the factor fibres of $N_\rho$
 are hyperfinite $\mathrm{II}_1$ factors.
\end{enumerate}
In (C), density is density of the indicated bounded balls; selfadjoint
contractions may therefore be approximated by selfadjoint contractions.
The same conclusion for arbitrary bounded elements follows by rescaling.
The continuity assertion for $x\in R$ also follows directly from
$|\tau(x-y)|\le\|x-y\|_{2,u}$ and bounded-ball approximation by $B$;
for $x\in M$, it is part of the designated tracial completion structure.
The facial property of $X$ ensures that central reweightings of any
$\rho\in X$ in its GNS von Neumann algebra still belong to $X$;
this point will be used explicitly in the compression argument.
All coefficients below belong to $M$, and all newly constructed
projections and partial isometries belong to $R$.

We first prove linear projection equidistribution.  This also explains
why no projection-transfer assertion is an additional hypothesis.

\begin{lemma}[Global fractional projections]
\label{lem:global-fractional-projection}
Let $0<\lambda<1$, let $\mathcal X\subset M$ and
$\mathcal H\subset R$ be finite, and let $\eta>0$.
There is a projection $e\in R$ such that
\[
 \tau(e)=\lambda,\qquad
 |\tau(ex)-\lambda\tau(x)|<\eta\quad(x\in\mathcal X,\ \tau\in X),
 \qquad
 \|[e,h]\|_{2,u}<\eta\quad(h\in\mathcal H).
 \tag{\thesection.1}\label{eq:q-1}
\]
\end{lemma}

\begin{proof}
Approximate the elements of $\mathcal H$ in uniform $2$-norm by
elements of a fixed $D_k$.  For sufficiently large $n\ge k$, (AF)
gives approximants in $Z(D_n)$, in uniform $2$-norm, for the averages
of the finitely many elements of $\mathcal X$.

Write $D_n=\bigoplus_l \mathrm M_{d_l}$ and choose matrix units
$(e^l_{ij})$ in each summand.  By (P), choose a projection
$f_l\le e^l_{11}$ in $R$ with
$\tau(f_l)=\lambda\tau(e^l_{11})$ for all $\tau\in X$.
The strict inequalities needed in (P) hold because each nonzero
projection of $B$ has strictly positive trace at every point of $X$.
Set
\[
 e=\sum_l\sum_{i=1}^{d_l}e^l_{i1}f_l e^l_{1i}.
\]
This is a projection commuting with $D_n$, and
\[
 \tau(ez)=\lambda\tau(z)\qquad(z\in Z(D_n),\ \tau\in X).
 \tag{\thesection.2}\label{eq:q-2}
\]
In particular $\tau(e)=\lambda$.

For $a\in M$, commutation with $D_n$ and traciality give
$\tau(ea)=\tau(e\mathbb E_n(a))$ and
$\tau(a)=\tau(\mathbb E_n(a))$.
Consequently, for $z\in Z(D_n)$,
\[
 |\tau(ea)-\lambda\tau(a)|
 =|\tau((e-\lambda1)(\mathbb E_n(a)-z))|
 \le(1+\lambda)\|\mathbb E_n(a)-z\|_{2,u}.
 \tag{\thesection.3}\label{eq:q-3}
\]
The last inequality is tracial Cauchy--Schwarz, uniformly over
$\tau\in X$.  If $h_0\in D_k$ approximates $h\in R$, then
$\|[e,h]\|_{2,u}\le2\|h-h_0\|_{2,u}$.
Choosing the approximations sufficiently close proves \eqref{eq:q-1}.
\end{proof}

\begin{lemma}[Uniform $2$-rounding in a fixed corner]
\label{lem:corner-rounding}
Let $p\in R$ be a projection.  For every positive contraction
$a\in pRp$ and every $\epsilon>0$, there is a projection $r\in pRp$
such that
\[
 \|r-a\|_{2,u}\le 2\|a-a^2\|_{2,u}+\epsilon.
 \tag{\thesection.4}\label{eq:q-4}
\]
\end{lemma}

\begin{proof}
By real rank zero, approximate $a$ arbitrarily closely in norm by a
finite-spectrum positive contraction $b\in pRp$.
Let $r=1_{[1/2,1]}(b)$.  The scalar inequality
$|1_{[1/2,1]}(t)-t|\le2t(1-t)$, $0\le t\le1$, gives
$\|r-b\|_{2,u}\le2\|b-b^2\|_{2,u}$.
Since $\|b-b^2-(a-a^2)\|\le3\|a-b\|$, norm approximation yields \eqref{eq:q-4}.
\end{proof}

\begin{lemma}[Linear projection equidistribution in a fixed corner]
\label{lem:linear-projection-equidistribution}
Let $p\in R$ be a projection of constant trace $\beta>0$.
Given $0<t<\beta$, a finite set $\mathcal X\subset pMp$, and
$\eta>0$, there is a projection $r\le p$ in $R$ satisfying
\[
 \tau(r)=t,\qquad
 \left|\tau(rx)-\frac{t}{\beta}\tau(x)\right|<\eta
 \quad(x\in\mathcal X,\ \tau\in X).
 \tag{\thesection.5}\label{eq:q-5}
\]
\end{lemma}

\begin{proof}
Rescale to suppose the elements of $\mathcal X$ are contractions.
Choose a small $\zeta>0$ with $\zeta<t/2$, and put
$\lambda=(t-\zeta)/\beta\in(0,1)$.
Apply Lemma~\ref{lem:global-fractional-projection} to $\mathcal X\cup\{p\}$
and $\mathcal H=\{p\}$, with an arbitrarily small error tolerance.
For the resulting projection $e$, put $a=pep\in pRp$.
Then
\[
 \tau(a)\approx t-\zeta,
 \qquad
 \|a-a^2\|_{2,u}
 =\|pe(1-p)ep\|_{2,u}
 \le\|(1-p)ep\|_{2,u}
 \le\|[e,p]\|_{2,u}.
 \tag{\thesection.6}\label{eq:q-6}
\]
By Lemma~\ref{lem:corner-rounding}, the projection $r_0\le p$ obtained
by rounding $a$ can be made arbitrarily close to $a$ in uniform
$2$-norm.  Choose the errors in \eqref{eq:q-6} and in rounding so small that
\[
 t-\tfrac32\zeta<\tau(r_0)<t-\tfrac12\zeta
 \qquad(\tau\in X).
 \tag{\thesection.7}\label{eq:q-7}
\]
The function $f(\tau)=t-\tau(r_0)$ is continuous affine and satisfies
\[
 0<f(\tau)<\tau(p-r_0),
\]
the last strict inequality being exactly $t<\beta$.
Thus (P) gives a projection $g\le p-r_0$ in $R$ with trace profile $f$.
Set $r=r_0+g$.  Then $\tau(r)=t$ and
\[
 \|r-a\|_{2,u}
 \le\|r_0-a\|_{2,u}+\sqrt{3\zeta/2}.
 \tag{\thesection.8}\label{eq:q-8}
\]
For $x\in pMp$, traciality gives $\tau(ax)=\tau(ex)$, while
\[
 \left|\left(\lambda-\frac t\beta\right)\tau(x)\right|
 \le\zeta
\]
for contractions $x\in pMp$.
Use \eqref{eq:q-1}, \eqref{eq:q-8}, and tracial Cauchy--Schwarz, first choosing $\zeta$ small
and then the errors in \eqref{eq:q-6} small.  This proves \eqref{eq:q-5}.
\end{proof}

We next establish the pointwise compression estimate using only a
finite-matrix second-moment calculation.  In particular, no uniform
freeness theorem is being invoked.

\begin{lemma}[A finite-matrix compression calculation]
\label{lem:matrix-compression}
Let $F\cong \mathrm M_d$, $d\ge2$, be a unital matrix algebra in a finite
von Neumann algebra $(P,\rho)$, and let $q_0\in F$ have matrix rank
$1\le r<d$.
Put $\delta=r/d$, and integrate over $u\in U(F)$ with normalized Haar
measure.  Then
\[
 \int u q_0u^*\,x\,u q_0u^*\,du
 =a_d x+b_d \mathbb E_{F'\cap P}(x),
 \quad
 a_d=\frac{r^2-\delta}{d^2-1},\quad b_d=\delta-a_d.
 \tag{\thesection.9}\label{eq:q-9}
\]
In particular $0\le a_d\le\delta^2$.
For selfadjoint $x$,
\[
 \int\rho((u q_0u^* x u q_0u^*)^2)\,du
 =a_d\rho(x^2)+b_d\|\mathbb E_{F'\cap P}(x)\|_{2,\rho}^2.
 \tag{\thesection.10}\label{eq:q-10}
\]
\end{lemma}

\begin{proof}
Identify $P$ with $\mathrm M_d\,\overline\otimes\,(e_{11}Pe_{11})$.
The matrix-factor averaging map is equivariant for conjugation by
$U(d)$, hence is scalar on the scalar matrices and on the traceless
matrices.  Its value at $1$ is $\delta1$.
The trace of this linear operator on the $d^2$-dimensional matrix space
is $r^2$, because the map $X\mapsto q_0Xq_0$ has rank $r^2$ and is a
projection.  Thus its traceless eigenvalue is
$(r^2-\delta)/(d^2-1)$, giving \eqref{eq:q-9}, coefficientwise in the second factor.
Taking the trace after multiplication by $x$, and using
$\rho(\mathbb E_{F'\cap P}(x)x)=\|\mathbb E_{F'\cap P}(x)\|_2^2$ for selfadjoint $x$,
gives \eqref{eq:q-10}.
\end{proof}

\begin{lemma}[Pointwise small-rank compression]
\label{lem:pointwise-compression}
Let $p\in R$ have constant positive trace.  Normalize the traces in
the corner by $\tau_p=\tau/\tau(p)$.
Let $x_1,\ldots,x_m\in(pMp)_{\rm sa}$ be contractions with
$|\tau_p(x_i)|\le\kappa$ for every $\tau\in X$.
Fix a rational $0<\delta<1$ and a projection $r\le p$ in $R$ with
$\tau_p(r)=\delta$ for all $\tau$.
For every $\rho\in X$,
\[
 \inf_{u\in U(pRp)}
 \sum_{i=1}^m\rho_p(|u r u^*x_i u r u^*|)
 \le m\delta\kappa+m\delta^{3/2}.
 \tag{\thesection.11}\label{eq:q-11}
\]
\end{lemma}

\begin{proof}
We first justify the full-corner and central-mean assertions, without
requiring simultaneous pointwise evaluations of elements of $M$ in a
central decomposition.
Write $P=P_\rho$, $N=N_\rho$, and $Z_0=Z(P)=Z(N)$.
For $0\ne h\in(Z_0)_+$, the formula
\[
 \sigma_h(x)=\frac{\rho(h\pi_\rho(x))}{\rho(h)}\qquad(x\in M)
\]
defines a tracial state of $M$.  It belongs to $X$: it is dominated
by $c\rho$, where $c=\|h\|/\rho(h)\ge1$; if $c>1$, write
$\rho=c^{-1}\sigma_h+(1-c^{-1})\psi$ with $\psi\in \mathrm T(M)$ and use
that $X$ is a face, while $c=1$ gives $\sigma_h=\rho$.
On bounded sets this formula is uniform-$2$ continuous, since
\[
 \sigma_h(x^*x)\le
 \frac{\|h\|}{\rho(h)}\rho(x^*x)
 \le\frac{\|h\|}{\rho(h)}\|x\|_{2,u}^2.
 \tag{\thesection.12a}\label{eq:q-12a}
\]
Constant trace of $p$ therefore implies
$\rho(hp)=\beta\rho(h)$ for all $h\in(Z_0)_+$, where
$\beta=\tau(p)>0$.  Equivalently $\mathbb E_{Z_0}(p)=\beta1$.
Thus $p$ has full central support in both $N$ and $P$.

The full-corner identifications give
\[
 Z(pPp)=pZ_0=Z(pNp),\qquad
 (pNp)'\cap pPp=pZ_0.
 \tag{\thesection.12b}\label{eq:q-12b}
\]
Here is an explicit verification of the relative commutant assertion.
Finite von Neumann comparison, applied to the constant centre-valued
rank $\beta$, provides a finite frame $v_i\in N$ with
$v_i=v_ip$ and $\sum_i v_i v_i^*=1$.
Indeed, partition $1$ into finitely many projections of
centre-valued traces at most $\beta1$, compare each with $p$, and
choose the corresponding partial isometries with those mutually
orthogonal range projections.
For $y\in(pNp)'\cap pPp$, put $Y=\sum_i v_i y v_i^*$.
Inserting the frame once and commuting $y$ with
$p v_j^* n v_i p\in pNp$ shows $nY=Yn$ for every $n\in N$.
Also $pYp=y$, by the same frame identity.
Thus $Y\in N'\cap P=Z_0$, proving the second equality in \eqref{eq:q-12b}.
The same argument with $N=P$ proves the centre identification.

Now work with the normalized corner trace $\rho_p=\rho/\beta$, and
write $Z=pZ_0$ for its centre.
If $h'=ph\in Z_+$, with $h\in(Z_0)_+$, the small-mean hypothesis
applied to $\sigma_h$ gives
\[
 |\rho_p(h'x_i)|
 =\beta^{-1}|\rho(hx_i)|
 \le\kappa\rho(h)
 =\kappa\rho_p(h').
\]
Testing all positive central $h'$ proves
\[
 \|\mathbb E_Z(x_i)\|\le\kappa.
 \tag{\thesection.12}\label{eq:q-12}
\]

We use the standard structure of a finite hyperfinite von Neumann
algebra with separable predual and only $\mathrm{II}_1$ factor fibres
\cite[Chapter~XIV, Sections~1--2, and Chapter~XVI, Section~1]{Takesaki2003}:
it is $Z\,\overline\otimes\,\mathcal R$, where $\mathcal R$ is the
hyperfinite $\mathrm{II}_1$ factor.
Consequently it contains an increasing sequence of full matrix
algebras $F_j$ whose matrix sizes are divisible by the denominator of
$\delta$, such that
\[
 \left(\bigcup_jF_j\right)'\cap pP_\rho p=Z.
 \tag{\thesection.13}\label{eq:q-13}
\]
For clarity, take full matrix approximants of a copy of $\mathcal R$;
their commutant in $pP_\rho p$ equals the commutant of
$Z\vee\mathcal R=pN_\rho p$, because $Z$ is central in $pP_\rho p$.
The last commutant is $Z$ by (F).

Write $x_i=a_i+y_i$, where $a_i=\mathbb E_Z(x_i)$.
Then $\mathbb E_Z(y_i)=0$ and $\rho_p(y_i^2)\le\rho_p(x_i^2)\le1$.
On the Hilbert space $L^2(pP_\rho p,\rho_p)$, the inner product is
$\langle x,y\rangle=\rho_p(y^*x)$.  The trace-preserving conditional
expectations onto $F_j'\cap pP_\rho p$ are the orthogonal projections
onto the corresponding $L^2$-subspaces.  These projections decrease
to the orthogonal projection onto $L^2(Z,\rho_p)$, by
\eqref{eq:q-13}.  Thus the expectations converge in
$\|\cdot\|_{2,\rho_p}$ to $\mathbb E_Z$.
Choose a projection of rank fraction $\delta$ in $F_j$, and average its
unitary conjugates $q$ within $F_j$.
Equations \eqref{eq:q-9}--\eqref{eq:q-10}, followed by Cauchy--Schwarz and Jensen, give
\[
 \limsup_j\int \rho_p(|qy_iq|)\,dq\le\delta^{3/2}.
 \tag{\thesection.14}\label{eq:q-14}
\]
Here $\rho_p(q)=\delta$, so
$\rho_p(|qy_iq|)\le\sqrt\delta\,\|qy_iq\|_{2,\rho_p}$.
Since $a_i$ is central and $\|a_i\|\le\kappa$,
$\rho_p(|qa_iq|)\le\delta\kappa$.
Summing and taking one conjugate with cost at most the average proves
the desired bound with an arbitrarily small additional error in
$pN_\rho p$.

Both this $q$ and the image of $r$ have constant centre-valued rank
$\delta$, so they are unitarily conjugate in $pN_\rho p$.
For $r$, the centre-valued rank assertion follows by repeating the
dominated-central-trace argument above with its constant trace profile
$\tau(r)=\delta\beta$.
The unitaries of $pRp$ are $2$-dense in the unitaries of $pN_\rho p$.
One way to see the latter assertion is to write a von Neumann unitary
as $e^{ih}$ with bounded selfadjoint $h$, apply Kaplansky density to
$h$, and exponentiate.
The cost in \eqref{eq:q-11} is $2$-continuous on bounded sets.
Approximation therefore proves \eqref{eq:q-11} for unitaries in $pRp$.
\end{proof}

The next elementary avoidance fact is used when turning a convex
mixture of compressions into a single compression.

\begin{lemma}[Avoiding a small projection]
\label{lem:avoid-small-projection}
Let $a,P\in R$ be projections of constant traces $\delta$ and $s$,
where $s<\delta/2$.  For every $\nu>0$, there are a projection
$h\le1-P$ of constant trace $\delta/2$ and a unitary $U\in R$ such that
\[
 h\le UaU^*,\qquad \|U-1\|<\nu.
 \tag{\thesection.15}\label{eq:q-15}
\]
\end{lemma}

\begin{proof}
By real rank zero, the hereditary algebra generated by $aPa$ has an
approximate unit of projections.  Choose such a projection $e\le a$
with $\|(a-e)aPa(a-e)\|$ arbitrarily small.
For every trace, $\tau(e)\le\tau(P)$: in the corresponding finite
von Neumann algebra the support of $aPa$ is equivalent to a
subprojection of $P$.
Thus $\tau(a-e)\ge\delta-s>\delta/2$.
By (P), choose $t\le a-e$ of constant trace $\delta/2$.
It satisfies $\|Pt\|^2=\|tPt\|$ arbitrarily small.

The polar part of $(1-P)t$ is a partial isometry in $R$ with initial
projection $t$ and final projection $h\le1-P$, and $h$ is arbitrarily
close to $t$ in norm.  For completeness, the polar part is
\[
 (1-P)t\,[t(1-P)t]^{-1/2},
\]
where the inverse is taken in $tRt$; it exists when $\|Pt\|<1$.
The standard close-projection unitary, obtained from the polar part
of $ht+(1-h)(1-t)$, conjugates $t$ to $h$ and tends to $1$ as
$\|h-t\|\to0$.
Choose this unitary $U$ with $\|U-1\|<\nu$.
Since $t\le a$, \eqref{eq:q-15} follows.
\end{proof}

\begin{proposition}[Uniform relative $L^1$ quantization]
\label{prop:uniform-L1-quantization}
Let $p\in R$ have constant positive trace, let
$x_1,\ldots,x_m\in(pMp)_{\rm sa}$ be contractions, and suppose
\[
 |\tau(x_i)|\le\kappa\tau(p)\qquad(i,\ \tau\in X).
\]
For every $\epsilon>0$, there is a projection $q\le p$ in $R$ of
constant positive trace such that
\[
 \|qx_iq\|_{1,\tau}
 \le(2m\kappa+\epsilon)\tau(q)
 \qquad(i,\ \tau\in X).
 \tag{\thesection.16}\label{eq:q-16}
\]
In particular the small-mean constant can be taken to be $C_m=2m$.
\end{proposition}

\begin{proof}
Normalize traces in $pMp$ and write them simply as $\tau$ during this
proof, so $\tau(p)=1$.
Choose a small rational rank $\delta>0$ and a projection $r\le p$
in $R$ of constant rank $\delta$, using (P).
For $u\in U(pRp)$ put $q_u=uru^*$ and
\[
 g_u(\tau)=\sum_{i=1}^m\tau(|q_ux_iq_u|).
\]
These are continuous affine functions on the compact convex set $X$.
By Lemma~\ref{lem:pointwise-compression},
\[
 \sup_{\tau\in X}\inf_u g_u(\tau)
 \le m\kappa\delta+m\delta^{3/2}.
\]
Choose $\delta$ so small that $m\sqrt\delta<\epsilon_0$.
The elementary convex minimax theorem, applied to finite probability
mixtures of the $g_u$, gives unitaries $u_1,\ldots,u_L$ and weights
$\alpha_j>0$, $\sum_j\alpha_j=1$, such that, with $q_j=u_jru_j^*$,
\[
 \sum_j\alpha_j\sum_i\tau(|q_jx_iq_j|)
 \le(m\kappa+2\epsilon_0)\delta\qquad(\tau\in X).
 \tag{\thesection.17}\label{eq:q-17}
\]
Only the trace variable is required to be compact; the other convex
variable is the space of finitely supported probability measures.
Equivalently, \eqref{eq:q-17} follows by separating the convex hull of the
continuous affine costs from a strict constant upper bound and using
the barycentre of the separating probability measure on $X$.

Fix $0<c<1/4$.  We recursively construct orthogonal projections
$r_j\in R$ of constant ranks $c\alpha_j\delta$.
Let $P_{j-1}=\sum_{k<j}r_k$, whose rank is at most $c\delta$.
Apply Lemma~\ref{lem:avoid-small-projection} to $q_j$ and $P_{j-1}$.
It gives a projection $h_j\le p-P_{j-1}$ of rank $\delta/2$,
contained in a conjugate $q'_j=U_jq_jU_j^*$, where $U_j$ can be chosen
arbitrarily close to $p$ in the corner.
Choose the conjugating unitaries so close to the corner unit that,
uniformly over $\tau$ and $i$,
\[
 \left|\tau(|q'_jx_iq'_j|)-\tau(|q_jx_iq_j|)\right|<\zeta,
 \tag{\thesection.18}\label{eq:q-18}
\]
where $\zeta>0$ will be chosen after the finite mixture is fixed.
Indeed the difference in the two trace norms is at most the operator
norm of $q'_jx_iq'_j-q_jx_iq_j$.

Apply Lemma~\ref{lem:linear-projection-equidistribution} in $h_j$
with ratio $2c\alpha_j$ to the positive coefficients
\[
 h_j|q'_jx_iq'_j|h_j,
 \qquad h_jx_i r_kx_i h_j\quad(k<j),
 \tag{\thesection.19}\label{eq:q-19}
\]
using error at most $\zeta$ in each linear test.
The resulting $r_j\le h_j$ has
rank $c\alpha_j\delta$ and is orthogonal to all previous $r_k$.
For a selfadjoint $H$ and a projection $e$,
\[
 \tau(|eHe|)\le\tau(e|H|),
 \tag{\thesection.20}\label{eq:q-20}
\]
as follows by writing $H=H_+-H_-$ and using the trace-norm triangle
inequality.  Thus \eqref{eq:q-17}--\eqref{eq:q-20} give for each $i$
\[
 \sum_j\|r_jx_ir_j\|_{1,\tau}
 \le2c(m\kappa+2\epsilon_0)\delta+(2c+L)\zeta.
 \tag{\thesection.21}\label{eq:q-21}
\]
The cross requirements in \eqref{eq:q-19} give
\[
 \tau(r_jx_i r_kx_i)
 \le2c\alpha_j\tau(r_k)+\zeta
 \qquad(k<j).
\]
Put $q=\sum_jr_j$, of constant rank $c\delta$.
By orthogonality of matrix blocks in tracial $L^2$,
\[
 \left\|qx_iq-\sum_jr_jx_ir_j\right\|_{2,\tau}^{\,2}
 =2\sum_{k<j}\tau(r_jx_i r_kx_i)
 \le2c^2\delta+L(L-1)\zeta.
 \tag{\thesection.22}\label{eq:q-22}
\]
Here we used
$\sum_{k<j}\alpha_j\alpha_k\le1/2$.
The operator in \eqref{eq:q-22} is supported on $q$, so its $1$-norm is at most
$\sqrt{\tau(q)}$ times its $2$-norm.
Combining \eqref{eq:q-21}--\eqref{eq:q-22} and dividing by $c\delta$ yields
\[
 \frac{\|qx_iq\|_{1,\tau}}{\tau(q)}
 \le2m\kappa+4\epsilon_0+
 \frac{(2c+L)\zeta}{c\delta}
 +\left(2c+\frac{L(L-1)\zeta}{c\delta}\right)^{1/2}.
\]
Choose $\epsilon_0$ and $c$ so that
$4\epsilon_0+\sqrt{2c}<\epsilon/2$, and then choose $\zeta>0$
small enough to obtain \eqref{eq:q-16}.
\end{proof}

The following local construction combines quantization with the
cancellation of same-sign words.  Taking the absolute values
\emph{after} the $q_0$-compression is essential to preserve the mean
term when thinning.

\begin{lemma}[Quantization with analytic cancellation]
\label{lem:quantization-analytic}
Let $p\in R$ have constant trace $\beta>0$.
Let $z_1,\ldots,z_m\in pMp$ be contractions satisfying
$|\tau(z_i)|\le\mu\beta$ for every $\tau\in X$.
Fix finitely many templates
\[
 b_0w^\sigma b_1\cdots w^\sigma b_k,
 \qquad 1\le k\le n,\quad\sigma\in\{1,-1\},
 \tag{\thesection.23}\label{eq:q-23}
\]
whose coefficients are contractions in $pMp$.
For every $\eta,b>0$ there are a projection $q\le p$ in $R$ of
constant trace $0<\Delta<b$ and a unitary $w\in qRq$ such that
\[
 \|qz_iq\|_{1,\tau}\le(C\mu+\eta)\Delta,
 \qquad
 |\tau(b_0w^\sigma b_1\cdots w^\sigma b_k)|\le\eta\Delta,
 \tag{\thesection.24}\label{eq:q-24}
\]
where one can take $C=4C_{2m}=16m$.
If necessary pad the $z_i$-list by zeros so that $m\ge1$.
\end{lemma}

\begin{proof}
Apply Proposition~\ref{prop:uniform-L1-quantization} to
$\operatorname{Re}z_i,\operatorname{Im}z_i$, with error $\epsilon_0>0$.
It gives a projection $q_0\le p$ of constant trace $r>0$ such that,
for every one of these selfadjoint contractions $x$,
\[
 \tau(|q_0xq_0|)\le(C_{2m}\mu+\epsilon_0)r.
 \tag{\thesection.25}\label{eq:q-25}
\]
Set $N=n+1$, and choose $0<c<1/2$ small enough that $cr<b$.
We shall select $N$ orthogonal projections $r_j\le q_0$, each of
constant trace $d_0=cr/N$.

For every template \eqref{eq:q-23}, put the coefficients
$b_1,\ldots,b_{k-1},b_kb_0$ in a list $\mathcal B$.
In particular the cyclic wrap coefficient is included.
At step $j$ put
\[
 h_j=q_0-\sum_{i<j}r_i,
 \qquad
 \theta_j=\frac{d_0}{r-(j-1)d_0}\le\frac{2c}{N}.
\]
Use linear projection equidistribution in $h_j$, with ratio
$\theta_j$, for the finite positive coefficient list
\[
 h_j|q_0xq_0|h_j,
 \quad h_j b r_i b^*h_j,
 \quad h_j b^*r_i b h_j
 \quad(b\in\mathcal B,\ i<j),
 \tag{\thesection.26}\label{eq:q-26}
\]
and tolerance $\zeta>0$.
It gives $r_j\le h_j$ of rank $d_0$.
Put $q=\sum_jr_j$, so $\Delta=cr$.

For $H=q_0xq_0$, \eqref{eq:q-26} and \eqref{eq:q-20} imply
\[
 \tau(|qxq|)\le\tau(q|H|)
 \le2c\tau(|H|)+N\zeta.
\]
Splitting $z_i$ into real and imaginary parts, \eqref{eq:q-25} yields
\[
 \frac{\|qz_iq\|_{1,\tau}}{\Delta}
 \le4C_{2m}\mu+4\epsilon_0+\frac{2N\zeta}{cr}.
 \tag{\thesection.27}\label{eq:q-27}
\]
This is the step that prevents amplification of the mean by $1/c$.

For $b\in\mathcal B$, put
$D_b=qbq-\sum_jr_jbr_j$.
The two cross orientations in \eqref{eq:q-26} imply
\[
 \|D_b\|_{2,\tau}^{\,2}\le2c^2r+N^2\zeta,
 \qquad
 \frac{\|D_b\|_{1,\tau}}{\Delta}
 \le\left(2c+\frac{N^2\zeta}{cr}\right)^{1/2}.
 \tag{\thesection.28}\label{eq:q-28}
\]
For the first inequality sum the orthogonal blocks $r_i b r_j$,
$i\ne j$; for each unordered pair, \eqref{eq:q-26} bounds the two squared
$2$-norms by $2\theta_jd_0+2\zeta$.
The second inequality follows because $D_b$ is supported on $q$.

By (P), the equal-rank $r_j$ admit matrix units $(e_{ij})$ in $R$
with $e_{jj}=r_j$.
Let $w=e_{21}+e_{32}+\cdots+e_{N,N-1}+e_{1N}$.
It is a unitary in $qRq$ that cyclically permutes the $r_j$.
After cyclically moving $b_0$ to the end of \eqref{eq:q-23}, replace each
$qbq$ by its block diagonal $\sum_jr_jbr_j$.
The block-diagonal map $b\mapsto\sum_jr_jbr_j$ is contractive,
so all replacement factors are contractions.
The resulting product has degree $\sigma k$ modulo $N$ in the cyclic
matrix-block grading.  As $1\le k<N$, its diagonal blocks and its
trace vanish.
Telescoping the replacements gives a trace error at most
\[
 n\left(2c+\frac{N^2\zeta}{cr}\right)^{1/2}\Delta.
 \tag{\thesection.29}\label{eq:q-29}
\]
First choose $\epsilon_0$ small, then $c$ small, and finally $\zeta$
small, to make the residual terms in \eqref{eq:q-27} at most $\eta$ and \eqref{eq:q-29}
at most $\eta\Delta$.
All choices are finite, and $\Delta$ can be smaller than any prescribed
positive $b$.  This proves \eqref{eq:q-24}.
\end{proof}

\begin{proposition}[Finite weakened-word cancellation]
\label{prop:finite-weakened-word}
Fix finitely many formal words
\[
 W(t)=a_0t^{\epsilon_1}a_1\cdots t^{\epsilon_k}a_k,
 \qquad 1\le k\le n,\quad\epsilon_i\in\{1,-1\},
 \tag{\thesection.30}\label{eq:q-30}
\]
whose coefficients are contractions in $M$.
For every $\varepsilon>0$ there is $\kappa_0>0$, depending only on
$\varepsilon$ and the finite word patterns (including their number),
such that the following implication holds:
\[
 \sup_{\tau\in X}|\tau(a_i)|\le\kappa_0
 \quad\hbox{whenever }\epsilon_i\ne\epsilon_{i+1}
 \quad\Longrightarrow\quad
 \exists t\in U(R)\ \sup_{\tau\in X}|\tau(W(t))|<\varepsilon
 \tag{\thesection.31}\label{eq:q-31}
\]
for every prescribed word.
No condition is imposed on the outer coefficients, or on coefficients
between equal signs.
\end{proposition}

\begin{proof}
We give the finite closure, constants, maximal-patching argument, and
limit argument explicitly.

For every original word \eqref{eq:q-30}, include every interval of letter
positions $i\le j$, with either inherited or unit boundary
coefficients:
\[
 b_Lt^{\epsilon_i}a_i\cdots a_{j-1}t^{\epsilon_j}b_R,
 \qquad b_L\in\{1,a_{i-1}\},\quad b_R\in\{1,a_j\}.
 \tag{\thesection.32}\label{eq:q-32}
\]
Call the resulting finite indexed family $\mathcal F$, retaining
repeated templates even when some coefficients coincide.
It contains the original family and is closed under the same interval
operation.  Each of its internal coefficients is an original internal
coefficient at the same pair of signs.  Its indexed sign patterns
therefore depend only on those of the original words.

Let $m\ge1$ be the number of pairs of positions $i<j$ with opposite
signs, summed over all words in $\mathcal F$, using a dummy zero entry
if this number is zero.  This is an upper bound for all gap-list
cardinalities below, independent of any subsequently chosen partial
unitary.  Pad by zeros to keep this cardinality fixed.
Let $C=4C_{2m}$ be the constant in
Lemma~\ref{lem:quantization-analytic}, and put
\[
 B_0=2^n,\qquad D=B_0C,\qquad K=\lceil2D\rceil+1.
 \tag{\thesection.33}\label{eq:q-33}
\]
It suffices to consider $0<\varepsilon\le1$.
Set
\[
 \gamma=\frac{\varepsilon}{4n},\quad S=1-\gamma,\quad
 L=\gamma^{-K}-1,\quad
 a=\kappa_0=\frac{\varepsilon}{8L},\quad
 \eta=\frac{a}{2B_0}.
 \tag{\thesection.34}\label{eq:q-34}
\]
Suppose the opposite-sign coefficient traces are bounded by
$\kappa\le\kappa_0$.  Define the increasing convex budget
\[
 F(s)=(\kappa+a)((1-s)^{-K}-1),\qquad 0\le s\le S.
 \tag{\thesection.35}\label{eq:q-35}
\]

Consider normal partial unitaries $v\in R$ whose support
$f=v^*v=vv^*$ has constant trace $s\le S$, and satisfying
\[
 |\tau(W(v))|\le F(s)\qquad(W\in\mathcal F,\ \tau\in X).
 \tag{\thesection.36}\label{eq:q-36}
\]
Here $v^{-1}$ means $v^*$.
Order them by orthogonal extension: $v\preceq u$ if $u=v+w$, where
$w$ is a normal partial unitary with support orthogonal to that of $v$.
The zero partial unitary belongs to this partially ordered set.

Every chain has an upper bound.  To verify this, let $s$ be the
supremum of its support ranks.
If a member attains this rank, it is an upper bound: comparable
members of equal rank differ by a partial unitary of zero uniform
$2$-norm and therefore agree.
Otherwise choose an increasing sequence $v_j$ in the chain with
support ranks $s_j\uparrow s$.  Orthogonal extension gives
\[
 \|v_j-v_i\|_{2,u}^2=s_j-s_i\qquad(j\ge i).
 \tag{\thesection.37}\label{eq:q-37}
\]
By (C), $v_j$ converges in the unit ball of $R$ to $v$.
Adjoints and multiplication on bounded sets are uniformly
$2$-continuous.  Passing the partial-unitary identities to the limit
shows that $v$ is a normal partial unitary of constant support rank
$s$.  Uniform tracial continuity of each finite word and continuity of
$F$ show that \eqref{eq:q-36} passes to the limit.
For any chain member $w$ of rank less than $s$, some $v_j$ extends it;
passing $v_j\operatorname{supp}(w)=w=
\operatorname{supp}(w)v_j$ to the limit shows that $v$ extends $w$.
This proves the upper-bound assertion.  Zorn's lemma therefore gives
a maximal $v$.

Suppose its rank is $s<S$, and put $p=1-f$ and $\beta=1-s$.
For a word in $\mathcal F$ and opposite-sign positions $i<j$, form
\[
 G=a_i v^{\epsilon_{i+1}}a_{i+1}\cdots
       v^{\epsilon_{j-1}}a_{j-1}.
 \tag{\thesection.38}\label{eq:q-38}
\]
If $j=i+1$, this means $G=a_i$.
When $G$ contains a variable letter, $G$ is an interval word in
$\mathcal F$, and so $|\tau(G)|\le F(s)$.
Otherwise $|\tau(G)|\le\kappa$, because the two endpoint signs are
opposite.
Furthermore $v^{\epsilon_i}Gv^{\epsilon_j}$ is the larger interval
word from position $i$ to $j$, with unit boundaries, and belongs to
$\mathcal F$.  The orientation matters: because
$\epsilon_j=-\epsilon_i$,
\[
 \tau(fG)=\tau(v^{\epsilon_i}Gv^{\epsilon_j}).
\]
Consequently
\[
 |\tau(pGp)|
 =|\tau(G)-\tau(fG)|\le\kappa+2F(s).
 \tag{\thesection.39}\label{eq:q-39}
\]

The finite coefficient list for quantization consists of all these
$pGp$, padded by zeros, and its normalized mean bound is
\[
 \mu=\frac{\kappa+2F(s)}{1-s}.
 \tag{\thesection.40}\label{eq:q-40}
\]
For the analytic list, formally expand each $W(v+w)$.
For each nonempty set of new-$w$ positions whose signs all agree,
write its trace as
\[
 \tau(w^\sigma B_1\cdots w^\sigma B_l),\qquad l\le n.
 \tag{\thesection.41}\label{eq:q-41}
\]
The $B_j\in pMp$ are the compressed gaps of old-$v$ letters and
coefficients.  The last $B_j$ includes the cyclic wrap through both
outer coefficients.  They are contractions, and there are only
finitely many such templates.  Their values depend on $v$, but their
number and degrees are bounded by the fixed word patterns.

Apply Lemma~\ref{lem:quantization-analytic} with local error $\eta$,
mean bound \eqref{eq:q-40}, and rank bound $S-s$.
It gives $q\le p$ of constant rank $0<\Delta<S-s$ and a unitary
$w\in qRq$.  Put $u=v+w$, a normal partial unitary of rank $s+\Delta$.

Fix a degree-$k$ word in $\mathcal F$ and expand at $v+w$.
Apart from its old term, there are at most $2^k-1$ terms.
If the new letters in a term all have one sign, \eqref{eq:q-41} bounds its trace
by $\eta\Delta$, including the case of a single new letter.
Otherwise choose two consecutive new letters with opposite signs.
Their intervening old gap is some $G$ in \eqref{eq:q-38}.
The term contains $w^{\epsilon_i}G w^{-\epsilon_i}$, whose trace norm
equals $\|qGq\|_{1,\tau}$; all surrounding factors are contractions.
It follows that its trace has absolute value at most
$(C\mu+\eta)\Delta$.
Thus for every word in $\mathcal F$,
\[
 |\tau(W(u))|
 \le F(s)+\left(
       \frac{D\kappa+2DF(s)}{1-s}+B_0\eta
      \right)\Delta.
 \tag{\thesection.42}\label{eq:q-42}
\]

The derivative of the budget is
$F'(s)=K(F(s)+\kappa+a)/(1-s)$, and hence
\[
 F'(s)-\frac{D\kappa+2DF(s)}{1-s}
 =\frac{(K-2D)F(s)+(K-D)\kappa+Ka}{1-s}
 \ge a.
 \tag{\thesection.43}\label{eq:q-43}
\]
Since $B_0\eta=a/2$, convexity of $F$ and \eqref{eq:q-42}--\eqref{eq:q-43} imply
$|\tau(W(u))|\le F(s+\Delta)$.
Thus $u$ is an admissible strict extension, contradicting maximality.
The maximal rank must be $S$.

Finally set $t=v+(1-f)\in U(R)$.
Telescoping over the at most $n$ variable positions, using
$t-v=t^*-v^*=1-f$, gives
\[
 |\tau(W(t)) -\tau(W(v))|\le n\tau(1-f)=n\gamma.
 \tag{\thesection.44}\label{eq:q-44}
\]
This is a telescoping estimate, not a $2^n$-term expansion.
By \eqref{eq:q-34}--\eqref{eq:q-35},
$F(S)\le2aL=\varepsilon/4$ and $n\gamma=\varepsilon/4$.
Thus $|\tau(W(t))|\le\varepsilon/2<\varepsilon$ uniformly in $\tau$,
as required.
\end{proof}

\section{Comparison storage and finite capture}
\label{sec:audited-nonunital-capture}

We write $\mathrm C^*(S)$ without adjoining a unit.  If $D$ is nonunital,
$\widetilde D$ denotes its minimal unitization.  The uniform tracial
completion is always the completion of bounded balls.  Its designated
traces are the extensions of the specified traces on $D$.  The
uniform tracial seminorm on $D$ need not be faithful.  All series,
operator norms, and functional calculi in this section are taken in
$D$ or its minimal unitization before applying the canonical map to
the completion.

\begin{proposition}[Comparison storage without a unit]
\label{prop:audited-nonunital-storage}
Let $D$ be a $\sigma$-unital $\mathrm C^*$-algebra with a nonempty family $X$ of
tracial states.  Suppose that
\[
 \mathrm d_\tau(a)<\mathrm d_\tau(b)\quad(\tau\in X)
 \quad\Longrightarrow\quad a\precsim b,
 \qquad a,b\in D_+.
\]
Let $a_i\in D_+$ have norm at most $2$, and let $k_i,\ell_i\in D_+$ be
contractions satisfying
\[
 \ell_i a_i=a_i,\qquad \ell_i k_i=\ell_i,\qquad
 \mathrm d_\tau(k_i)\le\alpha_i\quad(\tau\in X),\qquad
 \sum_{i\ge1}\alpha_i<\tfrac14.
\]
For prescribed $0<\rho_i\le1$, there are $H,K\in D_{\mathrm{sa}}$ and
$V_i\in \mathrm C^*(H,K)$ such that
\[
                    \|V_i^*V_i-a_i\|<\rho_i\qquad(i\ge1).
\]
In particular, the decoding algebra in this statement has no artificially
adjoined unit.
\end{proposition}

\begin{proof}
Fix a strictly positive contraction $h\in D$.  Every $\tau\in X$ extends
to a tracial state on $\widetilde D$.  In its nondegenerate tracial GNS
representation, $h$ has support $1$.  Consequently, for every positive
contraction $f\in\widetilde D$,
\begin{equation}\label{eq:audited-rank-fhf}
                         \mathrm d_\tau(fhf)=\mathrm d_\tau(f).
\end{equation}
Indeed, $\ker(h^{1/2}f)=\ker(f)$ in this representation; the positive
operators $fhf=(h^{1/2}f)^*(h^{1/2}f)$ and $f$ therefore have the same
support.  This argument uses the norm-one extension of $\tau$ to
$\widetilde D$, and does not introduce any other trace on the unitization.

Put $R=2$ and $\eta_i=\rho_i/[4(2R+1)]$.  Choose positive numbers
$w_i,\mu_i,\lambda_i,\delta_i$ recursively.  First choose $w_i\le2^{-i}$
such that, for $j<i$,
\begin{equation}\label{eq:audited-w-budget}
                 w_i<2^{-(i-j)}w_j\eta_j/(6R).
\end{equation}
Next choose $\mu_i\le2^{-i}$ such that, for $j<i$,
\begin{equation}\label{eq:audited-mu-budget}
                 2\mu_i<2^{-(i-j)}\lambda_j\eta_j/(6R).
\end{equation}
Choose $0<\lambda_i\le\mu_i$ sufficiently small that
\begin{equation}\label{eq:audited-lambda-budget}
 R\sqrt{3\lambda_i/\mu_i}
 +\frac{R}{w_i}\sum_{j<i}w_j
       \left(\sqrt{3\lambda_i/\lambda_j}
                    +\sqrt{3\lambda_i/\mu_j}\right)<\eta_i/3.
\end{equation}
Finally choose $0<\delta_i<\lambda_i\eta_i/(6R)$.  These choices depend
only on prescribed errors and the preceding numerical choices.

We construct contractions $e_i\in D_+$ and elements $r_i\in D$ with
\begin{equation}\label{eq:audited-witness-properties}
 e_i r_i=r_i,\quad r_i k_i=r_i,\quad \|r_i\|<R,\quad
 \mathrm d_\tau(e_i)\le\alpha_i,\quad
 \|r_i^*r_i-a_i\|<\rho_i/2.
\end{equation}
At stage $i$, set
\[
 S_i=\sum_{j<i}(\lambda_j e_j+\mu_jk_j)+\mu_i k_i,
 \qquad f_i=(1-S_i/\delta_i)_+\in\widetilde D,
 \qquad b_i=f_i h f_i\in D_+.
\]
Support subadditivity and the inductive bounds give
\[
 \mathrm d_\tau(S_i)\le\sum_{j<i}\alpha_j+\sum_{j\le i}\alpha_j<\tfrac12.
\]
Since $f_i=1$ on the kernel of $S_i$, \eqref{eq:audited-rank-fhf} gives
\begin{equation}\label{eq:audited-large-carrier}
 \mathrm d_\tau(b_i)=\mathrm d_\tau(f_i)\ge1-\mathrm d_\tau(S_i)>\tfrac12
       >\alpha_i\ge \mathrm d_\tau(a_i).
\end{equation}
The last inequality follows from the two mask identities in the
statement: $s(a_i)\le1_{\{1\}}(\ell_i)\le s(\ell_i)
\le1_{\{1\}}(k_i)\le s(k_i)$ in the bidual.
Comparison gives $a_i\precsim b_i$.  Choose $y_i\in D$ with
$\|y_i^*b_i y_i-a_i\|<\rho_i/4$ and set
$r_i^0=b_i^{1/2}y_i\ell_i$.  Then
\[
 r_i^0k_i=r_i^0,\qquad
 \|(r_i^0)^*r_i^0-a_i\|<\rho_i/4,\qquad \|r_i^0\|<R,
\]
and
\[
 \mathrm d_\tau(r_i^0(r_i^0)^*)
 =\mathrm d_\tau((r_i^0)^*r_i^0)\le \mathrm d_\tau(\ell_i)\le\alpha_i.
\]
Here comparison witnesses can be chosen in $D$ even when Cuntz
subequivalence was defined in the stabilization: compress a witness
to its $(1,1)$ matrix entry.

Put $c_i=r_i^0(r_i^0)^*$.  Choose continuous functions
$\chi_i,\gamma_i:[0,\infty)\to[0,1]$, zero near zero, with
$\gamma_i\chi_i=\chi_i$ and
\[
                 \sup_{t\ge0}t(1-\chi_i(t)^2)<\rho_i/4.
\]
Set $r_i=\chi_i(c_i)r_i^0$ and $e_i=\gamma_i(c_i)$.  Functional
calculus gives all of \eqref{eq:audited-witness-properties}.  For the
Gram estimate, use
\[
 (r_i^0)^*\chi_i(c_i)^2r_i^0
  =\bigl[t\chi_i(t)^2\bigr]\bigl((r_i^0)^*r_i^0\bigr).
\]
The range support of $b_i$ is contained in that of $f_i$, hence in
$1_{[0,\delta_i)}(S_i)$.  Therefore
\begin{equation}\label{eq:audited-S-small}
                              \|S_i r_i\|\le\delta_i R.
\end{equation}
This completes the induction.  Notice that replacing the unital
carrier $(1-S_i/\delta_i)_+$ by $f_i h f_i$ preserves both its rank
bound and its range control, which are precisely the two required
properties.

Define norm-convergent series
\[
 H=\sum_{i\ge1}(\lambda_i e_i+\mu_i k_i),\qquad
 K=\sum_{i\ge1}w_i(r_i+r_i^*).
\]
The numerical choices imply
\begin{equation}\label{eq:audited-tail-bounds}
 \sum_{j>i}w_j< w_i\eta_i/(6R),\qquad
 \sum_{j>i}(\lambda_j+\mu_j)<\lambda_i\eta_i/(6R).
\end{equation}
Choose a continuous $\chi:[0,\infty)\to[0,1]$ equal to one on $[0,2]$
and zero on $[3,\infty)$, and put
$F_i=\chi(H/\lambda_i)\in\widetilde D$.  The continuous function
\[
 \theta_i(t)=\frac{1-\chi(t/\lambda_i)}{t-\lambda_i},
\]
defined as zero near $\lambda_i$, has norm at most $1/\lambda_i$.
Since
\[
 (1-F_i)r_i=\theta_i(H)(H-\lambda_i)r_i,
\]
\eqref{eq:audited-S-small} and \eqref{eq:audited-tail-bounds} give
\begin{equation}\label{eq:audited-main-extraction}
 \|F_i r_i-r_i\|
 \le\frac{R}{\lambda_i}
       \left(\delta_i+\sum_{j>i}(\lambda_j+\mu_j)\right)<\eta_i/3.
\end{equation}
For the earlier witnesses, let $p_j,q_j$ be the left and right support
projections of $r_j$.  The identities $e_jr_j=r_j$ and $r_jk_j=r_j$
give $e_j\ge p_j$ and $k_j\ge q_j$.  Thus
$H\ge\lambda_jp_j$ and $H\ge\mu_jq_j$.  Since
$F_iHF_i\le3\lambda_i1$, we obtain
\begin{equation}\label{eq:audited-earlier-extraction}
 \|F_ir_j\|\le R\sqrt{3\lambda_i/\lambda_j},\qquad
 \|F_ir_j^*\|\le R\sqrt{3\lambda_i/\mu_j}.
\end{equation}
No commutation between $H$ and the witnesses or masks is used.

Set $V_i=w_i^{-1}F_iK$.  This element belongs to $\mathrm C^*(H,K)$, because
\[
                  F_iK=K+(F_i-1)K,
\]
and $F_i-1$ is a continuous function of $H$ vanishing at zero.
Splitting its defining series into the main term, its own adjoint,
earlier terms and later terms gives
\[
\begin{split}
 \|V_i-r_i\|\le{}&\|F_ir_i-r_i\|+R\sqrt{3\lambda_i/\mu_i}\\
 &+\frac{R}{w_i}\sum_{j<i}w_j
       \left(\sqrt{3\lambda_i/\lambda_j}
                    +\sqrt{3\lambda_i/\mu_j}\right)
       +\frac{2R}{w_i}\sum_{j>i}w_j<\eta_i.
\end{split}
\]
The three groups are bounded by $\eta_i/3$ using, respectively,
\eqref{eq:audited-main-extraction}, \eqref{eq:audited-lambda-budget},
and \eqref{eq:audited-tail-bounds}.  Finally
\[
 \|V_i^*V_i-r_i^*r_i\|
 \le(2R+\eta_i)\eta_i\le(2R+1)\eta_i=\rho_i/4.
\]
Together with \eqref{eq:audited-witness-properties}, this proves the
asserted error bound.
\end{proof}

\begin{proposition}[Finite capture by actual elements]
\label{prop:audited-finite-capture}
Let $D$ be $\sigma$-unital and let $X$ be a nonempty family of
tracial states such that $D$ has comparison by $X$ as in
Proposition~\ref{prop:audited-nonunital-storage}.
Let $\mathcal M$ be the bounded uniform tracial completion of $D$ for
$X$, and write $j_D:D\to\mathcal M$ for its canonical homomorphism.
Given $m\ge1$ and selfadjoint contractions $x_1,\ldots,x_m\in\mathcal M$,
there is a sub-$\mathrm C^*$-algebra $C_0\subseteq D$ generated by $m+2$
selfadjoint elements such that each $x_j$ is a uniform $2$-norm limit
of the images under $j_D$ of selfadjoint contractions in $C_0$.
The map $j_D$ need not be injective.
\end{proposition}

\begin{proof}
Fix summable positive sequences $\epsilon_n\le1$ and
$\rho_{j,n}^{\pm}\le1$, with
$\sum_n(\rho_{j,n}^++\rho_{j,n}^-)<\infty$ for every $j$.
Choose positive support budgets $\alpha_{j,n}^{\pm}$ whose total sum
is less than $1/4$.  Put
$\beta_n=\min_{1\le j\le m,\ \sigma\in\{+,-\}}\alpha_{j,n}^{\sigma}$.
Choose numbers $\gamma_n>0$ tending to zero such that
\[
 \gamma_n+\gamma_{n+1}<\tfrac12\epsilon_n\sqrt{\beta_n}.
\]
For example, it suffices to choose $\gamma_n$ smaller than $2^{-n}$
and than one eighth of each of
$\epsilon_n\sqrt{\beta_n}$ and, when $n\ge2$,
$\epsilon_{n-1}\sqrt{\beta_{n-1}}$.
Choose selfadjoint contractions $a_{j,n}\in D$ satisfying
$\|j_D(a_{j,n})-x_j\|_{2,X}<\gamma_n$, and put
\[
 d_{j,n}=a_{j,n+1}-a_{j,n},\qquad
 a_{j,n}^{\pm}=(\pm d_{j,n}-\epsilon_n)_+.
\]
These contractive approximants exist by bounded-ball density, or
by applying selfadjoint clipping to a representing sequence, even
when $j_D$ has a kernel.  We have
\begin{equation}\label{eq:audited-difference-budget}
 4\epsilon_n^{-2}\|d_{j,n}\|_{2,X}^2<\beta_n.
\end{equation}

Choose continuous $f_n,g_n:\mathbb R\to[0,1]$ such that $f_n=0$ below
$\epsilon_n/2$ and $f_n=1$ above $3\epsilon_n/4$, whereas $g_n=0$
below $3\epsilon_n/4$ and $g_n=1$ above $\epsilon_n$.  Put
$k_{j,n}^{\pm}=f_n(\pm d_{j,n})$ and
$\ell_{j,n}^{\pm}=g_n(\pm d_{j,n})$.  These elements belong to $D$;
both functions vanish at zero.  They have the required mask identities,
and spectral Markov inequality gives
\[
 \mathrm d_\tau(k_{j,n}^{\pm})
 \le4\epsilon_n^{-2}\|d_{j,n}\|_{2,\tau}^2
 <\alpha_{j,n}^{\pm}.
\]
Apply Proposition~\ref{prop:audited-nonunital-storage} to this single
countable family.  It gives $H,K\in D_{\mathrm{sa}}$ and positive
$E_{j,n}^{\pm}\in \mathrm C^*(H,K)$ with
$\|E_{j,n}^{\pm}-a_{j,n}^{\pm}\|<\rho_{j,n}^{\pm}$.

Define
\[
 r_{j,n}=d_{j,n}-a_{j,n}^{+}+a_{j,n}^{-},\qquad
 e_{j,n}=(a_{j,n}^{+}-E_{j,n}^{+})
                         -(a_{j,n}^{-}-E_{j,n}^{-}).
\]
Scalar functional calculus gives $\|r_{j,n}\|\le\epsilon_n$ and
$\|e_{j,n}\|\le\rho_{j,n}^{+}+\rho_{j,n}^{-}$.  Therefore
\[
 c_j=a_{j,1}+\sum_{n\ge1}(r_{j,n}+e_{j,n})\in D_{\mathrm{sa}}
\]
is a norm-convergent series.  Let $C_0=\mathrm C^*(c_1,\ldots,c_m,H,K)$ and
\[
 P_{j,N}=c_j+\sum_{n<N}(E_{j,n}^{+}-E_{j,n}^{-})\in C_0.
\]
The identity $r_{j,n}+e_{j,n}=d_{j,n}-E_{j,n}^{+}+E_{j,n}^{-}$ gives
the exact telescoping formula
\[
 P_{j,N}-a_{j,N}=\sum_{n\ge N}(r_{j,n}+e_{j,n}).
\]
Its right side tends to zero in the operator norm of $D$.
Consequently $j_D(P_{j,N})\to x_j$ uniformly in $2$-norm, and the
$P_{j,N}$ are uniformly norm bounded in $D$.
Clip them to $[-1,1]$.  Clipping does not increase the $2$-norm
distance to a selfadjoint contraction, separately at every trace;
the clipped elements remain in $C_0$ because the clipping function
vanishes at zero.  They give the asserted contractive approximants.
\end{proof}

\begin{remark}[Location of captured elements]
Proposition~\ref{prop:audited-finite-capture} applies to the bounded
uniform completion of the diagonal algebra.  It does not assert that
an arbitrary element of the uniform tracial ultrapower is a uniform
$2$-norm limit of diagonal elements.  Such an ultrapower element need
not admit the uniformly $2$-Cauchy sequence used in the proof.
\end{remark}

\section{Uniform relative generation from finite word estimates}
\label{sec:uniform-relative-generation}

The argument in this section converts a finite uniform word estimate
into relative generation in uniform tracial norm.  The fixed algebra
$C_0$ may be nonunital and need not contain the auxiliary AF algebra
used to obtain that estimate.  All unitary witnesses, however, must
belong to $1+C_0$.  We give the finite perturbation and the simultaneous
Baire argument in full.

\newtheorem{rghypothesis}{Hypothesis}[section]
\newtheorem{rglemma}[rghypothesis]{Lemma}
\newtheorem{rgproposition}[rghypothesis]{Proposition}
\newtheorem{rgtheorem}[rghypothesis]{Theorem}
\newtheorem{rgcorollary}[rghypothesis]{Corollary}
\theoremstyle{remark}
\newtheorem{rgremark}[rghypothesis]{Remark}
\theoremstyle{plain}

\subsection{Setting and the precise finite input}

Let $A$ be a separable, possibly nonunital $\mathrm C^*$-algebra, and let
$X=\mathrm T(A)$ be its nonempty compact space of normalized bounded tracial
states.  Compactness is an explicit assumption when $A$ is nonunital.
Let $C_0\subseteq A$ be a fixed sub-$\mathrm C^*$-algebra.  Write $\widetilde A$
for the minimal unitization of $A$, with $\widetilde A=A$ when $A$ is
unital, and put
\begin{align*}
 \widetilde C_0&=\mathrm C^*(1,C_0)\subseteq\widetilde A,\\
 \mathcal U_1(A)&=\{u\in\mathcal U(\widetilde A):u-1\in A\},\\
 \mathcal U_1(C_0)&=\{v\in\mathcal U(\widetilde C_0):v-1\in C_0\}.
\end{align*}
Thus $\mathcal U_1(A)=\mathcal U(A)$ when $A$ is unital.  Both groups
carry the operator norm topology.  The group $\mathcal U_1(A)$ is a
closed subset of $\widetilde A$, and hence is a complete separable
metric space.

Let $\mathcal M$ be a unital tracial completion equipped with a unital
$*$-homomorphism $j:\widetilde A\to\mathcal M$.  Every $\tau\in X$ is
understood to have a designated tracial state on $\mathcal M$ whose
pullback by $j$ is the normalized extension of $\tau$ to
$\widetilde A$.  We assume that, for each fixed
$x\in\mathcal M$, the function $\tau\mapsto\tau(x)$ on $X$ is continuous
and affine.  We write
\[
 \|x\|_{2,\tau}=\tau(x^*x)^{1/2},\qquad
 \|x\|_{2,X}=\sup_{\tau\in X}\|x\|_{2,\tau}.
\]
The map $j$ need not be injective.  For elements of $\widetilde A$,
these formulas denote the corresponding tracial seminorms; inside
products with elements of $\mathcal M$ we implicitly apply $j$.
All operator norms, algebra memberships, and functional calculi in
$A$ or $\widetilde A$ refer to the actual algebras, before applying $j$.
The compact trace space throughout is $X$, not the full trace space of
$\mathcal M$.  In particular, the scalar quotient trace on a
unitization is not added to $X$.

Here is the exact finite input needed in this section.  A word in a new
unitary letter $t$ is written in expanded form as
\begin{equation}\label{eq:rg-word}
 x_0t^{\epsilon_1}x_1\cdots t^{\epsilon_m}x_m,
 \qquad \epsilon_j\in\{-1,1\},\quad m\geq1.
\end{equation}
Only gaps between successive, opposite signs are required to have small
trace.  The two outside coefficients, and all gaps between equal signs,
are arbitrary contractions.

\begin{rghypothesis}[Finite weakened-word estimate]
\label{hyp:rg-finite-words}
For every pair of positive integers $q,n$ and every $\varepsilon>0$ there
is a number $\kappa(q,n,\varepsilon)>0$ with the following property.
Consider at most $q$ words of the form \eqref{eq:rg-word}, with at most
$n$ occurrences of $t^{\pm1}$ in each word and with all coefficients in
the unit ball of $\mathcal M$.  Suppose that every internal coefficient
$x_j$ lying between opposite signs satisfies
\[
 \sup_{\tau\in X}|\tau(x_j)|\leq\kappa(q,n,\varepsilon).
\]
Then some $t\in\mathcal U_1(C_0)$, evaluated in each word as $j(t)$,
makes the trace of every prescribed word have absolute value at most
$\varepsilon$, simultaneously for all $\tau\in X$.
The number $\kappa(q,n,\varepsilon)$ is independent of the actual
coefficients.
\end{rghypothesis}

If the finite-word lemma is initially formulated for a fixed list of
sign patterns, the displayed formulation follows by taking the minimum
of its moduli over the finitely many lists with these size bounds.
Cardinality dependence must be retained.  No linear dependence on
$\varepsilon$ is assumed.

On uniformly norm-bounded sets, multiplication and taking adjoints are
continuous for $\|\cdot\|_{2,X}$.  In particular, for contractions,
\begin{equation}\label{eq:rg-telescope}
 \|a_1\cdots a_m-b_1\cdots b_m\|_{2,X}
 \leq\sum_{j=1}^m\|a_j-b_j\|_{2,X}.
\end{equation}
Consequently, any finite collection of strict moment inequalities is
retained, with an arbitrarily small additional error, when its unitary
letters are replaced sufficiently closely in uniform $2$-norm.

\begin{rglemma}[Transferring external AF witnesses]
\label{lem:rg-external-transfer}
Let $B\subseteq\mathcal R\subseteq\mathcal M$ be unital
sub-$\mathrm C^*$-algebras with the same unit as $\mathcal M$, where $B$ is AF.
Suppose that the unit ball of $B$ is dense in the unit ball of
$\mathcal R$ for $\|\cdot\|_{2,X}$ and that
$B\subseteq\overline{j(C_0)}^{\|\cdot\|_{2,X}}$.
Then every unitary of $\mathcal R$ is a uniform $2$-norm limit of
images of elements of $\mathcal U_1(C_0)$.  Consequently, a finite
weakened-word estimate with witnesses in $\mathcal U(\mathcal R)$ implies
Hypothesis~\ref{hyp:rg-finite-words}, after decreasing its prescribed
output error before transferring the witness.
\end{rglemma}
\begin{proof}
Fix $v\in\mathcal U(\mathcal R)$.  Approximate $v$ in uniform $2$-norm
by an element of $B$, and then approximate that element in norm by an
element of a finite-dimensional unital subalgebra of $B$.  Invertibles
are dense in this finite-dimensional algebra, so we obtain an
invertible $b\in B$ arbitrarily close to $v$ in uniform $2$-norm.
Write $b=w|b|$, where $w\in\mathcal U(B)$.  For every $\tau\in X$,
\[
 \operatorname{Re}\tau(v^*w|b|)\leq\tau(|b|)
\]
and expansion of squared norms gives
$\|b-w\|_{2,\tau}\leq\|b-v\|_{2,\tau}$.  Thus
\[
 \|w-v\|_{2,X}\leq2\|b-v\|_{2,X}.
\]
This proves that $\mathcal U(B)$ is uniformly $2$-norm dense in
$\mathcal U(\mathcal R)$.

The same finite-dimensional approximation and polar decomposition
approximate a unitary of $B$ in norm by a finite-stage unitary.  Every
finite-dimensional unitary has the form $e^{ih}$ for a self-adjoint
element $h$ of that stage.  By the capture assumption, choose
self-adjoint $h_m\in C_0$ such that $j(h_m)$ converges to $h$ in
uniform $2$-norm.
These approximants may be chosen with $\|h_m\|\leq\|h\|$.
Indeed, if $R=\|h\|>0$, clip any self-adjoint approximant $x$ to
$q=\operatorname{clip}_R(x)\in C_0$.  With
$p_+=(x-R)_+$ and $p_-=(-x-R)_+$, traciality gives
\[
 \tau((x-q)(q-h))
 =\tau(p_+(R-h))+\tau(p_-(R+h))\geq0,
\]
so $\|q-h\|_{2,\tau}\leq\|x-h\|_{2,\tau}$.  The case $h=0$ is
immediate.  The Duhamel formula and traciality now give
\[
 \|j(e^{ih_m})-e^{ih}\|_{2,X}\leq\|j(h_m)-h\|_{2,X}.
\]
Hence $e^{ih_m}\in\mathcal U_1(C_0)$ proves the asserted density.
For a finite family of word tests, apply the estimate in
$\mathcal R$ with output error $\varepsilon/2$, then approximate its
witness closely enough that \eqref{eq:rg-telescope} changes each
trace by less than $\varepsilon/2$.  Its input modulus is unchanged.
\end{proof}

\subsection{A finite uniform perturbation}

\begin{rgproposition}\label{prop:rg-uniform-correlation}
Assume Hypothesis~\ref{hyp:rg-finite-words}.  Given $u,z\in\mathcal U_1(A)$,
$\xi\in\widetilde A$, and $0<\alpha<1/2$, there are $u_\ell\in\mathcal U_1(A)$ and
$W_\ell\in\mathcal U(\mathrm C^*(1,C_0,u_\ell))$, for $\ell\geq1$, such that
\begin{equation}\label{eq:rg-norm-perturbation}
 \|u_\ell-u\|\leq\frac{2\sqrt2\alpha}{\sqrt\ell}
\end{equation}
and
\begin{equation}\label{eq:rg-correlation}
 \left|\tau(\xi^*W_\ell)
       -\frac\alpha2\tau(\xi^*z)\right|
 \leq\frac{4\alpha^2}{1-2\alpha}\|\xi\|_{2,\tau}+r_\ell
 \qquad(\tau\in X),
\end{equation}
where $r_\ell\geq0$ and $r_\ell\to0$.
\end{rgproposition}

\begin{proof}
Fix $\ell$.  All accuracies used below are finite accuracies, chosen
before the unitary letters are constructed.

\emph{Formal labels.}
Let $G$ be the abstract free group on $s_1,\ldots,s_\ell$, and set
\[
 g_j=s_1\cdots s_j,\qquad p_{kj}=g_kg_j^{-1},\qquad
 P=\{e\}\cup\{g_j^{-1}:1\leq j\leq\ell\}
          \cup\{p_{kj}:1\leq k,j\leq\ell\}.
\]
The set $P$ consists of labels, with repetitions removed.  We use a
second free group $F_P$, whose free generators are $V_p$, $p\in P$.
The symbol $V_e$ is a free generator, not the identity of $F_P$.
Moreover, $V_{p^{-1}}$ and $V_p^{-1}$ are different symbols whenever
both labels occur; no relation between them is imposed.

The diagonal labels $p_{kk}$ are all $e$.  If $k\ne j$, the reduced
word for $p_{kj}$ has $k$ positive letters followed by $j$ negative
letters, with no cancellation at the change of sign.  Thus all
off-diagonal labels $p_{kj}$ are distinct and different from $e$.
The labels $g_j^{-1}$ are also pairwise distinct.  Every nonidentity
$p^{-1}q$, with $p,q\in P$, has reduced $G$-length at most $4\ell$.

Write $\varphi$ for the canonical trace on $C_r^*(F_P)$, and let
\[
 c=\frac\alpha{2\ell},\qquad
 A_k=c\sum_{j=1}^\ell(V_{p_{kj}}-V_{p_{kj}}^{-1}),\qquad
 T=c\sum_{k,j=1}^\ell(V_{p_{kj}}+V_{p_{kj}}^{-1}).
\]
If $C_d$ denotes total degree $d$ in $\prod_{k=1}^\ell\exp(A_k)$,
then the absolute values of its reduced group coefficients are
coefficientwise bounded by the coefficients of $T^d$.  To see this
before group reduction, a term in $C_d$ has a nondecreasing sequence
of row indices and coefficient $1/(m_1!\cdots m_\ell!)\leq1$.
Every such term occurs in $T^d$; the latter also contains all other
row-index sequences.  Taking absolute values and then reducing words
preserves this domination.

The label multiplicities just checked give
\[
 \|T\|_2^2
 =c^2\{2\ell^2+2\ell(\ell-1)\}
 =\alpha^2\left(1-\frac1{2\ell}\right)\leq\alpha^2.
\]
The length-one case of Haagerup's inequality
\cite[Lemma~1.4]{Haagerup1979} gives $\|T\|\leq2\alpha$.
Therefore
\begin{equation}\label{eq:rg-free-degree}
 \|C_d\|_{2,\varphi}\leq\|T^d\|_{2,\varphi}
 \leq(2\alpha)^d.
\end{equation}
Also, the sum of the absolute group coefficients of $C_d$ is at most
$(\alpha\ell)^d/d!$.

Choose $0<\delta\leq\ell^{-1}$, and then choose $N\geq2$ such that
\begin{equation}\label{eq:rg-tail-choice}
 \sum_{d>N}\frac{(\alpha\ell)^d}{d!}\leq\delta,
 \qquad L_N=\sum_{d=2}^N\frac{(\alpha\ell)^d}{d!}.
\end{equation}
Set
\[
 b_\ell=\frac{\sqrt2\alpha}{\sqrt\ell},\qquad
 c_\ell=2b_\ell.
\]
Choose $k\geq1$ so that
\begin{equation}\label{eq:rg-k-choice}
 c_\ell^2\,2^{-2k}\leq\tfrac12\ell^{-4},
 \qquad D=\max\{2N,2k\}.
\end{equation}
Finally choose a positive moment tolerance $\eta$, as small as needed
below.  In particular, we require
\begin{equation}\label{eq:rg-eta-choices}
 \sqrt\eta L_N\leq\delta,
 \qquad
 \eta\alpha^{2k}c_\ell^{-2k+2}\leq\tfrac12\ell^{-4}.
\end{equation}

\emph{The finite tests and their order of construction.}
We construct $t_0,\ldots,t_\ell\in\mathcal U_1(C_0)$, and now
evaluate the formal $G$-words by
\[
 s_1=t_0,\qquad s_{j+1}=ut_j\ (1\leq j<\ell),\qquad
 a=g_\ell u,\quad W=at_\ell,\quad y=zW^*,\quad Y_p=pyp^{-1}.
\]
Let $S=\{p_{kj}:1\leq k,j\leq\ell\}$ and put
$M_\xi=\max\{1,\|\xi\|\}$.  The final requirements are
\begin{align}
 |\tau(Q(Y))|&\leq\eta
 &&\text{for every nonidentity reduced $F_P$-word $Q$ of length $\leq D$},
 \label{eq:rg-Y-tests}\\
 |\tau(\xi^*W)|&\leq M_\xi\eta,
 \label{eq:rg-special0}\\
 |\tau(\xi^*Y_pW)|&\leq M_\xi\eta
 &&(p\in S\setminus\{e\}),\label{eq:rg-specialplus}\\
 |\tau(\xi^*Y_p^{-1}W)|&\leq M_\xi\eta
 &&(p\in S),\label{eq:rg-specialminus}
\end{align}
simultaneously for all $\tau\in X$.

Here are all the centering conditions needed for the final letter
$t_\ell$.  Group a reduced $F_P$-word into powers
$Y_{p_1}^{m_1}\cdots Y_{p_r}^{m_r}$ with successive labels different.
Since $y=zt_\ell^{-1}a^{-1}$, expansion into $t_\ell^{\pm1}$ has
opposite-sign gaps only of the forms
\begin{equation}\label{eq:rg-final-gaps}
 a^*p_i^{-1}p_{i+1}a
 \quad\text{or}\quad
 z^*p_i^{-1}p_{i+1}z.
\end{equation}
Their traces are those of nonidentity $G$-words of length at most
$4\ell$.  Outside coefficients are unrestricted in
Hypothesis~\ref{hyp:rg-finite-words}, so no cyclic end-gap condition
is needed.  The special tests have the expansions
\begin{align*}
 \xi^*W&=\xi^*a t_\ell,\\
 \xi^*Y_pW&=\xi^*pz\,t_\ell^{-1}
                 [a^*p^{-1}a]t_\ell,\\
 \xi^*Y_p^{-1}W&=\xi^*pa\,t_\ell
                 [z^*p^{-1}a]t_\ell.
\end{align*}
The first and third are analytic words.  In the second, $p\ne e$,
so its single opposite-sign gap again has the trace of a nonidentity
$G$-word of length at most $4\ell$.  For $p=e$ the omitted test is
instead the exact identity $Y_eW=z$.

We specify the earlier tolerances backwards, before choosing any
letters.  The final test count is at most
$\sum_{m=1}^D(2|P|)^m+2\ell^2+1$, and each final test has at most
$\max\{D,2\}$ occurrences of $t_\ell^{\pm1}$.
After dividing the special tests by $M_\xi$, all coefficients are
contractions.  Choose $\theta_\ell>0$ below the input-gap modulus
from Hypothesis~\ref{hyp:rg-finite-words} for this final list and
output error $\eta/2$.

For $r=\ell,\ldots,2$, consider all nonidentity reduced words of
length at most $L=4\ell$ on $s_1,\ldots,s_r$.  For the words involving
$s_r$, substitute $s_r=ut_{r-1}$ and expand its powers.  Every
opposite-sign gap is $c_0$ or $u^*c_0u$, where $c_0$ is a nonidentity
reduced word on the earlier generators of length at most $L$.
The number of tests, their lengths, and their coefficient norm bounds
are bounded solely in terms of $r,L$.  Let
$\kappa_r(\theta_r)>0$ be the corresponding structural input modulus,
and choose
\begin{equation}\label{eq:rg-backward}
 0<\theta_{r-1}\leq\min\{\theta_r,\kappa_r(\theta_r)\}.
\end{equation}
At the first stage, the tests are the nonzero powers of $s_1=t_0$
of length at most $L$; all are analytic and require no centered gaps.
Choose $t_0$ with these traces at most $\theta_1$.  Applying the
finite-word hypothesis successively constructs the remaining earlier
letters and ensures that all nonidentity words on $s_1,\ldots,s_r$
of length at most $L$ have trace at most $\theta_r$.  Words not
involving the new generator retain the preceding bound because
$\theta_{r-1}\leq\theta_r$.
The final application constructs $t_\ell$ and gives
\eqref{eq:rg-Y-tests}--\eqref{eq:rg-specialminus} with the smaller
tolerance $\eta/2$ in place of $\eta$.
This order of choices uses no coefficient-dependent modulus and no
accuracy-independent linear modulus.

All letters already lie in $\mathcal U_1(C_0)$, and all evaluated
unitary words lie in $\mathcal U_1(A)$.  In particular,
$Y_p-1\in A$ and $\operatorname{Im}(Y_p)\in A$ for every label $p$.

\emph{Remainder and first-order correlation.}
Define
\[
 h=\frac\alpha\ell\sum_{j=1}^\ell
          \operatorname{Im}(Y_{g_j^{-1}}),\qquad
 u'=e^{ih}u,\qquad H_k=g_khg_k^{-1}.
\]
Then $h=h^*\in A$, $\|h\|\leq\alpha$, and direct multiplication gives
\begin{equation}\label{eq:rg-product-factorization}
 W'=t_0\prod_{k=1}^\ell(u't_k)
    =\left(\prod_{k=1}^\ell e^{iH_k}\right)W,
 \qquad
 iH_k=c\sum_j(Y_{p_{kj}}-Y_{p_{kj}}^{-1}).
\end{equation}
All products are ordered by increasing indices from left to right.

For any group polynomial $Q$ of degree at most $N$, the moment tests
give
\begin{equation}\label{eq:rg-moment-transfer}
 \left|\tau(Q(Y)^*Q(Y))-\varphi(Q(V)^*Q(V))\right|
 \leq\eta\|Q\|_{\ell^1(F_P)}^2.
\end{equation}
Indeed, every nonidentity group word in $Q^*Q$ has length at most
$2N$; identity coefficients contribute no error.  Applying this to
$R_N=\sum_{d=2}^NC_d$, and using
\eqref{eq:rg-free-degree}--\eqref{eq:rg-eta-choices}, yields
\[
 \|R_N(Y)\|_{2,X}
 \leq\frac{4\alpha^2}{1-2\alpha}+\delta.
\]
In the actual algebra, each row in
\eqref{eq:rg-product-factorization} has norm at most $\alpha$.
Thus \eqref{eq:rg-tail-choice} bounds the operator norm of its
total-degree tail by $\delta$.  We obtain
\begin{equation}\label{eq:rg-uniform-remainder}
 \left\|\prod_ke^{iH_k}-1
       -c\sum_{k,j}(Y_{p_{kj}}-Y_{p_{kj}}^{-1})\right\|_{2,X}
 \leq\frac{4\alpha^2}{1-2\alpha}+2\delta.
\end{equation}

The label $e$ occurs $\ell$ times in the positive first-order sum,
and $Y_eW=z$.  All remaining first-order correlations are among
\eqref{eq:rg-specialplus}--\eqref{eq:rg-specialminus}.  Therefore,
for every $\tau\in X$,
\begin{equation}\label{eq:rg-before-clip}
 \left|\tau(\xi^*W')-\frac\alpha2\tau(\xi^*z)\right|
 \leq\left(\frac{4\alpha^2}{1-2\alpha}+2\delta\right)
           \|\xi\|_{2,\tau}
       +M_\xi\eta(1+\alpha\ell).
\end{equation}

\emph{Norm-small clipping.}
The free polynomial
\[
 h_{\mathrm f}=\frac\alpha{2i\ell}
       \sum_j(V_{g_j^{-1}}-V_{g_j^{-1}}^{-1})
\]
has $2\ell$ distinct length-one terms.  Haagerup's inequality gives
$\|h_{\mathrm f}\|\leq b_\ell$.  Its coefficient $\ell^1$-norm is
$\alpha$, so the moment tests through degree $2k$ imply
\[
 \tau(h^{2k})\leq\varphi(h_{\mathrm f}^{2k})+\eta\alpha^{2k}
 \leq b_\ell^{2k}+\eta\alpha^{2k}.
\]
Let $\widehat h$ be the continuous spectral clipping of $h$ to
$[-c_\ell,c_\ell]$.  The scalar estimate
\[
 |s-\operatorname{clip}_{c_\ell}(s)|^2
 \leq c_\ell^{-2k+2}|s|^{2k}\qquad(s\in\mathbb R)
\]
and \eqref{eq:rg-k-choice}--\eqref{eq:rg-eta-choices} give
\begin{equation}\label{eq:rg-clip-error}
 \|h-\widehat h\|_{2,X}^2
 \leq c_\ell^2\,2^{-2k}
       +\eta\alpha^{2k}c_\ell^{-2k+2}
 \leq\ell^{-4}.
\end{equation}
Put
\[
 u_\ell=e^{i\widehat h}u,\qquad
 W_\ell=t_0\prod_{k=1}^\ell(u_\ell t_k).
\]
The clipping function vanishes at zero, so $\widehat h\in A$.
These elements belong to $\mathcal U_1(A)$ and
$\mathcal U(\mathrm C^*(1,C_0,u_\ell))$, respectively.
The norm estimate \eqref{eq:rg-norm-perturbation} follows from
$\|\widehat h\|\leq c_\ell$ in $A$.
Since $h$ and $\widehat h$ commute, scalar functional calculus gives
$\|e^{ih}-e^{i\widehat h}\|_{2,X}\leq\ell^{-2}$.
Telescoping the $\ell$ occurrences of the perturbed unitary yields
\[
 \|W_\ell-W'\|_{2,X}\leq\ell^{-1}.
\]
In choosing $\eta$ before constructing the letters, impose also
$M_\xi\eta(1+\alpha\ell)\leq\ell^{-1}$.
Equation~\eqref{eq:rg-before-clip} now proves
\eqref{eq:rg-correlation}, for example with
$r_\ell=(3\|\xi\|+1)/\ell$, uniformly over $\tau$.
\end{proof}

\subsection{Hypographs and moving traces}

We give the topological step explicitly, because pointwise Baire
arguments in separate trace fibres would not suffice.

\begin{rglemma}\label{lem:rg-hypograph}
Let $S$ be a complete metric space, $Y$ a nonempty compact metric space, and
$f:S\times Y\longrightarrow[0,1]$ an upper semicontinuous function.
For $s\in S$ set
\[
 H(s)=\{(y,t)\in Y\times[0,1]:0\leq t\leq f(s,y)\}.
\]
There is a dense $G_\delta$ subset of $S$ at every point of which
$H$ is continuous in the Hausdorff metric.
\end{rglemma}
\begin{proof}
Each $H(s)$ is nonempty and compact.  Joint upper semicontinuity of
$f$ gives a closed graph for $H$.  Compactness of $Y\times[0,1]$
then gives upper semicontinuity of the set-valued map $H$.
For every fixed $w\in Y\times[0,1]$, the function
\[
 d_w(s)=\operatorname{dist}(w,H(s))
\]
is bounded and lower semicontinuous.  A bounded lower semicontinuous
real-valued function on a complete metric space has a dense
$G_\delta$ set of continuity points.  For completeness, in any
nonempty open set choose a point whose value is within $\varepsilon$
of the supremum on that set; lower semicontinuity supplies a smaller
open set on which the oscillation is at most $2\varepsilon$.
The Baire category theorem, applied for a sequence of $\varepsilon$
tending to zero, proves the assertion.

Take a countable dense set of $w$'s and intersect their continuity
sets.  All functions $w\mapsto\operatorname{dist}(w,H(s))$ are
$1$-Lipschitz.  A finite net in the compact target shows that,
at a common continuity point, their convergence is uniform in $w$.
For nonempty compact sets $K,L$ in a metric space,
\[
 d_{\mathrm H}(K,L)
 =\sup_w|\operatorname{dist}(w,K)-\operatorname{dist}(w,L)|.
\]
Thus $H$ is Hausdorff-continuous at every such point.
\end{proof}

\begin{rgtheorem}\label{thm:rg-uniform-relative-generation}
Assume Hypothesis~\ref{hyp:rg-finite-words}.  The set of
$u\in\mathcal U_1(A)$ for which the unit ball of $\mathrm C^*(1,C_0,u)$ is dense
in the unit ball of $\widetilde A$ for $\|\cdot\|_{2,X}$ contains a dense
$G_\delta$ subset of $\mathcal U_1(A)$, with its norm topology.
\end{rgtheorem}
\begin{proof}
For $z\in\mathcal U_1(A)$ define
\begin{equation}\label{eq:rg-defect}
 f_z(u,\tau)=\inf_{b\in \mathrm C^*(1,C_0,u)}\|z-b\|_{2,\tau}^2.
\end{equation}
This function takes values in $[0,1]$.  Expressing the infimum using
all fixed finite $*$-polynomials in $u,u^*$ with coefficients in $\widetilde C_0$
shows that $f_z$ is jointly upper semicontinuous on
$\mathcal U_1(A)\times X$, where the first space has its norm topology
and the second its weak-$*$ topology.  For each fixed polynomial,
the squared error is jointly continuous, and an infimum of continuous
functions is upper semicontinuous.

Both $\mathcal U_1(A)$ and $X$ are metrizable; the former is complete
for its norm metric and the latter is compact.  Apply
Lemma~\ref{lem:rg-hypograph} to $f_z$.  Intersect its dense
$G_\delta$ sets for a countable norm-dense family
$\{z_j\}\subset\mathcal U_1(A)$.  Fix $u$ in the resulting dense
$G_\delta$ set.

We claim that $f_z(u,\tau)=0$ for every test $z=z_j$ and every
$\tau\in X$.  Suppose instead that
$f_z(u,\tau_0)=d^2>0$.  Choose
\[
 0<\alpha<\min\{1/8,d/128\},\qquad
 0<\varepsilon<\alpha d^2/64.
\]
There is a fixed polynomial $b=b(u)$ with coefficients in $\widetilde C_0$ such
that
\begin{equation}\label{eq:rg-near-projection}
 \|z-b\|_{2,\tau_0}^2\leq d^2+\varepsilon^2.
\end{equation}
Put $\xi=z-b$.  In the GNS von Neumann algebra of $\tau_0$, let
$N_0=\pi_{\tau_0}(\mathrm C^*(1,C_0,u))''$ and let $\mathbb E_{N_0}$ be the
trace-preserving conditional expectation.  The distance in
\eqref{eq:rg-defect} equals the distance to $L^2(N_0,\tau_0)$.
Pythagoras therefore gives
\begin{equation}\label{eq:rg-xi-data}
 \|\xi\|_{2,\tau_0}\leq\sqrt{d^2+\varepsilon^2},\qquad
 |\tau_0(\xi^*z)-d^2|\leq\varepsilon.
\end{equation}
For the second inequality, write $b=\mathbb E_{N_0}z+e$, where
$\|e\|_{2,\tau_0}\leq\varepsilon$, and use
$\|\mathbb E_{N_0}z\|_{2,\tau_0}\leq1$.

Apply Proposition~\ref{prop:rg-uniform-correlation} to these fixed
$u,z,\xi,\alpha$.  Since $u_\ell\to u$ in norm and $u$ is a
Hausdorff-continuity point of the hypograph, the point
$(\tau_0,d^2)$ in that hypograph has approximants
$(\tau_\ell,t_\ell)$ in the hypographs at $u_\ell$ satisfying
\begin{equation}\label{eq:rg-moving-trace}
 \tau_\ell\longrightarrow\tau_0,
 \qquad f_z(u_\ell,\tau_\ell)\geq t_\ell\longrightarrow d^2.
\end{equation}
Joint upper semicontinuity also gives
$\limsup_\ell f_z(u_\ell,\tau_\ell)\le d^2$, so these defects converge
to $d^2$.
Define $b_\ell=b(u_\ell)$ using the same polynomial.  Then
$b_\ell\to b$ in norm.  In the GNS algebra of $\tau_\ell$, put
$N_\ell=\pi_{\tau_\ell}(\mathrm C^*(1,C_0,u_\ell))''$.  Pythagoras and
\eqref{eq:rg-near-projection}--\eqref{eq:rg-moving-trace} give
\begin{align*}
 \|\mathbb E_{N_\ell}z-b_\ell\|_{2,\tau_\ell}^2
 &=\|z-b_\ell\|_{2,\tau_\ell}^2-f_z(u_\ell,\tau_\ell)\\
 &\leq\varepsilon^2+o(1).
\end{align*}
The first term on the right converges because $b_\ell\to b$ in
norm and $\tau_\ell\to\tau_0$ weak-$*$.
Furthermore,
\[
 \mathbb E_{N_\ell}\xi
 =(\mathbb E_{N_\ell}z-b_\ell)+\mathbb E_{N_\ell}(b_\ell-b),
\]
so contractivity gives
$\|\mathbb E_{N_\ell}\xi\|_{2,\tau_\ell}\leq\varepsilon+o(1)$.
Since $W_\ell\in N_\ell$ is unitary,
\begin{equation}\label{eq:rg-small-correlation}
 |\tau_\ell(\xi^*W_\ell)|\leq\varepsilon+o(1).
\end{equation}

On the other hand, the \emph{uniform} estimate
\eqref{eq:rg-correlation} can be evaluated at these moving traces.
Using \eqref{eq:rg-xi-data} gives
\begin{align*}
 \liminf_\ell|\tau_\ell(\xi^*W_\ell)|
 &\geq\frac\alpha2(d^2-\varepsilon)
       -\frac{4\alpha^2}{1-2\alpha}\sqrt{d^2+\varepsilon^2}\\
 &>\frac{\alpha d^2}{8}>\varepsilon.
\end{align*}
For the displayed numerical bound, use
$4\alpha^2/(1-2\alpha)\leq8\alpha^2$,
$\sqrt{d^2+\varepsilon^2}\leq2d$, and the chosen bounds on
$\alpha,\varepsilon$.  This contradicts
\eqref{eq:rg-small-correlation}.  Thus the claim holds.

Norm density of the test unitaries extends this conclusion to every
element of $\mathcal U_1(A)$.  Every unitary of $\widetilde A$ is a
scalar multiple of one in $\mathcal U_1(A)$: in the nonunital case,
divide by its scalar image in $\widetilde A/A$.  Since every element
of a unital $\mathrm C^*$-algebra is a linear combination of unitaries, for
every $a\in\widetilde A$ we have
\begin{equation}\label{eq:rg-pointwise-density}
 \inf_{b\in \mathrm C^*(1,C_0,u)}\|a-b\|_{2,\tau}^2=0
 \qquad(\tau\in X).
\end{equation}
Let $C=\mathrm C^*(1,C_0,u)$ and regard $C$ as a real locally convex vector
space with its norm topology.  The function
$(\tau,b)\mapsto\|a-b\|_{2,\tau}^2$ is continuous and affine in
$\tau$ and continuous and convex in $b$.  Sion's minimax theorem
\cite[Theorem~4.2]{Sion1958}, with the compact convex variable $X$, therefore gives
\begin{equation}\label{eq:rg-sion}
 \inf_{b\in C}\sup_{\tau\in X}\|a-b\|_{2,\tau}^2
 =\sup_{\tau\in X}\inf_{b\in C}\|a-b\|_{2,\tau}^2=0.
\end{equation}
No compactness of the approximating space $C$ is required.

For clarity, bounded approximation follows from
\eqref{eq:rg-sion}; it is not an additional assumption.  If
$\|a\|\leq1$ and $b\in C$, set
$\widehat b=b f(|b|)$, where $f(t)=1$ for $0\leq t\leq1$ and
$f(t)=1/t$ for $t\geq1$.  Then $\widehat b\in C$ is a contraction.
For each trace, write $b=v|b|$ in the corresponding finite von
Neumann algebra and put $P=(|b|-1)_+$.  Since
$b-\widehat b=vP$,
\[
 \operatorname{Re}\tau((b-\widehat b)^*(\widehat b-a))
 =\tau(P)-\operatorname{Re}\tau(Pv^*a)\geq0.
\]
Expanding $\|b-a\|_{2,\tau}^2$ now gives
$\|\widehat b-a\|_{2,\tau}\leq\|b-a\|_{2,\tau}$.
Thus the contraction ball of $C$ is dense in that of $\widetilde A$ in the
uniform $2$-norm, as asserted.
\end{proof}

\begin{rgcorollary}[Removing the adjoined unit]
\label{cor:rg-remove-unit}
Assume Hypothesis~\ref{hyp:rg-finite-words}, and assume that the unit of
$\mathcal M$ belongs to the uniform $2$-norm closure of $j(C_0)$.
There is $h=h^*\in A$ such that the unit ball of the possibly nonunital
algebra $C=\mathrm C^*(C_0,h)\subseteq A$ is dense in the unit ball of $A$ for
$\|\cdot\|_{2,X}$.  If $C_0$ is generated by $m$ self-adjoint elements,
then $C$ is generated by $m+1$ self-adjoint elements.
\end{rgcorollary}
\begin{proof}
Choose a unitary $u$ from the dense set in
Theorem~\ref{thm:rg-uniform-relative-generation} such that
$\|u-1\|<1$.  Continuous functional calculus with the principal
logarithm gives $h=-i\log u=h^*\in A$, because the logarithm vanishes
at $1$.  Moreover,
\[
 \mathrm C^*(1,C_0,u)=\mathrm C^*(1,C_0,h)=\mathrm C^*(1,C).
\]
The uniform $2$-norm closure of the linear space $j(C)$ contains
$1_{\mathcal M}$ by the hypothesis on $C_0$.  It therefore contains
$j(C+\mathbb C1)$, and hence $j(\mathrm C^*(1,C))$, since $j$ is
contractive for operator norms.  The theorem now implies that every
$a\in A$ is a uniform $2$-norm limit of elements of $C$.  For a
contraction $a$, clip each approximant $b\in C$ to
$b f(|b|)$ using the function in the proof of the theorem.  Although
$f(0)=1$, the product belongs to $C$, because
$f(|b|)\in C+\mathbb C1$ and $b\in C$.  The same trace inequality shows
that clipping decreases the error.  This gives the stated density of
unit balls.  The generator count is immediate from the definition of
$C$.
\end{proof}

\begin{rgremark}[Dependency boundary]\label{rem:rg-dependencies}
The only structural input in this section is the finite-word estimate
with actual witnesses in $\mathcal U_1(C_0)$.  An external AF algebra
inside the uniform tracial closure of $j(C_0)$ can supply those witnesses
by the transfer argument above.  If that AF algebra has unit $1$,
the extra unit-capture hypothesis in
Corollary~\ref{cor:rg-remove-unit} holds automatically.  The proof
uses neither the topology of $\partial\mathrm T(A)$ nor nuclearity.  The
passages from this tracially dense algebra to norm generation and
then to single generation require separate compression and encoding
arguments.
\end{rgremark}

\section{A uniform finite generator bound}
\label{sec:seven-generators}

\begin{proposition}[Seven selfadjoint generators from an external AF core]
\label{prop:audited-seven-generators}
Let $D$ be separable and $\sigma$-unital with compact nonempty $\mathrm T(D)$,
and assume comparison by these tracial states for positive elements
of $D$.  Suppose its bounded uniform tracial completion $\mathcal M_D$
is unital and has a unital AF subalgebra $B$ with the following two
properties:
\begin{enumerate}
\item $B$ has two selfadjoint generators, without adjoining its unit.
\item The conclusion of Proposition~\ref{prop:finite-weakened-word} holds with coefficients in $\mathcal M_D$, witnesses in the uniform $2$-closure of $B$, and designated traces $\mathrm T(D)$, with moduli independent of coefficient values.
\end{enumerate}
Then $D$ has at most seven selfadjoint norm generators.
\end{proposition}

\begin{proof}
The canonical map $D\to\mathcal M_D$ need not be injective.
Whenever a closure in $\mathcal M_D$ is mentioned below, it is the
closure of the canonical image.  All generators and comparison
witnesses are actual elements of $D$.

Scale the two selfadjoint AF generators to contractions.  By
Proposition~\ref{prop:audited-finite-capture}, a subalgebra
$C_0\subseteq D$ generated by four selfadjoint elements has bounded
uniform $2$-closure containing $B$.  It also contains the external
unit in its bounded uniform $2$-closure, because the AF generators
generate that unit without adjoining it.  Concretely, evaluate
nonunital $*$-polynomials in bounded approximants and then apply
the radial clipping $b\mapsto b f(|b|)$ from the proof of
Corollary~\ref{cor:rg-remove-unit}.  This preserves $C_0$ and
does not increase the uniform $2$-distance to a contraction.

By Lemma~\ref{lem:rg-external-transfer}, finite-word witnesses in
the completion of $B$ transfer to exact unitaries in $1+C_0$.
Thus Hypothesis~\ref{hyp:rg-finite-words} holds for this fixed
$C_0\subseteq D$.

Apply Theorem~\ref{thm:rg-uniform-relative-generation} to the complete norm-metric group
$\mathcal U(1+D)$.  The perturbation stays in this group: its
selfadjoint perturbing term is a real linear combination of $(Y-Y^*)/i$ with $Y\in1+D$,
so lies in $D$, and the clipping function vanishes at zero.
The good unitaries are norm dense; choose one, $u$, with
$\|u-1\|<1$.  Its principal logarithm $z=-i\log u$ belongs to
$D_{\mathrm{sa}}$, and
$\mathrm C^*(1,u)=\mathrm C^*(1,z)$ in $\widetilde D$.

The relative-generation conclusion gives uniform $2$-density of
$\mathrm C^*(1,C_0,u)$ in $\mathcal M_D$.  Its adjoined unit is itself a
bounded uniform $2$-limit from $C_0$.  Hence, by
Corollary~\ref{cor:rg-remove-unit}, $C=\mathrm C^*(C_0,z)\subseteq D$
has selfadjoint contraction ball uniformly $2$-dense in that of
$D$.  It has five selfadjoint generators.  Clipping here does not
adjoin a unit, because clipping to $[-1,1]$ vanishes at zero.

Choose a norm-dense sequence $(a_j)$ in the selfadjoint unit ball of
$D$, and put $\epsilon_k=2^{-k}$.  Give each sign of each pair $(j,k)$
the budget
\[
 \alpha_{j,k}^{\pm}=2^{-(j+k+4)},\qquad
       \sum_{j,k\ge1,\,\sigma\in\{+,-\}}\alpha_{j,k}^\sigma=\tfrac18.
\]
By the selfadjoint contraction-ball density just proved, choose
$c_{j,k}\in C_{\mathrm{sa}}$ with $\|c_{j,k}\|\le1$ so that, for
$x_{j,k}=a_j-c_{j,k}$,
\[
 4\epsilon_k^{-2}\|x_{j,k}\|_{2,\mathrm T(D)}^2
                 <\min\{\alpha_{j,k}^+,\alpha_{j,k}^-\}.
\]
The positive cuts $a_{j,k}^{\pm}=(\pm x_{j,k}-\epsilon_k)_+$
have norm at most $2$.  Apply to $\pm x_{j,k}$ the same mask functions
as in Proposition~\ref{prop:audited-finite-capture}: the outer mask is
zero on $(-\infty,\epsilon_k/2]$ and one on
$[3\epsilon_k/4,\infty)$, and the inner mask is zero on
$(-\infty,3\epsilon_k/4]$ and one on $[\epsilon_k,\infty)$.
These masks lie in $D$, satisfy the two exact mask identities, and
their outer support ranks are less than $\alpha_{j,k}^{\pm}$ by
spectral Markov inequality.  Thus a single application of
Proposition~\ref{prop:audited-nonunital-storage} gives two selfadjoint
elements $H,K\in D$ and positive elements
$E_{j,k}^{\pm}\in\mathrm C^*(H,K)$ with
$\|E_{j,k}^{\pm}-a_{j,k}^{\pm}\|<\epsilon_k$.
The scalar identity
\[
 \|x-(x-\epsilon_k)_++(-x-\epsilon_k)_+\|\le\epsilon_k
 \qquad(x=x^*)
\]
shows, in the actual operator norm of $D$, that
\[
 \|a_j-(c_{j,k}+E_{j,k}^+-E_{j,k}^-)\|<3\epsilon_k.
\]
Letting $k\to\infty$ gives $D=\mathrm C^*(C,H,K)$, which proves
the asserted seven-selfadjoint bound.  In particular, this last
step also encodes elements annihilated by every designated trace;
no injectivity of $D\to\mathcal M_D$ is required.
\end{proof}

\section{Compression using a matrix cone and a reserved matrix position}
\label{sec:cone-compression}

We use $\operatorname{Her}_A(a)=\overline{aAa}$ for $a\in A_+$.
All $\mathrm C^*$-algebras generated in this section are generated without
automatically adjoining a unit.  A positive invertible element will
recover the ambient unit explicitly.

\begin{lemma}[Full hereditary bimodules]
\label{lem:full-hereditary-bimodule}
Let $A$ be a $\mathrm C^*$-algebra, let $B,E\subseteq A$ be nonzero hereditary
sub-$\mathrm C^*$-algebras, and put
\[
 X=\overline{\operatorname{span}}\,BAE.
\]
If $v\in X$ and $v^*v$ is full in $E$, then
\[
 X=\overline{\operatorname{span}}\,BvE.
\]
Consequently, if $v\in A$ satisfies $vv^*\in B$ and $v^*v$ is
full in $E$, then $\mathrm C^*(B,E,v)$ contains
$\overline{\operatorname{span}}BAE$ and
$\overline{\operatorname{span}}EAB$.
\end{lemma}

\begin{proof}
The closed linear span
\[
 J=\overline{\operatorname{span}}\,Ev^*BvE
\]
is a closed two-sided ideal of $E$.  An approximate unit
$(b_\lambda)$ of $B$ satisfies $b_\lambda v\to v$ in norm.
Multiplying $v^*b_\lambda v$ on both sides by an approximate unit
of $E$ shows that $v^*b_\lambda v\in J$, and hence $v^*v\in J$.
Fullness of $v^*v$ in $E$ gives $J=E$.
For $x\in X$, one has $xEv^*B\subseteq B$: this follows
by approximating $x$ and $v$ by finite sums in $BAE$ and using
$BAB\subseteq B$.  Hence
\[
 x(Ev^*BvE)\subseteq BvE.
\]
The right action of $E$ on $X$ is nondegenerate.  Approximating a
right approximate unit from $J=E$ therefore gives
$x\in\overline{\operatorname{span}}BvE$.  The reverse inclusion is
immediate.  For the last assertion, approximate units of $B$ and
$E$ show that $v\in\overline{BAE}$ whenever $vv^*\in B$ and
$v^*v\in E$; then apply the first assertion and take adjoints.
\end{proof}

\begin{proposition}[Compression by a matrix cone]
\label{prop:matrix-cone-compression}
Let $A$ be a unital $\mathrm C^*$-algebra with
$\mathrm T(A)\ne\varnothing$.  Assume comparison by ordinary
tracial states for positive elements of $A$:
\[
 \mathrm d_\tau(a)<\mathrm d_\tau(b)\quad(\tau\in\mathrm T(A))
 \quad\Longrightarrow\quad a\precsim b.
\]
Let $N\ge1$ and $n\ge3$ satisfy $n-1\ge N$.  Suppose there is a
completely positive contractive order-zero map
\[
 \phi\colon\mathrm M_n\longrightarrow A,
 \qquad t=\phi(1_n),
\]
such that, for some $0<\varepsilon<1/(64n)$,
\[
 \tau(t)>1-\varepsilon
 \qquad(\tau\in\mathrm T(A)).
\]
Put $\delta=1/16$ and suppose that
\[
 D=\operatorname{Her}_A\bigl((\phi(e_{11})-\delta)_+\bigr)
\]
is generated, without adjoining a unit, by at most $N$
selfadjoint elements.  Then $A$ is singly generated.
\end{proposition}

\begin{proof}
Write $H_t=\operatorname{Her}_A(t)$.  The supporting
$*$-homomorphism of $\phi$ from
\cite[Theorem~2.3]{WinterZacharias2009} gives matrix units, denoted again by
$e_{ij}$, in the multiplier algebra of $H_t$, with
\[
 \phi(e_{ij})=e_{ij}t=te_{ij}.
\]
Their sum $p=\sum_i e_{ii}$ is the support of $t$ in $A^{**}$.
For completeness, they act as multipliers on $H_t$: on the dense subspace $tAt$, multiplication on either side stays in $H_t$, since $e_{ij}t=\phi(e_{ij})\in H_t$.  These matrix units need not belong to $A$.  Whenever
$f\in \mathrm C([0,1])$ satisfies $f(0)=0$, however,
\begin{equation}\label{eq:cone-membership}
 e_{ij}f(t)\in H_t\subseteq A.
\end{equation}
For example, this follows first when $f(s)=s g(s)$, and then by
uniform approximation for arbitrary continuous $f$ vanishing at
zero.  The matrix units also act as multipliers on
$H_\delta=\operatorname{Her}_A((t-\delta)_+)$.  Thus
\[
 H_\delta\cong\mathrm M_n(D),
 \qquad
 D=e_{11}H_\delta e_{11},
\]
with the isomorphism given by
$(x_{ij})\mapsto\sum_{ij}e_{i1}x_{ij}e_{1j}$.
Every such expression belongs to $A$, since its entries belong to
$H_\delta$ and the $e_{ij}$ are multipliers of that algebra.

Put $r_0=n-1$, $\gamma=1/4$, and
\[
 d=(t-\delta)_+,
 \qquad d_i=e_{ii}d,
 \qquad w=e_{1n}(t-\gamma)_+.
\]
\[
 B_0=\sum_{i,j<n}e_{i1}De_{1j}\cong\mathrm M_{r_0}(D).
\]
Here $B_0$ is a hereditary subalgebra of $A$, and the displayed
finite matrix sum is closed.  Let
\[
 E_0=\operatorname{Her}_A\bigl(e_{nn}(t-\gamma)_+\bigr).
\]
We claim that
\begin{equation}\label{eq:good-carrier}
 B:=\mathrm C^*(B_0,w)
 =\operatorname{Her}_A\left(
    \sum_{i<n}d_i+e_{nn}(t-\gamma)_+\right).
\end{equation}
Indeed, $ww^*=e_{11}(t-\gamma)_+^2\in B_0$, and
\[
 \overline{w^*B_0w}
 =\operatorname{Her}_A(w^*w)=E_0.
\]
The first equality follows because $B_0$ is hereditary and
contains the hereditary algebra generated by $ww^*$.
The standard identity
\[
 \overline{w^*\operatorname{Her}(ww^*)w}
 =\operatorname{Her}(w^*w)
\]
follows directly from functional calculus.
Thus $E_0\subseteq \mathrm C^*(B_0,w)$.
Both $B_0$ and $E_0$ are nonzero: $\tau(t)>1-\varepsilon>\gamma$
forces $(t-\gamma)_+\ne0$, and the supporting matrix units make
all its diagonal slices simultaneously nonzero.  Moreover,
$w^*w$ is strictly positive in $E_0$.  Lemma
\ref{lem:full-hereditary-bimodule} recovers their two cross
bimodules from $w$.  Their support projections are orthogonal,
and the hereditary algebra on the sum of those supports is the
closed span of the two diagonal algebras and the two cross
bimodules.  This proves \eqref{eq:good-carrier}.

Choose continuous functions $g_0,g_1\colon[0,1]\to[0,1]$ such that
\[
 \begin{array}{c|cc}
       &\text{zero on}&\text{one on}\\ \hline
 g_0 &[0,\delta]&[2\delta,1]\\
 g_1 &[0,\gamma]&[2\gamma,1].
 \end{array}
\]
Set
\[
 s=\sum_{i<n}e_{ii}g_0(t)+e_{nn}g_1(t)\in B,
 \qquad r=1-s,
\]
\[
 b=e_{nn}(t-3/4)_+\in B.
\]
All the terms belong to $A$ by \eqref{eq:cone-membership};
membership in $B$ follows from \eqref{eq:good-carrier}.
The element $s$ is a positive contraction.  In every tracial
GNS closure, $r$ vanishes on the spectral subspace where
$t\ge1/2$.  Hence
\begin{equation}\label{eq:cone-residual-rank}
 \mathrm d_\tau(r)
 \le\tau(1_{[0,1/2)}(t))
 \le2\tau(1-t)<2\varepsilon.
\end{equation}
Traciality and the supporting matrix units give
\[
 n \mathrm d_\tau(b)
 =\tau(1_{(3/4,1]}(t))
 \ge\tau(t)-\tfrac34>\tfrac14-\varepsilon.
\]
One can justify the first equality using continuous cutoff
functions vanishing at zero and then monotone convergence in the
finite GNS algebra.  Our choice of $\varepsilon$ implies
\[
 2\varepsilon<\frac{1/4-\varepsilon}{n}.
\]
Comparison therefore gives $r\precsim b$.

Use the elementary exact-cutdown consequence of Cuntz
subequivalence: if $a\precsim b$ and $\eta>0$, then there is
$v\in A$ with $v^*v=(a-\eta)_+$ and
$vv^*\in\operatorname{Her}_A(b)$.  Indeed, choose $x\in A$ with
$\|a-x^*bx\|<\eta$.  The norm-perturbation cutdown factorization gives
$y\in A$ such that $(a-\eta)_+=y^*x^*bxy$; put
$v=b^{1/2}xy$.  Although the intermediate factors need not have a
uniform norm bound, this final witness satisfies
$\|v\|^2=\|(a-\eta)_+\|$.
With $\eta=1/4$ we obtain
\begin{equation}\label{eq:cone-link}
 v^*v=q:=(r-\tfrac14)_+,
 \qquad vv^*\in\operatorname{Her}_A(b)\subseteq B.
\end{equation}
Notice that $br=0$, hence $bq=0$.  The initial and final support
projections of $v$ are orthogonal, so $v^2=0$.
Also
\begin{equation}\label{eq:cone-invertible-cover}
 s+q=s+(1-s-\tfrac14)_+\ge\tfrac34\,1_A.
\end{equation}

Choose selfadjoint contraction generators $x_1,\ldots,x_N$ of $D$,
rescaling and padding with zero elements if necessary.  Choose distinct
numbers $\lambda_1,\ldots,\lambda_n$ in $[1,2]$ and a continuous function
$F\colon[0,1]\to[0,1]$ which is zero near zero and one on
$[\delta,1]$.  Define
\[
 L=\sum_{i=1}^n\lambda_i e_{ii}F(t),
 \qquad h=1_A+L,
\]
and
\begin{align*}
 k={}&\sum_{i=1}^N e_{i1}x_i e_{1i}
       +\sum_{i=1}^{n-2}
          \bigl(e_{i,i+1}d+e_{i+1,i}d\bigr)\\
    &+(w+w^*)+(v+v^*).
\end{align*}
These are selfadjoint elements of $A$.  In particular,
$e_{ij}d\in A$ by \eqref{eq:cone-membership}, and the payload
terms belong to $H_\delta$ by the matrix identification above.
Put $C=\mathrm C^*(h,k)$.  Since $h\ge1_A$, continuous functional calculus
with a function vanishing at zero and equal to one on
$[1,\|h\|]$ gives $1_A\in C$.
The construction also has bounds independent of $n$ and $N$:
\[
 1_A\le h\le3\,1_A,\qquad
 \|k\|\le1+2\|d\|+\|w\|+\|v\|
 \le1+\tfrac{15}{8}+\tfrac34+\sqrt{\tfrac34}<5.
\]
Here the diagonal payload has norm at most one, the consecutive
matrix edges have norm at most $2\|d\|$, and the orthogonal initial
and final supports give
$\|w+w^*\|=\|w\|$ and $\|v+v^*\|=\|v\|$.

Let $\mu_i=1+\lambda_i$, and choose real Lagrange polynomials
$\ell_i$ satisfying $\ell_i(\mu_j)=\delta_{ij}$.
Although the spectrum of $h$ may contain intervals between these
labels, its action on the hereditary carrier has exact labels:
for $x\in e_{i1}De_{1j}$,
\begin{equation}\label{eq:exact-plateau-labels}
 hx=\mu_i x,
 \qquad xh=\mu_j x.
\end{equation}
The same identities hold on the relevant corners involving
$w$, since its support lies where $t>\gamma>\delta$.
Moreover, the final range condition in \eqref{eq:cone-link}
implies
\[
 hv=\mu_n v,
 \qquad v^*h=\mu_n v^*.
\]
For example, approximate $v$ on the left by an approximate unit
of $\operatorname{Her}_A(b)$ and use that $h$ has the constant
label $\mu_n$ on this hereditary algebra.

It follows that for $i,j<n$ the cuts
$\ell_i(h)k\ell_j(h)$ kill both $v$ and $v^*$: the former is
killed on the left, and the latter on the right.  They also kill
$w$ and $w^*$, since $w$ has right label $\mu_n$ and $w^*$ has
left label $\mu_n$.  Consequently $C$ contains
\[
 X_i=e_{i1}x_i e_{1i}\quad(1\le i\le N),
 \qquad
 z_i=e_{i,i+1}d\quad(1\le i\le n-2).
\]
The roots of $z_i z_i^*$ and $z_i^*z_i$ give $d_i$ and
$d_{i+1}$, so all $d_i$, $i<n$, belong to $C$.
For $m\ge1$ put $f_m(u)=u/(u+m^{-1})$ and
\[
 z_{i,m}=z_i(d_{i+1}+m^{-1}1_A)^{-1}
        =e_{i,i+1}f_m(d)\in C.
\]
For $2\le i<n$ their products give
\[
 R_{i,m}=z_{1,m}\cdots z_{i-1,m}
        =e_{1i}f_m(d)^{i-1}\in C.
\]
Since $d_1$ is strictly positive in $D$,
\[
 R_{i,m}X_iR_{i,m}^*\longrightarrow x_i
 \quad\text{in norm}\qquad(2\le i\le N).
\]
Together with $X_1=x_1$, this shows $D\subseteq C$ in the first
corner.  Set $R_{1,m}=f_m(d_1)$.  For every $x\in D$ and
$i,j<n$, the approximate-unit property then gives
\[
 R_{i,m}^*xR_{j,m}\longrightarrow e_{i1}xe_{1j}
 \quad\text{in norm}.
\]
Thus
\begin{equation}\label{eq:recover-first-carrier}
 B_0\subseteq C.
\end{equation}

Choose $a_0\colon[0,1]\to[0,1]$ continuous, zero on
$[0,2\delta]$, and one on $[\gamma,1]$, and put
$a=e_{11}a_0(t)$.  This is a continuous function, vanishing at
zero, of $d_1$; hence $a\in B_0\subseteq C$.  On the first
matrix position, $q$ can be nonzero only where
$g_0(t)<3/4$, which is contained in $\{t<2\delta\}$.
Therefore $aq=0$, while $aw=w$.  In particular,
\[
 \|av^*\|^2=\|a v^*v a\|=\|aqa\|=0.
\]
Using the exact labels, $a\ell_1(h)=a$, and these two identities, we obtain
\begin{equation}\label{eq:recover-last-bridge}
 a\ell_1(h)k\ell_n(h)=w.
\end{equation}
Indeed, the core terms contribute only $w$; the $v$ term is
killed by the left label, and the $v^*$ term is killed by $a$.
Equations \eqref{eq:recover-first-carrier} and
\eqref{eq:recover-last-bridge} give $B=\mathrm C^*(B_0,w)\subseteq C$.
Subtracting the now-recovered core terms from $k$ yields
$v+v^*\in C$.

Let $(u_m)$ be the functional-calculus approximate unit of
the hereditary algebra generated by $b$.  Then $u_m\in B\subseteq C$,
$u_mv\to v$ in norm by $vv^*\in\operatorname{Her}_A(b)$, and
$u_mv^*=0$ because $bq=0$.  Hence
\[
 u_m(v+v^*)\longrightarrow v\quad\text{in norm},
\]
and $v\in C$.

It remains to verify that $B$ and $v$ generate $A$; the argument
also records why no support projection must be recovered in norm.
If $q=0$, \eqref{eq:cone-invertible-cover} makes $s$ an invertible
positive element of the hereditary algebra $B$, so $B=A$.
Suppose $q\ne0$, and put $E=\operatorname{Her}_A(q)$.
The identity
\[
 \overline{v^*\operatorname{Her}_A(vv^*)v}
 =\operatorname{Her}_A(v^*v)=E
\]
shows $E\subseteq \mathrm C^*(B,v)$.  Since $v^*v=q$ is strictly
positive in $E$, Lemma
\ref{lem:full-hereditary-bimodule} shows that this algebra also
contains $\overline{BAE}$ and $\overline{EAB}$.
Put $z=s+q$.  By \eqref{eq:cone-invertible-cover}, $z$ is
invertible.  Heredity gives
\[
 zAz\subseteq B+\overline{BAE}+\overline{EAB}+E
       \subseteq \mathrm C^*(B,v).
\]
Functional calculus, using functions vanishing at zero, gives
$1_A,z^{-1}\in \mathrm C^*(B,v)$, and thus
\[
 A=z^{-1}(zAz)z^{-1}\subseteq \mathrm C^*(B,v).
\]
Therefore $A=\mathrm C^*(B,v)\subseteq C=\mathrm C^*(h,k)$.
The reverse inclusion is automatic, so
$A=\mathrm C^*(h,k)=\mathrm C^*(h+ik)$.
\end{proof}

\begin{remark}[Scope of the compression proposition]
The proposition uses actual elements of $A$ throughout its encoder.
The supporting matrix units are used only as multiplier notation;
every occurrence in a generator or a recovered element is multiplied
by an element of an appropriate hereditary algebra.  No support
projection is claimed to lie in $A$.  The remaining input for a
global application is a uniform finite selfadjoint generator bound
for the selected hereditary coefficient algebra $D$.
\end{remark}

\section{Proof of the single-generation theorem}
\label{sec:main-proof}

We apply the preceding construction to a matrix cone in a unital
copy of the Jiang--Su algebra.

\begin{theorem}\label{thm:z-comparison}\label{thm:pure-comparison}
Let $A$ be a separable unital $\mathrm C^*$-algebra with
$\mathrm T(A)\ne\varnothing$ and strict comparison by ordinary
tracial states.  If $\mathcal Z$ embeds unitally in $A$, then $A$
is singly generated.
\end{theorem}

\begin{proof}[Proof of Theorems~\ref{thm:main} and~\ref{thm:z-comparison}]
Fix a unital copy of $\mathcal Z$ in $A$.  Choose $n=8$ and
$0<\varepsilon<1/(64n)$.  By Proposition~\ref{prop:z-cone}, there
is a c.p.c.\ order-zero map
$\phi:\mathrm M_n\to\mathcal Z\subseteq A$ such that, writing
$t=\phi(1_n)$,
\[
 \tau(t)>1-\varepsilon\qquad(\tau\in\mathrm T(A)).
\]
Set $\delta=1/16$ and
\[
 d=(\phi(e_{11})-\delta)_+,\qquad D=\Her_A(d).
\]
The element $d$ is nonzero and full in $A$, and has the same
strictly positive rank $c_\delta$ at every trace.
Lemma~\ref{lem:audited-constant-rank-corner} gives a compact nonempty
normalized trace space for $D$, affinely identified with
$\mathrm T(A)$.  Comparison of positive elements in $D$ agrees
with comparison in $A$, so $D$ inherits comparison by these
tracial states.

Proposition~\ref{prop:z-cone} gives almost divisibility of the
strictly positive class $[d]\in\mathrm{Cu}(D)$.
Corollary~\ref{cor:audited-regular-hereditary-completion}
shows that the normalized uniform tracial completion of $D$ is
regular, factorial, and of type $\mathrm{II}_1$.
Its designated traces form a metrizable Choquet simplex by the
trace identification above, and its unit ball is separable in the
uniform $2$-norm because $D$ is norm separable.

Apply Theorem~\ref{thm:external-af-core} to this completion.  It
supplies an external unital simple AF algebra $B$ with the
required common-centre and averaging properties.
Proposition~\ref{prop:finite-weakened-word} supplies finite
weakened-word estimates for coefficients in that completion.
The algebra $B$ has two selfadjoint generators
\cite[Theorem~7.3]{Thiel2021}.  Shifting one of them by a
sufficiently large scalar multiple of the unit makes it positive
and invertible; functional calculus then recovers the unit without
adjoining it.  Proposition~\ref{prop:audited-seven-generators}
therefore gives at most seven selfadjoint norm generators for $D$.

Proposition~\ref{prop:matrix-cone-compression} now applies with
$N=7$ and $n=8$.  Its coefficient algebra is precisely the $D$
chosen above, and $n-1=N$.  Thus there are $h,k\in A_{\sa}$ such
that $A=\mathrm C^*(h+ik)$.
\end{proof}

\begin{corollary}\label{cor:intrinsic}
Let $A$ be a separable unital pure $\mathrm C^*$-algebra of stable
rank one.  Suppose that every lower semicontinuous $2$-quasitrace
on $A$ is an ordinary trace.  Then $A$ is singly generated.
\end{corollary}

\begin{proof}
The zero algebra is immediate, so suppose $A\ne0$.
Lemma~\ref{lem:intrinsic-state-comparison} gives
$\mathrm T(A)\ne\varnothing$ and strict comparison by ordinary
tracial states.  Pureness gives almost unperforation of
$\mathrm{Cu}(A)$ and almost divisibility of $[1_A]$.
Proposition~\ref{prop:structural-z} gives a unital embedding of
$\mathcal Z$ into $A$.  Apply Theorem~\ref{thm:z-comparison}.
\end{proof}

\begin{corollary}\label{cor:simple-comparison}
Let $A$ be a simple separable unital $\mathrm C^*$-algebra with
$\mathrm T(A)\ne\varnothing$ and strict comparison by ordinary
tracial states.  Then $A$ is singly generated.
\end{corollary}

\begin{proof}
If $A$ is finite dimensional, it is isomorphic to $\mathrm M_m$.
A diagonal matrix with distinct positive diagonal entries and the
selfadjoint sum of adjacent matrix units generate $\mathrm M_m$;
for $m=1$ the unit alone suffices.

Suppose that $A$ is infinite dimensional.
Lemma~\ref{lem:simple-state-comparison-quasitraces} gives strict
comparison by the full ordinary trace cone, which is the convention
in \cite[Corollary~B]{Thiel2026}.  That result makes $A$ pure.
Its tracial states are
faithful, so $A$ is not purely infinite.  Lin's stable rank theorem
\cite[Corollary~4.5]{Lin2024} gives stable rank one.
Proposition~\ref{prop:structural-z} supplies a unital copy of
$\mathcal Z$, and Theorem~\ref{thm:z-comparison} applies.
\end{proof}

\begin{corollary}\label{cor:nowhere-scattered}
Let $A$ be a separable unital nowhere scattered
$\mathrm C^*$-algebra of stable rank one.  If $A$ has strict
comparison by the full cone of lower semicontinuous ordinary
traces, then $A$ is singly generated.
\end{corollary}

\begin{proof}
Ng--Robert's comparison theorem
\cite[Theorem~3.6 and Remark~3.7]{NgRobert2016} gives strict
comparison by $2$-quasitraces and shows that every lower
semicontinuous $2$-quasitrace is an ordinary trace.  Strict
comparison implies controlled comparison, so
\cite[Theorem~5.12]{APTV2024} shows that $A$ is pure.
Corollary~\ref{cor:intrinsic} applies.
\end{proof}

\begin{remark}\label{rem:factor-connection}
Corollary~\ref{cor:simple-comparison} also recovers the finite-factor conclusion
of \cite{OpenAI2026}.  Indeed, let $M$ be a $\mathrm{II}_1$ factor
with separable predual and normalized trace $\tau$.  A countable
witness construction gives a weakly dense separable simple unital
sub-$\mathrm C^*$-algebra $A\subseteq M$ with real rank zero, unique trace
$\tau|_A$, and strict comparison by ordinary traces.

Here are the ingredients of that construction.  Start with $1$
and a countable $2$-norm-dense subset of the unit ball of $M$.
At each stage, adjoin finite-spectrum norm approximants, unitaries
witnessing norm Dixmier averaging
\cite[Theorem~6.4.1 and Corollary~6.4.2]{AnantharamanPopa}, and
bounded approximate minimizers for
\[
 F_R(a,b)=\inf_{\substack{x\in\mathrm M_m(M)\\\|x\|\le R}}
                    \|a-x^*bx\|.
\]
For these choices use every accuracy $1/k$, $k\ge1$.  For the last
step take positive pairs from countable norm-dense lists at every
matrix level and integers $R\ge1$.  Adjoin the matrix entries of
the witnesses and repeat these choices for norm-dense lists in each
enlarged algebra.  The estimate
\[
 |F_R(a,b)-F_R(a',b')|
 \le\|a-a'\|+R^2\|b-b'\|
\]
ensures that, in the resulting algebra $A$, the infimum restricted
to $\mathrm M_m(A)$ equals the infimum in $\mathrm M_m(M)$.
For positive $a,b\in\mathrm M_m(A)$ with
$\mathrm d_\tau(a)<\mathrm d_\tau(b)$, choose $\delta>0$ so that
$\mathrm d_\tau(a)<\mathrm d_\tau((b-\delta)_+)$.
Put $p_\delta=1_{(\delta,\infty)}(b)$.
Projection comparison in the factor gives a partial isometry $w$
with $w^*w$ equal to the support of $a$ and $ww^*\le p_\delta$.
Then $x=(b|_{p_\delta})^{-1/2}wa^{1/2}\in\mathrm M_m(M)$ is
bounded and satisfies $x^*bx=a$.  Choose an integer $R>\|x\|$.
Equality of the
two $F_R$ infima gives $a\precsim b$ in $\mathrm M_m(A)$.
For stabilization, put $p_j=1_M\otimes\sum_{i=1}^j e_{ii}$.
For $a,b\in(A\otimes\mathcal K)_+$ with
$\mathrm d_\tau(a)<\mathrm d_\tau(b)$, one has
\[
 \mathrm d_\tau(p_jap_j)\le\mathrm d_\tau(a),\qquad
 \mathrm d_\tau(p_jbp_j)\uparrow\mathrm d_\tau(b).
\]
The second assertion follows by replacing $p_jbp_j$ with the
Cuntz-equivalent increasing sequence $b^{1/2}p_jb^{1/2}$ and
using normality of the matrix-amplified trace.
Choose $l$ with
$\mathrm d_\tau(p_lbp_l)>\mathrm d_\tau(a)$.
Comparison at finite matrix levels, after padding with zero blocks,
gives $p_jap_j\precsim p_lbp_l\precsim b$ for every $j$.
Since $p_jap_j\to a$ in norm, this proves $a\precsim b$.
The spectral approximants give real rank zero.
Dixmier averaging gives both uniqueness of the trace and
simplicity: every nonzero $a\in A_+$ averages in norm to the
nonzero scalar $\tau(a)1$.

Corollary~\ref{cor:simple-comparison} now gives
$A=\mathrm C^*(g)$ and therefore $M=\mathrm W^*(g)$.  This observation uses the
$\mathrm C^*$-algebra theorem directly, while the relative-generation
argument in its proof adapts the method of
\cite{OpenAI2026}.
\end{remark}

\begin{remark}
In the nonsimple setting, pointwise comparison by normalized
tracial states and uniform-ratio comparison by the full lower
semicontinuous trace cone must be distinguished.
Theorem~\ref{thm:z-comparison} assumes the former.
Corollary~\ref{cor:intrinsic} obtains it from the latter using stable
rank one and the existence of traces on nonzero unital quotients,
as proved in Lemma~\ref{lem:intrinsic-state-comparison}.
No faithfulness assumption is imposed on the individual tracial
states in these nonsimple results.
\end{remark}

\section{Dense generators in real rank zero}
\label{sec:rr0-density}

We first remove simplicity from the spectral-refinement argument.
The point is to choose a connecting element that generates the whole
rectangular bimodule, rather than merely a nonzero connecting element.
As before, $\mathrm C^*(S)$ denotes generation without automatically
adjoining a unit.

\begin{lemma}[Full rectangular bridges]
\label{lem:rr0-full-rectangular-bridges}
Let $A$ be a separable unital $\mathrm C^*$-algebra of real rank zero,
and let $p,q\in A$ be projections.  For every $x\in pAq$ and
$\varepsilon>0$, there is $z\in pAq$ such that
\[
 \|z-x\|<\varepsilon,\qquad
 pAq=\overline{\operatorname{span}}\,(pAp)z(qAq).
\]
The assertion includes the case $pAq=0$, when $z=0$.
\end{lemma}

\begin{proof}
We first check that closed-range elements are dense in every
rectangular corner $PAQ$, where $P,Q$ are projections in $A$.
Fix $u\in PAQ$.  Given small numbers $0<\kappa<\delta$, choose a
positive finite-spectrum element $d\in QAQ$ with
$\|d-u^*u\|<\kappa$.  This is possible because hereditary
subalgebras have real rank zero \cite{BrownPedersen1991}.
The spectral projection $e=1_{(\delta,\infty)}(d)$ belongs to $QAQ$,
and
\[
 e u^*u e\ge(\delta-\kappa)e,\qquad
 \|u-ue\|^2
 =\|(Q-e)u^*u(Q-e)\|\le\delta+\kappa.
\]
If $e=0$, put $y=0$ and $v=0$.  Otherwise set
\[
 y=ue,\qquad v=ue(eu^*ue)^{-1/2},
\]
where the inverse is taken in the unital corner $eAe$.
Then $y=v|y|$, $v^*v=e$, and $f=vv^*\le P$ is a projection.
Moreover, $|y|$ is invertible in $eAe$.
Letting $\delta,\kappa$ tend to zero proves the asserted density.
We shall use this construction also when one of $P,Q$ is zero.

Put $X=pAq$, and suppose $X\ne0$.
Choose a closed-range $y_0\in X$ with
$\|y_0-x\|<\varepsilon/2$.  Write
$y_0=v_0|y_0|$ and let $e_0=v_0^*v_0$, $f_0=v_0v_0^*$;
if $y_0=0$, take $v_0=e_0=f_0=0$.
Choose a sequence $(x_n)_{n\ge1}$ in the closed unit ball of $X$
obtained by repeating every member of a countable norm-dense
subset infinitely often.
Inductively, once $e_j,f_j$ are chosen for $0\le j<n$, put
\[
 P_{n-1}=p-\sum_{j<n}f_j,
 \qquad Q_{n-1}=q-\sum_{j<n}e_j,
 \qquad r_n=P_{n-1}x_nQ_{n-1}.
\]
Choose a closed-range element $y_n\in P_{n-1}AQ_{n-1}$ with
$\|y_n-r_n\|<1/n$.  Its polar decomposition in $A$ gives
$y_n=v_n|y_n|$, with
\[
 e_n=v_n^*v_n\le Q_{n-1},\qquad
 f_n=v_nv_n^*\le P_{n-1}.
\]
For a zero $y_n$ all three elements are taken to be zero.
Thus the $e_n$ are pairwise orthogonal, as are the $f_n$.
Choose strictly positive scalars $c_n$ with
$\sum_{n\ge1}c_n<\varepsilon/2$, and define the norm-convergent sum
\[
 z=y_0+\sum_{n\ge1}c_nv_n.
\]
Then $z\in X$ and $\|z-x\|<\varepsilon$.

Let $I$ be a closed two-sided ideal of $A$ containing $z$.
The orthogonality on both sides gives
\[
 z^*z=y_0^*y_0+\sum_{n\ge1}c_n^2e_n.
\]
Consequently $e_n\in I$ for $n\ge1$, since each positive summand
is dominated by $z^*z$.  Also $y_0^*y_0\in I$; if $e_0\ne0$,
invertibility of $|y_0|$ in $e_0Ae_0$ gives a number $\alpha>0$
with $y_0^*y_0\ge\alpha e_0$, so $e_0\in I$ as well.
It follows that $v_n,f_n,y_n\in I$ for every $n\ge0$.
If $\pi:A\to A/I$ is the quotient map, then
$\pi(r_n)=\pi(x_n)$, and hence
\[
 \|\pi(x_n)\|=\|\pi(r_n)\|
       \le\|r_n-y_n\|<1/n.
\]
Every member of the chosen dense subset occurs at arbitrarily
large indices, so $\pi(X)=0$.
Applying this observation to $I=\operatorname{Ideal}_A(z)$ proves
\begin{equation}\label{eq:rr0-bridge-ideal}
 \operatorname{Ideal}_A(z)=\operatorname{Ideal}_A(X).
\end{equation}

We finish by passing from ideals to the rectangular bimodule.
Set $B=pAp$, $E=qAq$, and $K=\operatorname{Ideal}_A(X)$.
The hereditary ideal correspondence gives
\[
 J_X:=\overline{\operatorname{span}}X^*X=E\cap K.
\]
Indeed, $J_X$ is an ideal of $E$ and generates $K$ as an ideal of
$A$.  The same correspondence and
\eqref{eq:rr0-bridge-ideal} show that
\[
 J:=\overline{\operatorname{span}}Ez^*BzE
   =\operatorname{Ideal}_E(z^*z)
   =E\cap K=J_X.
\]
For the first equality with the generated ideal, use $p\in B$
and $z^*bz\le\|b\|z^*z$ for $b\in B_+$.
The right action of $J_X$ on $X$ is nondegenerate: if $(a_\lambda)$
is an approximate unit of $J_X$, then
$\|u-ua_\lambda\|^2
 =\|(q-a_\lambda)u^*u(q-a_\lambda)\|\to0$ for $u\in X$.
Since $XEz^*B\subseteq B$, we have
\[
 XJ\subseteq\overline{\operatorname{span}}BzE.
\]
Thus $X\subseteq\overline{\operatorname{span}}BzE$.
The reverse inclusion is immediate.
\end{proof}

\begin{proposition}\label{prop:rr0-density-criterion}
Let $A$ be a separable unital $\mathrm C^*$-algebra of real rank zero.
Suppose that $pAp$ is singly generated for every nonzero projection
$p\in A$.  Then
\[
 \operatorname{Gen}(A)=\{g\in A:\mathrm C^*(g)=A\}
\]
is a norm-dense $G_\delta$ subset of $A$.
Equivalently, $A$ has generator rank at most one.
\end{proposition}

\begin{proof}
If $A=0$, the conclusion is immediate, so assume $A\ne0$.
Fix $a,b\in A_{\sa}$ and $\varepsilon>0$, and choose
$0<\eta<\varepsilon/6$.
Real rank zero gives a finite partition
$1=\sum_{i=1}^m p_i$ by nonzero projections and real numbers
$\lambda_i$ such that
\[
 \left\|a-\sum_i\lambda_i p_i\right\|<\eta.
\]
Approximate $p_i b p_i$ within $\eta$ by a finite-spectrum
selfadjoint element of the real-rank-zero corner $p_iAp_i$.
Relabelling its spectral projections and discarding zero projections gives a
refinement $1=\sum_{\ell=1}^Lq_\ell$, with
$q_\ell\le p_{i(\ell)}$, and real numbers $\mu_\ell$ satisfying
\[
 \|q_\ell bq_\ell-\mu_\ell q_\ell\|<\eta.
\]
The block-diagonal element
\[
 \Delta=\sum_{\ell=1}^L
       (\mu_\ell q_\ell-q_\ell bq_\ell)
\]
has norm less than $\eta$.

Choose pairwise distinct nonzero real numbers $\alpha_\ell$ with
$|\alpha_\ell-\lambda_{i(\ell)}|<\eta$.
By hypothesis there are selfadjoint contractions
$H_\ell,K_\ell\in q_\ell Aq_\ell$ such that
$q_\ell Aq_\ell=\mathrm C^*(H_\ell,K_\ell)$.
Choose $0<\delta<\eta$ so that the intervals
$[\alpha_\ell-\delta,\alpha_\ell+\delta]$ are pairwise disjoint
and avoid zero, and put
\[
 h=\sum_{\ell=1}^L(\alpha_\ell q_\ell+\delta H_\ell).
\]
Continuous functional calculus, with functions vanishing at zero,
recovers $1$ and every $q_\ell$ from $h$.  Consequently
\[
 H_\ell=\delta^{-1}(q_\ell hq_\ell-\alpha_\ell q_\ell)
            \in\mathrm C^*(h).
\]

For every unordered pair $\ell<j$, use
Lemma~\ref{lem:rr0-full-rectangular-bridges} to choose
$z_{\ell j}\in q_\ell Aq_j$ such that
\[
 \|z_{\ell j}-q_\ell bq_j\|<\eta/L^2,
 \qquad
 q_\ell Aq_j=
 \overline{\operatorname{span}}(q_\ell Aq_\ell)z_{\ell j}(q_jAq_j).
\]
When $q_\ell Aq_j=0$, choose $z_{\ell j}=0$.
Set
\[
\begin{split}
 w&=\sum_{\ell<j}
    \bigl(z_{\ell j}-q_\ell bq_j
          +z_{\ell j}^*-q_jbq_\ell\bigr),\\
 k&=b+\Delta+\delta\sum_{\ell=1}^L K_\ell+w,
 \qquad C=\mathrm C^*(h,k).
\end{split}
\]
The sum defining $w$ has $L(L-1)$ terms, so $\|w\|<\eta$;
for $L=1$ it is zero.
Its diagonal compressions vanish, and hence
\[
 K_\ell=\delta^{-1}(q_\ell kq_\ell-\mu_\ell q_\ell)\in C.
\]
Thus $C$ contains every $q_\ell Aq_\ell$.
Also $q_\ell kq_j=z_{\ell j}\in C$ for $\ell<j$.
The full-bridge property now puts every $q_\ell Aq_j$ and its
adjoint in $C$.  Since $\sum_\ell q_\ell=1$, we obtain $C=A$.
Finally,
\[
 \|h-a\|<3\eta,\qquad \|k-b\|<3\eta,
\]
so $h+ik$ generates $A$ and
$\|(h+ik)-(a+ib)\|<6\eta<\varepsilon$.

For the $G_\delta$ assertion, choose a norm-dense sequence $(a_j)$
in $A$.  Then
\[
 \operatorname{Gen}(A)
 =\bigcap_{j,m\ge1}
 \left\{g\in A:
       \operatorname{dist}(a_j,\mathrm C^*(g))<1/m\right\}.
\]
Each set in braces is open: its strict inequality is witnessed by
a fixed $*$-polynomial in $g,g^*$ with zero constant term, whose
evaluation is norm continuous.  The generator-rank formulation is
\cite[Theorem~3.4 and Remark~3.7]{Thiel2021}.
\end{proof}

\begin{lemma}\label{lem:rr0-intrinsic-corners}
Let $A$ be a unital pure $\mathrm C^*$-algebra of stable rank one.
Suppose every lower semicontinuous $2$-quasitrace on $A$ is a trace.
Every nonzero projection corner $pAp$ has these same three properties.
\end{lemma}

\begin{proof}
Let $I=\operatorname{Ideal}_A(p)$.
The full hereditary inclusion $pAp\subseteq I$ identifies
$\mathrm{Cu}(pAp)$ with the ideal
$J=\mathrm{Cu}(I)$ of $\mathrm{Cu}(A)$.
Almost unperforation passes to this ideal.
For almost divisibility, take $x'\ll x$ in $J$ and $n\ge1$.
An ambient witness $z$ satisfies $nz\le x$ and
$x'\le(n+1)z$; since $z\le nz\le x$ and $J$ is hereditary,
$z\in J$.  Thus $pAp$ is pure.

For completeness, stable rank one passes to this corner by a
Schur-complement argument.  If $p=1$ there is nothing to prove.
For $a\in pAp$, approximate $a+(1-p)$ by invertible elements
$u\in A$.  The block $(1-p)u(1-p)$ is invertible in
$(1-p)A(1-p)$ for sufficiently close approximants.
The Schur complement
\[
 pup-pu(1-p)\bigl((1-p)u(1-p)\bigr)^{-1}(1-p)up
\]
is invertible in $pAp$ and converges to $a$.

Let $\rho$ be a lower semicontinuous $2$-quasitrace on $pAp$.
The quasitrace--functional correspondence
\cite[Proposition~4.2 and Theorem~4.4]{ERS2011} transports
$\mathrm d_\rho$ to a functional $\lambda:J\to[0,\infty]$.
Extend it to $\mathrm{Cu}(A)$ by
\[
 \widetilde\lambda(x)=
 \begin{cases}
  \lambda(x),&x\in J,\\
  \infty,&x\notin J.
 \end{cases}
\]
This extension is additive and order preserving because $J$ is
an ideal.  It preserves increasing suprema because $J$ is closed
under such suprema: an increasing sequence with supremum outside
$J$ must already have a term outside $J$.
Thus $\widetilde\lambda$ is a functional on $\mathrm{Cu}(A)$.
By the same correspondence it defines a lower semicontinuous
$2$-quasitrace $\widetilde\rho$ on $A$ whose restriction to $pAp$
is $\rho$.
In fact the inverse formula
\[
 \widetilde\rho(a)
   =\int_0^\infty
        \widetilde\lambda([(a-t)_+])\,\mathrm dt
\]
shows the restriction assertion directly.
The hypothesis makes $\widetilde\rho$ a trace, so $\rho$ is a trace.
The extension is allowed to be unbounded or not densely finite;
this is why the full lower semicontinuous quasitrace hypothesis
is used.
\end{proof}

\begin{theorem}\label{thm:rr0-intrinsic-density}
Let $A$ be a separable unital pure $\mathrm C^*$-algebra with real
rank zero and stable rank one.  Suppose that every lower
semicontinuous $2$-quasitrace on $A$ is an ordinary trace.
Then $\operatorname{Gen}(A)$ is a norm-dense $G_\delta$ subset of $A$.
Equivalently, $A$ has generator rank at most one.
\end{theorem}

\begin{proof}
The zero algebra is immediate.  Every nonzero projection corner
$pAp$ satisfies the intrinsic hypotheses by
Lemma~\ref{lem:rr0-intrinsic-corners}, and is separable and unital.
Corollary~\ref{cor:intrinsic} makes it singly generated.
Apply Proposition~\ref{prop:rr0-density-criterion}.
\end{proof}

\begin{corollary}\label{cor:rr0-simple-density}
Let $A$ be a simple separable unital $\mathrm C^*$-algebra with
$\mathrm T(A)\ne\varnothing$, real rank zero, and strict comparison
by tracial states.  Then $\operatorname{Gen}(A)$ is a norm-dense
$G_\delta$ subset of $A$, and $A$ has generator rank at most one.
\end{corollary}

\begin{proof}[Proof of Corollaries~\ref{cor:rr0-density}
                 and~\ref{cor:rr0-simple-density}]
If $A$ is finite dimensional, it is a full matrix algebra.
Every nonzero projection corner is again a full matrix algebra
and is singly generated, so
Proposition~\ref{prop:rr0-density-criterion} applies.
Otherwise Lemma~\ref{lem:simple-state-comparison-quasitraces}
gives comparison by the full ordinary trace cone and the
quasitrace hypothesis.  Hence $A$ is pure by
\cite[Corollary~B]{Thiel2026} and has stable rank one by
\cite[Corollary~4.5]{Lin2024}.
Theorem~\ref{thm:rr0-intrinsic-density} now applies.
\end{proof}

\begin{remark}
The density theorem uses single generation of all projection
corners.  A prescribed unital $\mathcal Z$-embedding into $A$
does not by itself supply unital $\mathcal Z$-embeddings into
nonfull corners, and traces on such corners need not extend
boundedly to $A$.
Thus the argument above does not establish density under the
hypotheses of Theorem~\ref{thm:main} and real rank zero alone.
Under the stronger assumption of $\mathcal Z$-stability, density
for arbitrary separable real-rank-zero algebras is already known
from \cite[Theorem~5.3]{Thiel2022Generators}, without a comparison
hypothesis.
\end{remark}

\section{Polar generators and topological free entropy dimension}
\label{sec:entropy}

The matrix-cone construction also controls the topological free
entropy dimension of its generator.  The control comes from the
constant spectral labels of the first selfadjoint element and the
sparse matrix positions occupied by the second.  We first specify
the microstate convention and the change of coordinates used below.

Let $x=(x_1,\ldots,x_\ell)$ be a selfadjoint tuple in a unital
$\mathrm C^*$-algebra, and fix
$R>\max_j\|x_j\|$.  For a finite set $\mathcal F$ of complex
noncommutative polynomials in $\ell$ selfadjoint variables and
$\sigma>0$, put
\[
 \begin{split}
 \Gamma_R^{\mathrm{top}}(x;q,\mathcal F,\sigma)
 =\{X\in(\mathrm M_q^{\sa})^\ell:\;&
       \max_j\|X_j\|\le R,\\[-2pt]
 &\bigl|\|P(X)\|-\|P(x)\|\bigr|<\sigma
                         \quad(P\in\mathcal F)\}.
 \end{split}
\]
For a bounded set $\Omega\subseteq(\mathrm M_q^{\sa})^\ell$, let
$N_\infty(\Omega,\omega)$ be its covering number for balls of radius
$\omega$ in the maximum operator norm.  We allow the centres to
belong to the ambient matrix space.  Requiring centres in $\Omega$
changes the covering radius by at most a factor of two and gives
the same dimension.  With $\log0=-\infty$, set
\begin{equation}\label{eq:entropy-fixed-scale}
 \mathbb K_R(x;\omega)
 =\inf_{\mathcal F,\sigma>0}
   \limsup_{q\to\infty}\frac1{q^2}
    \log N_\infty\bigl(
       \Gamma_R^{\mathrm{top}}(x;q,\mathcal F,\sigma),\omega\bigr).
\end{equation}
Then
\begin{equation}\label{eq:entropy-dimension-definition}
 \delta_{\mathrm{top}}(x)
 =\limsup_{\omega\downarrow0}
       \frac{\mathbb K_R(x;\omega)}{|\log\omega|}.
\end{equation}
The value is independent of the choice of
$R>\max_j\|x_j\|$.  Indeed, the coordinate polynomials already
bound every sufficiently accurate microstate by any such $R$.
These are the operator-norm microstate definitions of
\cite{VoiculescuTopological}.  For a possibly nonselfadjoint element
$g$, we use the notation
\[
 \delta_{\mathrm{top}}(g)
   =\delta_{\mathrm{top}}(\re g,\im g),
 \qquad
 \chi_{\mathrm{top}}(g)
   =\chi_{\mathrm{top}}(\re g,\im g).
\]
Here $\chi_{\mathrm{top}}(x)$ is defined using the same norm
microstates by replacing the covering expression in
\eqref{eq:entropy-fixed-scale} with
\[
 \limsup_{q\to\infty}
 \left(\frac1{q^2}\log\operatorname{vol}
          \Gamma_R^{\mathrm{top}}(x;q,\mathcal F,\sigma)
                  +\frac\ell2\log q\right),
\]
and then taking the infimum over $\mathcal F,\sigma$.
The volume is Lebesgue measure for the real Hilbert space structure
given by the unnormalized Hilbert--Schmidt inner product.

We shall also use the elementary implication
\begin{equation}\label{eq:entropy-dimension-volume}
 \delta_{\mathrm{top}}(x)<\ell
                    \quad\Longrightarrow\quad
                         \chi_{\mathrm{top}}(x)=-\infty.
\end{equation}
Indeed, a ball of radius $\omega$ in the maximum operator norm
is contained in a Hilbert--Schmidt ball of radius
$\sqrt{\ell q}\,\omega$ in real dimension $\ell q^2$.
The Euclidean ball volume formula and Stirling's estimate give
\[
 \limsup_{q\to\infty}
 \left(\frac1{q^2}\log\operatorname{vol}B_\infty(0,\omega)
                                  +\frac\ell2\log q\right)
                         \le\ell\log\omega+C_\ell.
\]
Covering a norm-microstate set by these balls and then taking the
infimum over its tests gives
$\chi_{\mathrm{top}}(x)\le
\mathbb K_R(x;\omega)+\ell\log\omega+C_\ell$.
If $\delta_{\mathrm{top}}(x)<d<\ell$, then for all sufficiently
small $\omega$ the right side is at most
$(\ell-d)\log\omega+C_\ell$, which tends to $-\infty$.

A separable unital $\mathrm C^*$-algebra is MF if it admits a unital
embedding into a quotient
$\prod_j\mathrm M_{d_j}/\bigoplus_j\mathrm M_{d_j}$, where the
product consists of bounded sequences and the direct sum consists
of sequences converging to zero in norm.  Such an algebra has a
tracial state, obtained by taking an ultralimit of normalized matrix
traces.  If a unital MF algebra has no nonzero finite-dimensional
representations, then every finite selfadjoint generating tuple
satisfies
\begin{equation}\label{eq:entropy-mf-lower}
                    \delta_{\mathrm{top}}(x)\ge1.
\end{equation}
This is \cite[Theorem~5.3(i)]{MFEntropyLowerBound}.  A unital
embedding of $\mathcal Z$ into $A$ excludes nonzero
finite-dimensional representations of $A$: restricting one to
$\mathcal Z$ would give a nonzero finite-dimensional
representation of the simple infinite-dimensional algebra
$\mathcal Z$.  Thus \eqref{eq:entropy-mf-lower} applies to every
MF algebra in Theorem~\ref{thm:z-comparison}.  The MF assumption
is needed for this lower bound; the upper estimate below does not
require matricial approximation.

\begin{lemma}[Polar coordinates and topological free entropy dimension]
\label{lem:polar-coordinates}
Let $A$ be a unital $\mathrm C^*$-algebra, and let $h,k\in A_{\sa}$
satisfy $h\geq 1_A$ and $A=\mathrm C^*(h,k)$, where the generated
algebra is not assumed to contain a unit in advance.  Then
$g=h+ik$ is invertible, its polar decomposition
\[
 g=u|g|,\qquad u=g(g^*g)^{-1/2},
\]
has $u\in\mathcal U(A)$, and $A=\mathrm C^*(u,|g|)$.

For every $\alpha>0$ with $\alpha\|k\|<\pi$, the element
\[
 G_\alpha=e^{i\alpha k}h
\]
satisfies
\begin{equation}\label{eq:polar-generator-coordinates}
 |G_\alpha|=h,\qquad
 \operatorname{phase}(G_\alpha)=e^{i\alpha k},\qquad
 \mathrm C^*(G_\alpha)=A.
\end{equation}
Moreover,
\begin{equation}\label{eq:polar-dimension-preservation}
 \delta_{\mathrm{top}}(\operatorname{Re}G_\alpha,
                         \operatorname{Im}G_\alpha)
 =\delta_{\mathrm{top}}(h,k).
\end{equation}
Here the same convention for empty norm-microstate spaces is used
on both sides.
\end{lemma}

\begin{proof}
The selfadjoint element $a=h^{-1/2}kh^{-1/2}$ gives the factorization
\[
 g=h^{1/2}(1_A+ia)h^{1/2}.
\]
The middle factor is invertible by continuous functional calculus,
and $\|(1_A+ia)^{-1}\|\leq1$.  Hence $g$ is invertible and
$\|g^{-1}\|\leq1$.  In particular, $|g|\geq1_A$, and the displayed
formula for $u$ defines a unitary in $A$.  Functional calculus in
$\mathrm C^*(g)$ recovers $|g|$, $1_A$, $|g|^{-1}$, and $u$.
Conversely $g=u|g|$, while $h$ and $k$ are its real and imaginary
parts.  This proves the first assertion, including the assertion
about generation without adjoining a unit.

For the alternative generator, direct multiplication gives
$G_\alpha^*G_\alpha=h^2$.  Thus its modulus and unitary phase are
as in \eqref{eq:polar-generator-coordinates}.  Its modulus again
recovers the unit and its inverse.  Since
$\alpha\sigma(k)\subset(-\pi,\pi)$, the principal argument is
continuous on the spectrum of $e^{i\alpha k}$ and gives
\[
 k=\alpha^{-1}\operatorname{Arg}(e^{i\alpha k}).
\]
It follows that $\mathrm C^*(G_\alpha)=\mathrm C^*(h,k)=A$.

We prove \eqref{eq:polar-dimension-preservation} by comparing
norm-microstate covering numbers under these explicit coordinate
maps.  In particular, no invariance under an arbitrary change of
generators is used.  Choose constants
\[
 0<c<1,\qquad M>\|h\|,\qquad
 \alpha\|k\|<\theta<\pi.
\]
For every matrix size $q$, put
\begin{align*}
 \mathcal D_q
   &=\{(H,K)\in(\mathrm M_q(\mathbb C)_{\sa})^2:
          c1_q\leq H\leq M1_q,\ \|\alpha K\|\leq\theta\},\\
 \mathcal E_q
   &=\{Z\in\mathrm M_q(\mathbb C):
          c1_q\leq|Z|\leq M1_q,\
          \sigma(U_Z)\subseteq\{e^{it}:|t|\leq\theta\}\},
 \qquad U_Z=Z|Z|^{-1}.
\end{align*}
The maps
\begin{equation}\label{eq:polar-matrix-coordinate-maps}
 \Phi_q(H,K)=e^{i\alpha K}H,\qquad
 \Psi_q(Z)=\bigl(|Z|,\alpha^{-1}\operatorname{Arg}(U_Z)\bigr)
\end{equation}
are inverse bijections between $\mathcal D_q$ and $\mathcal E_q$.
We identify a complex matrix $Z$ with its selfadjoint pair
$(\operatorname{Re}Z,\operatorname{Im}Z)$; the corresponding norms
satisfy
\[
 \max\{\|\operatorname{Re}Z\|,\|\operatorname{Im}Z\|\}
 \leq\|Z\|
 \leq2\max\{\|\operatorname{Re}Z\|,\|\operatorname{Im}Z\|\}.
\]

Both maps in \eqref{eq:polar-matrix-coordinate-maps} are Lipschitz
for these operator-norm metrics, with constants independent of $q$.
For $\Phi_q$, the integral formula for the difference of two
exponentials of selfadjoint matrices gives
\[
 \|e^{i\alpha K}H-e^{i\alpha K'}H'\|
 \leq\|H-H'\|+\alpha M\|K-K'\|.
\]
For completeness, if $Z,Z'\in\mathcal E_q$ and
$P=|Z|$, $P'=|Z'|$, the identity
\[
 P(P-P')+(P-P')P'=Z^*Z-Z'^*Z'
\]
and its integral solution using $e^{-sP}$ and $e^{-sP'}$ give
\[
 \|P-P'\|\leq\frac{1}{2c}\|Z^*Z-Z'^*Z'\|
 \leq\frac Mc\|Z-Z'\|.
\]
The identity $P^{-1}-P'^{-1}=P^{-1}(P'-P)P'^{-1}$ then shows
that $Z\mapsto U_Z=ZP^{-1}$ is uniformly Lipschitz on these
domains.  Finally, the principal logarithm is holomorphic on a
fixed neighborhood of the compact arc
$\{e^{it}:|t|\leq\theta\}$.  A fixed contour in that neighborhood,
the resolvent identity, and the normality of $U_Z$ show that
$U_Z\mapsto\operatorname{Arg}(U_Z)=-i\log U_Z$ is Lipschitz with
a constant depending only on $\theta$.  All these estimates are
independent of the matrix size.

We explain why the coordinate maps apply to the norm-microstate
spaces in the definition of $\delta_{\mathrm{top}}$.  Fix norm
cutoffs larger than the actual coordinate norms for each tuple;
coordinate norm tests make their precise values immaterial.
Sufficiently
accurate finite norm-polynomial tests for $(h,k)$ force its matrix
microstates into $\mathcal D_q$.  Indeed, norm tests for $H$, $K$,
and $H-S1_q$, with a fixed $S>\|h\|$, enforce the upper bounds
and the positive lower bound, using the strict margins in the
choices of $c$, $M$, and $\theta$.

The corresponding assertion for
$(\operatorname{Re}G_\alpha,\operatorname{Im}G_\alpha)$ also follows
from finitely many norm-polynomial tests.  First test $Z$ and
$1_q-Z^*Z/S^2$, where $S>\|G_\alpha\|$.  Since
$G_\alpha^*G_\alpha=h^2\geq1_A$, these tests, with sufficiently
small tolerance, force $c1_q\leq|Z|\leq M1_q$.  On this bounded
positive-gap set, $Z\mapsto U_Z$ is uniformly approximable by
noncommutative $*$-polynomials in $Z,Z^*$, by approximating
$s\mapsto s^{-1/2}$ on a fixed positive interval.  The strict
inequality
\[
 \|e^{i\alpha k}-1_A\|
 <2\sin(\theta/2)
\]
therefore allows one further finite polynomial norm test to force
$\|U_Z-1_q\|\leq2\sin(\theta/2)$.  Since $U_Z$ is unitary,
this is precisely the required phase-arc condition.  Thus these
microstates lie in $\mathcal E_q$.

On $\mathcal D_q$ the exponential in $\Phi_q$ is uniformly
approximable by polynomials, independently of $q$.  On
$\mathcal E_q$, the square root, inverse square root, and
principal argument in $\Psi_q$ are likewise uniformly
approximable by $*$-polynomials: use ordinary polynomial
approximation on the fixed positive spectral interval and
Laurent-polynomial approximation on the fixed unitary arc.
Consequently, any prescribed finite collection of norm-microstate
conditions for either tuple pulls back to a finite collection for
the other tuple, with any prescribed positive tolerance.  This
follows by first approximating the coordinate maps, then
substituting their polynomial approximants into the finitely
many test polynomials.  The approximation errors are uniform
in $q$.

Let $L\geq1$ be a common Lipschitz constant for the two coordinate
maps, with the above identification of complex matrices and
selfadjoint pairs.  Given any finite microstate conditions for one
tuple, strengthen the conditions for the other tuple so that its
image lies in the prescribed microstate space and the relevant
gap conditions hold.  Applying the inverse coordinate map to
a covering gives a covering with radius at most $2L$ times the
original radius.  The factor $2$ permits the centers of the
original covering balls to be replaced by points of the covered
set.  This yields the covering-number inequalities in both
directions, with the covering radius changed by the fixed factor
$2L$.  Write $x=(h,k)$ and
$y=(\operatorname{Re}G_\alpha,\operatorname{Im}G_\alpha)$,
and let $R_x,R_y$ be their chosen norm cutoffs.  Taking the
normalized $q^{-2}\log$ covering numbers, the limsup as
$q\to\infty$, and the infimum over finite tests gives
\[
 \mathbb K_{R_y}(y;2L\omega)\leq\mathbb K_{R_x}(x;\omega),
 \qquad
 \mathbb K_{R_x}(x;2L\omega)\leq\mathbb K_{R_y}(y;\omega).
\]
Finally let the covering radius tend to zero.
Multiplication of that radius by $2L$ does not change the
result, because
\[
 \frac{|\log((2L)\omega)|}{|\log\omega|}\longrightarrow1.
\]
This proves \eqref{eq:polar-dimension-preservation}.  The same
microstate transfers show that failure of matricial
approximation occurs simultaneously for the two tuples.
\end{proof}

\begin{remark}
The polar form above imposes no Haar or freeness condition.
In fact, the phase $u$ of the original strictly accretive
generator $g=h+ik$ cannot be Haar for any tracial state.
To see this, put $p=|g|$ and $B=\|g\|$.  If
$\lambda\in\sigma(u)$, unit approximate eigenvectors for $u$
in a faithful representation give
\[
 1\leq(\operatorname{Re}\lambda)
          \langle p\xi,\xi\rangle+o(1).
\]
Positivity of $p$ first forces $\operatorname{Re}\lambda>0$;
then $\langle p\xi,\xi\rangle\leq B$ gives
$\operatorname{Re}\lambda\geq B^{-1}$.  Hence
$\operatorname{Re}u\geq B^{-1}1_A$, whereas a Haar unitary has
zero trace.  No atomlessness or free independence of the polar
factors is asserted.
\end{remark}

We record the trace consequence of norm microstates that will be
used in the covering argument.  A continuous matrix expression
means an expression obtained from noncommutative polynomials by
finitely many algebraic operations and continuous functional
calculus on selfadjoint elements in fixed bounded intervals.
Such expressions are uniformly approximable by noncommutative
polynomials on a fixed norm-bounded set of matrix tuples, in all
matrix sizes.  This follows successively from polynomial
approximation of the scalar continuous functions.

\begin{lemma}\label{lem:entropy-trace-control}
Let $A=\mathrm C^*(1,x_1,\ldots,x_\ell)$, where the $x_j$ are
selfadjoint, and fix $R>\max_j\|x_j\|$.  Suppose that $\Psi$ is
a continuous matrix expression whose values are selfadjoint on
selfadjoint tuples, and
\[
             \tau(\Psi(x))\ge b\qquad(\tau\in\T(A)).
\]
For every $\eta>0$ there are finite norm-microstate tests
$\mathcal F,\sigma$ such that, in every matrix size,
\[
 \operatorname{tr}_q(\Psi(X))>b-\eta
 \quad\text{for }X\in
       \Gamma_R^{\mathrm{top}}(x;q,\mathcal F,\sigma).
\]
Here $\operatorname{tr}_q=q^{-1}\operatorname{Tr}$.
Norms of finitely many continuous matrix expressions can also be
prescribed within any given tolerance by finite norm-microstate
tests.
\end{lemma}

\begin{proof}
For the last assertion, approximate each expression uniformly by
a polynomial on the $R$-ball and include those polynomials in the
tests.  To prove the trace assertion, suppose it fails.  Enumerate
the polynomials with rational complex coefficients and choose
$R$-bounded tuples $X^{(j)}$ of sizes $q_j$ which satisfy the first
$j$ norm tests within $1/j$ but have
\[
       \operatorname{tr}_{q_j}(\Psi(X^{(j)}))\le b-\eta.
\]
Density of the enumerated polynomials and the norm bound give
\[
           \|P(X^{(j)})\|\longrightarrow\|P(x)\|
\]
for every polynomial $P$.  Thus $x_i\mapsto[(X_i^{(j)})_j]$
extends to a unital isometric $*$-homomorphism from $A$ into the norm
quotient of the bounded product of the $\mathrm M_{q_j}$ by the
norm-null sequences.  An ultralimit of the normalized matrix traces
is a tracial state on this quotient.  Its pullback $\tau$ to $A$
satisfies
\[
 \tau(\Psi(x))
   =\lim_{j\to\mathcal U}
        \operatorname{tr}_{q_j}(\Psi(X^{(j)}))\le b-\eta.
\]
Continuous functional calculus commutes with the quotient map, so
the equality applies to the stated expressions.  This contradicts
the hypothesis.
\end{proof}

\begin{proposition}\label{prop:entropy-cone}
Use the hypotheses and the actual elements $h,k$ constructed in
the proof of Proposition~\ref{prop:matrix-cone-compression}.
Write $m$ for the matrix size denoted by $n$ there; thus
$m\ge3$, $1\le N\le m-1$, and
$0<\varepsilon<1/(64m)$.  Then $g=h+ik$ is an invertible
generator of $A$, and
\begin{equation}\label{eq:entropy-cone-upper}
 \delta_{\mathrm{top}}(g)
 \le 2-
      \frac{m^2-m+2-N}{m^2}(1-2\varepsilon)^2,
 \qquad
 \chi_{\mathrm{top}}(g)=-\infty.
\end{equation}
If $A$ is $\mathrm{MF}$ and has no nonzero finite-dimensional
representations, then also $\delta_{\mathrm{top}}(g)\ge1$.  For the polar decomposition
$g=u|g|$, one has
\[
                 A=\mathrm C^*(u,|g|)=\mathrm C^*(g).
\]
\end{proposition}

\begin{proof}
Retain the notation of that construction, including its supporting
matrix units $e_{ij}$, the positive contraction $t$, and
$\delta=1/16$.  Thus
\[
 h=1+\sum_{i=1}^m\lambda_i e_{ii}F(t),\qquad
       \mu_i=1+\lambda_i,
\]
where the $\lambda_i$ are distinct positive numbers and
$F=1$ on $[\delta,1]$.  The diagonal terms of $k$ occur only in
positions $1,\ldots,N$.  Its other prescribed cone terms lie on
the edges of the tree
\begin{equation}\label{eq:entropy-tree}
 \mathcal E
  =\{\{i,i+1\}:1\le i\le m-2\}\cup\{\{1,m\}\}.
\end{equation}
The remaining term is $v+v^*$, where $v^*v=q$ and $q$ is
orthogonal to every continuous function of $t$ supported in
$(1/2,1]$.  In particular, this term disappears when it is cut
on both sides by such functions of $t$.

Fix $a$ with $1/2<a<1-\varepsilon$, and choose a continuous
function $c:[0,1]\to[0,1]$ which is zero on $[0,1/2]$ and one
on $[a,1]$.  Put
\begin{equation}\label{eq:entropy-beta-a}
              \beta_a=\frac{1-\varepsilon/(1-a)}m>0.
\end{equation}
For any $\rho>0$ small enough that the intervals
$(\mu_i-\rho,\mu_i+\rho)$ are pairwise disjoint, choose continuous
functions $0\le f_i\le1$ supported in these intervals and
satisfying $f_i(\mu_i)=1$.  The constant labels of $h$ on the
support of $c(t)$ give the exact identities
\begin{equation}\label{eq:entropy-positive-cuts}
 a_i:=f_i(h)c(t)^2f_i(h)=e_{ii}c(t)^2
                         \qquad(1\le i\le m).
\end{equation}
These identities do not require the supporting matrix units to
belong to $A$; the displayed products belong to $A$ by
\eqref{eq:cone-membership}.

Let $\mathcal J$ be the following set of ordered matrix positions:
\[
 \mathcal J
  =\{(i,j):i\ne j,\ \{i,j\}\notin\mathcal E\}
      \cup\{(i,i):N<i\le m\}.
\]
The construction of $k$ implies
\begin{equation}\label{eq:entropy-exact-zeros}
                    a_i k a_j=0\qquad((i,j)\in\mathcal J).
\end{equation}
Indeed, every cone term occupies one of the listed diagonal
positions or tree edges.  Also $c(t)q=0$, so $v c(t)=0$ and
$c(t)v^*=0$.  These equalities kill respectively $v$ and $v^*$
in every product in \eqref{eq:entropy-exact-zeros}.

For every positive integer $s$ and every $\tau\in\T(A)$,
\begin{equation}\label{eq:entropy-high-moments}
 \tau(a_i^s)
   =\frac1m\tau(c(t)^{2s})
   \ge\frac1m\tau(1_{[a,1]}(t))\ge\beta_a.
\end{equation}
The trace equalities for the matrix slices follow first for
continuous functions vanishing at zero and then for spectral
projections by monotone convergence in the finite tracial GNS
closure.  The last inequality follows from
\[
 (1-a)\tau(1_{[0,a)}(t))
                  \le\tau(1-t)<\varepsilon.
\]

We now give a uniform covering estimate.  Fix
$R>\max\{\|h\|,\|k\|\}$ and
$0<\eta<\beta_a/4$, and set $\theta=\beta_a-2\eta$.
Take $b=1/2$, and choose an integer $s$ with $b^s<\eta$.
For a covering radius $0<\omega<1$, choose
\[
 0<\rho<\frac{\omega}{16m}
\]
small enough that the spectral bands above are disjoint, and
choose the corresponding $f_i$.  Since
$A=\mathrm C^*(h,k)$, there is a selfadjoint polynomial $p$ with
$p(h,k)$ arbitrarily close to $t$.  Write
\[
 \zeta(r)=\min\{1,\max\{0,r\}\}\quad(r\in\R),\qquad
 T=\zeta(p(H,K)),\qquad
 A_i=f_i(H)c(T)^2f_i(H).
\]
The letters $T,A_i$ in this display denote actual matrix functions
of a prospective microstate $(H,K)$, not additional microstate
coordinates.  They satisfy $0\le T,A_i\le1$.

Choose $p$ sufficiently accurately that evaluation at $(h,k)$
makes $A_i$ as close as needed to $a_i$ in
\eqref{eq:entropy-positive-cuts}.  Lemma
\ref{lem:entropy-trace-control}, applied also to the powers $A_i^s$,
and the norm part of that lemma, applied to the finitely many
products in \eqref{eq:entropy-exact-zeros}, then give finite
tests $\mathcal F,\sigma$ for which every
\[
       (H,K)\in
       \Gamma_R^{\mathrm{top}}((h,k);q,\mathcal F,\sigma)
\]
satisfies
\begin{equation}\label{eq:entropy-model-cuts}
 \operatorname{tr}_q(A_i^s)>\beta_a-\eta,
 \qquad
 \|A_i K A_j\|<\frac{b^2\omega}{16m^2}
                       \quad((i,j)\in\mathcal J).
\end{equation}
Here one first makes the errors at $(h,k)$ smaller than half the
desired tolerances and then imposes the finite microstate tests.
The polynomial $p$, the functions $f_i$, and these tests may
depend on $\omega$.  This is permitted by the infimum in
\eqref{eq:entropy-fixed-scale}.

Let
\[
 R_i=1_{[b,1]}(A_i),\qquad r_i=\rank R_i\quad(1\le i\le m),
 \qquad R_0=1-\sum_{i=1}^mR_i,\quad r_0=\rank R_0.
\]
If $P_i=1_{(\mu_i-\rho,\mu_i+\rho)}(H)$, then
$R_i\le P_i$.  Thus the $R_i$ are pairwise orthogonal, and
\begin{equation}\label{eq:entropy-rank-and-band}
 \frac{r_i}q
   \ge\operatorname{tr}_q(A_i^s)-b^s>\theta,
 \qquad
                  \|(H-\mu_i)R_i\|\le\rho.
\end{equation}
The rank estimate uses only $0\le A_i\le1$ and the spectral
cut at $b$.  The band estimate does not require $R_i$ to commute
with $H$.

The restriction of $A_i$ to $R_i\C^q$ is invertible with inverse
of norm at most $b^{-1}$.  Therefore
\[
 \|R_i K R_j\|\le b^{-2}\|A_i K A_j\|
       <\frac{\omega}{16m^2}
                         \qquad((i,j)\in\mathcal J).
\]
Define selfadjoint matrices
\[
 H'=\sum_{i=1}^m\mu_iR_i+R_0 H R_0,
 \qquad
 K'=K-\sum_{(i,j)\in\mathcal J}R_i K R_j.
\]
The set $\mathcal J$ is closed under transposition, so $K'$ is
selfadjoint.  The preceding estimates give
\begin{equation}\label{eq:entropy-model-error}
 \|H-H'\|\le2m\rho<\frac\omega8,
 \qquad
 \|K-K'\|<\frac\omega{16}.
\end{equation}
For the first estimate, write the difference as
\[
 \sum_i(H-\mu_i)R_i
                 +\sum_iR_i(H-\mu_i)R_0
\]
and use \eqref{eq:entropy-rank-and-band}.
Also $\|H'\|\le R$ and $\|K'\|\le R+1$.

For fixed ranks $r_0,r_1,\ldots,r_m$, the flags of orthogonal
projections with these ranks form the homogeneous space
\[
 \mathcal V_{\boldsymbol r}
       =\mathcal U(q)\big/
             \bigl(\mathcal U(r_0)\times\cdots\times
                                      \mathcal U(r_m)\bigr).
\]
Zero-rank factors are omitted.  Its real dimension is
\[
                 d_F=q^2-r_0^2-\sum_{i=1}^m r_i^2.
\]
By \cite[Corollary~12]{Szarek1998}, its covering numbers in the
quotient metric induced by the operator norm are bounded by
$(C/\alpha)^{d_F}$, with an absolute constant $C$.  The same
upper bound, after changing $C$, applies if one uses the intrinsic
operator-norm metric.  Conjugation of a matrix of norm at most
$R+1$ changes it by at most $2(R+1)$ times the distance of the
conjugating unitaries.  Consequently a flag net at radius a fixed
multiple of $\omega/(R+1)$ suffices.

For a fixed standard flag, the free residual block of $H'$ has
real dimension $r_0^2$.  The allowed selfadjoint matrices $K'$
form a real linear space of dimension
\[
 d_K=q^2-\sum_{i=N+1}^m r_i^2
       -2\!\sum_{\substack{1\le i<j\le m\\
                             \{i,j\}\notin\mathcal E}}r_i r_j.
\]
Their norm-bounded balls admit nets of cardinality at most
$(C_R/\omega)^{d_K}$ by the elementary volume estimate for a
ball in a finite-dimensional normed space.  The residual block
of $H'$ has a net of cardinality at most
$(C_R/\omega)^{r_0^2}$.

These flag and fibre nets combine without adding parameters for
the block-diagonal unitary group.  Indeed, that group preserves
the allowed linear space of $K'$ and the norm balls in both
fibres.  A block-diagonal adjustment in comparing a flag with a
flag-net point can therefore be absorbed into its fibre
parameters.  We obtain a net for all model pairs with the given
ranks of cardinality at most
\begin{equation}\label{eq:entropy-net-exponent}
 \left(\frac{C_R}{\omega}\right)^D,
 \qquad
 D=2q^2-\sum_{i=1}^m r_i^2-\sum_{i=N+1}^m r_i^2
       -2\!\sum_{\substack{1\le i<j\le m\\
                              \{i,j\}\notin\mathcal E}}r_i r_j.
\end{equation}
Here and below fixed numerical factors in the covering radius
are absorbed into $C_R$.  This constant is independent of $q$,
$\omega$, the widths of the $f_i$, and the polynomial $p$.
Allowing dependence on the fixed integer $m$ would give the same
conclusion.

The tree has $m-1$ edges.  The number of forbidden ordered
off-diagonal positions is therefore $(m-1)(m-2)$.  From
\eqref{eq:entropy-rank-and-band} and
\eqref{eq:entropy-net-exponent} we get
\[
 \frac D{q^2}
 \le2-\bigl(m+(m-N)+(m-1)(m-2)\bigr)\theta^2
       =2-(m^2-m+2-N)\theta^2.
\]
There are at most $(q+1)^{m+1}$ rank profiles.  Combining their
nets with \eqref{eq:entropy-model-error}, and choosing the fixed
numerical radii so that the resulting balls have radius $\omega$,
gives
\[
 N_\infty\bigl(
    \Gamma_R^{\mathrm{top}}((h,k);q,\mathcal F,\sigma),\omega\bigr)
 \le(q+1)^{m+1}
       \left(\frac{C_R}{\omega}\right)^%
       {\left(2-(m^2-m+2-N)\theta^2\right)q^2}.
\]
The polynomial rank-profile factor disappears after division by
$q^2$ and passage to the limit.  The constant $C_R$ disappears
after division by $|\log\omega|$ and $\omega\downarrow0$.
It follows that
\[
 \delta_{\mathrm{top}}(h,k)
       \le2-(m^2-m+2-N)(\beta_a-2\eta)^2.
\]
First let $\eta\downarrow0$, and then let $a\downarrow1/2$.
These are auxiliary parameters for the fixed pair $(h,k)$.
Equation~\eqref{eq:entropy-beta-a} gives precisely the upper
estimate in \eqref{eq:entropy-cone-upper}.

Since $N\le m-1$ and $0<\varepsilon<1/(64m)$, the upper bound
in \eqref{eq:entropy-cone-upper} is strictly less than two.
Equation~\eqref{eq:entropy-dimension-volume} therefore gives
$\chi_{\mathrm{top}}(h,k)=-\infty$.

The invertibility and the polar assertion follow from
Lemma~\ref{lem:polar-coordinates}, since $h\ge1$.  The lower bound
under the stated MF hypothesis is \eqref{eq:entropy-mf-lower}.
\end{proof}

\begin{theorem}\label{thm:entropy-infimum}
Let $A$ be a separable unital $\mathrm{MF}$ $\mathrm C^*$-algebra
with strict comparison by ordinary tracial states.  Suppose that
$\mathcal Z$ embeds unitally in $A$.  Then
\begin{equation}\label{eq:entropy-infimum}
 \inf\{\delta_{\mathrm{top}}(g):g\in A,\ \mathrm C^*(g)=A\}=1.
\end{equation}
More precisely, for every $\eta>0$ there is an invertible
generator $g=u b$, where $u\in\mathcal U(A)$ and $b=|g|\ge1_A$,
such that $\mathrm C^*(u,b)=\mathrm C^*(g)=A$ and
\[
 1\le\delta_{\mathrm{top}}(g)<1+\eta,\qquad
 \chi_{\mathrm{top}}(g)=-\infty.
\]
\end{theorem}

\begin{proof}
The MF assumption supplies a tracial state, and the unital copy
of $\mathcal Z$ excludes nonzero finite-dimensional
representations.  We may replace the choice $n=8$ in the proof of
Theorem~\ref{thm:z-comparison} by any integer $m\ge8$.
For completeness, fix the given unital copy of $\mathcal Z$ in
$A$.  Proposition~\ref{prop:z-cone} supplies a c.p.c.\ order-zero
map $\phi_m:\mathrm M_m\to\mathcal Z\subseteq A$ with
\[
 \tau(\phi_m(1_m))>1-\varepsilon_m\qquad(\tau\in\T(A)),
 \qquad 0<\varepsilon_m<\frac1{64m}.
\]
Set
\[
 d_m=(\phi_m(e_{11})-1/16)_+,
 \qquad D_m=\Her_A(d_m).
\]
The element $d_m$ is full in $A$ and has constant positive rank
over $\T(A)$.  Indeed, it is nonzero in the simple unital copy
of $\mathcal Z$, and every tracial state on $A$ restricts to
the unique trace of that copy.  The class $[d_m]$ is almost
divisible, as the image of a soft class in
$\mathrm{Cu}(\mathcal Z)$.
The argument in Theorem~\ref{thm:z-comparison} therefore applies:
Lemma
\ref{lem:audited-constant-rank-corner} gives a compact nonempty
normalized trace space, the structural and external AF results
apply, and Proposition~\ref{prop:audited-seven-generators}
gives seven selfadjoint generators for $D_m$.  This bound is
independent of $m$.

Since $m-1\ge7$, Proposition
\ref{prop:matrix-cone-compression} now gives an invertible
generator $g_m=h_m+ik_m$ with $N=7$.  Proposition
\ref{prop:entropy-cone} gives
\[
 1\le\delta_{\mathrm{top}}(g_m)
    \le2-\left(1-\frac1m-\frac5{m^2}\right)
                         (1-2\varepsilon_m)^2,
 \qquad \chi_{\mathrm{top}}(g_m)=-\infty.
\]
The upper bound tends to one as $m\to\infty$.  Its lower
bound is valid for every generator of $A$, by
\eqref{eq:entropy-mf-lower}.  This proves
\eqref{eq:entropy-infimum}.  Choosing $m$ sufficiently large
gives a cone pair with
$\delta_{\mathrm{top}}(h_m,k_m)<\min\{1+\eta,2\}$.
Choose $\alpha>0$ with $\alpha\|k_m\|<\pi$ and set
\[
             g=e^{i\alpha k_m}h_m,
             \qquad u=e^{i\alpha k_m},\quad b=h_m.
\]
Lemma~\ref{lem:polar-coordinates} gives $|g|=b$, phase $u$,
generation, and
$\delta_{\mathrm{top}}(g)=\delta_{\mathrm{top}}(h_m,k_m)$.
Equation
\eqref{eq:entropy-dimension-volume} then gives
$\chi_{\mathrm{top}}(g)=-\infty$.
\end{proof}

\begin{remark}\label{rem:entropy-eight}
Already the choice $m=8$, $N=7$ in the single-generation proof
gives an invertible generator with
\[
 1\le\delta_{\mathrm{top}}(g)
      \le2-\frac{51}{64}(1-2\varepsilon)^2<1.21,
 \qquad \chi_{\mathrm{top}}(g)=-\infty
\]
under the hypotheses of Theorem~\ref{thm:entropy-infimum}.
Increasing the matrix size gives the infimum in Theorem~\ref{thm:entropy-infimum};
the argument does not assert that a generator attains that
infimum.
\end{remark}

\end{document}